\documentclass[11pt,reqno]{amsart}
\usepackage[a4paper,left=1in,right=1in,top=1in,bottom=1in]{geometry}

\usepackage{amsmath,amssymb,amsthm,enumerate,array,cite}
\usepackage{graphicx}
\usepackage{tikz}
\usetikzlibrary{arrows.meta,decorations.pathmorphing,decorations.markings,
  decorations.pathreplacing,patterns,calc,positioning}
\usepackage{caption}
\usepackage{hyperref}
\hypersetup{colorlinks=true,linkcolor=blue,citecolor=blue,urlcolor=blue}

\numberwithin{equation}{section}

\theoremstyle{plain}
\newtheorem{theorem}{Theorem}[section]
\newtheorem{lemma}[theorem]{Lemma}
\newtheorem{proposition}[theorem]{Proposition}
\newtheorem{corollary}[theorem]{Corollary}
\theoremstyle{definition}
\newtheorem{definition}[theorem]{Definition}
\theoremstyle{remark}
\newtheorem{remark}[theorem]{Remark}
\newtheorem{example}[theorem]{Example}

\let\le\leqslant

\let\ge\geqslant

\newcommand{\R}{\mathbb R}
\newcommand{\C}{\mathbb C}
\newcommand{\Dm}{\mathbb D}
\newcommand{\Om}{\Omega}
\newcommand{\nv}{\nu}
\newcommand{\kap}{\kappa}
\newcommand{\dd}{\mathrm d}
\newcommand{\dive}{\operatorname{div}}
\newcommand{\curl}{\operatorname{curl}}
\newcommand{\Per}{\operatorname{Per}}
\newcommand{\Area}{\operatorname{Area}}
\newcommand{\Bv}{\bar B}
\newcommand{\pv}{\bar p}
\newcommand{\Pv}{\bar P}
\newcommand{\cof}{\operatorname{cof}}
\newcommand{\lam}{\lambda}
\newcommand{\doi}[1]{\href{https://doi.org/#1}{\texttt{doi:#1}}}

\title[Free boundary MHD equilibria]{Free boundary MHD equilibria with surface tension}
\author{Chengchun Hao}
\address{State Key Laboratory of Mathematical Sciences, Academy of Mathematics and Systems Science, Chinese Academy of Sciences, Beijing 100190, China \newline
	\and
	School of Mathematical Sciences, University of Chinese Academy of Sciences, Beijing 100049, China
}
\email{hcc@amss.ac.cn}
\date{}

\begin{document}
\begin{abstract}
Free boundary equilibria of an ideal incompressible magnetohydrodynamic (MHD)
plasma, at rest, with surface tension and an exterior vacuum are studied.  Within
the frozen-in orbit of a reference magnetic field, every equilibrium is a critical
point of the magnetic plus surface energy over the space of volume-preserving
embeddings; the capillary boundary condition is the first-variation boundary term,
and no assumption on the connectivity of the boundary is required.

In the plane the frozen-in law fixes the distribution of the stream function and
turns the energy into a rearrangement functional.  For a single-axis reference
field a sharp lower bound is proved that depends only on that distribution and is
attained by an explicit area-preserving diffeomorphism, smooth across the axis,
which yields a classical equilibrium.  For Morse stream functions sharper bounds
governed by the level-set topology are obtained, with a strict gap at an isolated
simple saddle merge.  In three dimensions no such reduction is available: the
energy does not control the full gradient of the embedding, no scalar
rearrangement principle applies, and the axisymmetric reduction yields only a
conditional lower bound.

\bigskip
\noindent\textit{Keywords:} MHD, capillary free boundary, stream function,
variational principle, volume-preserving maps, rearrangement

\smallskip
\noindent\textit{MSC 2020:} 76W05, 76B45, 35R35
\end{abstract}
\maketitle

\section{Introduction}

A stationary state of the
incompressible ideal magnetohydrodynamic (MHD) free boundary problem, with the
fluid at rest, is a triple $(\Om,\Bv,\pv)$ consisting of a bounded smooth
domain $\Om\subset\R^3$, a field $\Bv$, and a \emph{hydrodynamic} pressure
$\pv$, such that
\begin{equation}\label{e:eq}
\dive\Bv=0,\qquad
(\curl\Bv)\times\Bv=\nabla\pv\quad\text{in }\Om,\qquad
\Bv\cdot\nv=0\quad\text{on }\partial\Om,
\end{equation}
together with the free boundary condition
\begin{equation}\label{e:bc}
\pv+\tfrac12|\Bv|^2=\kap H\quad\text{on }\partial\Om,
\end{equation}
where $\nv$ is the outward unit normal, $\kap>0$ is a coefficient of surface tension and $H$ is the mean curvature of
$\partial\Om$ with respect to $\nv$, normalized so that a ball of radius $r$
has $H=2/r$.  The quantity $\Pv:=\pv+\tfrac12|\Bv|^2$ is the total pressure,
and \eqref{e:eq} is equivalent to $(\Bv\cdot\nabla)\Bv=\nabla\Pv$.  We call
$(\Om,\Bv,\pv)$ an \emph{equilibrium}.

The condition \eqref{e:bc} is the balance of the total pressure against
surface tension across the interface, the exterior being vacuum at zero
pressure.  Its combination with the interior equations is what makes the
problem free: \eqref{e:eq} determines $\Bv$ and $\pv$ from the domain and
from the data along the field lines, while \eqref{e:bc} constrains the
domain itself.

The exterior is
vacuum at zero pressure and carries no magnetic field.  There is thus a
vacuum interface, but not a plasma--vacuum interface of the kind treated in
\cite{BerestyckiBrezis1980,Temam1975,FujiiHirai1983}, where a \emph{nonzero}
exterior vacuum field is present, the interface is a flux surface of the
total field, and the boundary is determined by the solution rather than by a
prescribed curvature; the distinction is the vanishing of the exterior field,
not the absence of a vacuum region.  With $\nv$ pointing into the vacuum, the
jump condition $\nv\times(\Bv_{\rm vac}-\Bv)=K$ gives the surface current
$K=\Bv\times\nv$.  A nonzero tangential trace of $\Bv$ therefore jumps to
zero across the interface and is carried by this surface current.  The
interface condition imposed is the balance of normal stresses, \eqref{e:bc},
in which the total pressure $\Pv$ balances the capillary force $\kap H$; no
tangential stress condition and no flux-surface condition is imposed.  Fix
$\kap>0$; every statement below holds for each such $\kap$.  At $\kap=0$
the lower bound of Proposition~\ref{p:elow} vanishes and the surface term
drops out of the energy; the energy does not control the full gradient
of admissible maps, already on the disk (Example~\ref{ex:cav}).

A prescribed domain overdetermines the problem: \eqref{e:bc} then becomes an
extra constraint on $(\Bv,\pv)$ rather than a consequence of the interior
equations.  For instance, an axisymmetric field tangent to a sphere vanishes
at the two poles, so at those two points \eqref{e:bc} forces $\pv=2\kap/R$
while elsewhere it reads $\pv=2\kap/R-\tfrac12|\Bv|^2$.

We write $\R^n$ for Euclidean space, $\dd x$ for Lebesgue measure, and
$|E|$ for the volume of a set $E$.  For a bounded set $E$ of finite
perimeter, $\Per(E)$ denotes its perimeter, $\partial^*E\subseteq\partial E$
its reduced boundary, and $\Area(\partial E)$ the surface area of its
boundary; thus $\Per(E)=\mathcal H^{n-1}(\partial^*E)$, where $\mathcal H^k$
is the $k$-dimensional Hausdorff measure (so $\mathcal H^1$ is arclength in
the plane).  We use $C^\infty$ for smooth
maps, $\cof M$ for the cofactor matrix of a matrix $M$,
and $D\eta$ for the differential of a map $\eta$.  A volume-preserving
embedding $\eta$ acts on a field by push-forward, written $\eta_*\Bv_0$ and
defined in \eqref{e:Omega}.  In the plane we set $x^\perp=(-x_2,x_1)$ for
$x=(x_1,x_2)\in\R^2$, write $e_\theta$ for the angular unit vector,
$\Dm=\{x\in\R^2:|x|<1\}$ for the open unit disk, and let $\mathcal G_2$
denote the space of
area-preserving embeddings of a planar domain, defined in
Section~\ref{s:2d}.  In the plane $\curl\Bv$ denotes the scalar
$\partial_1\Bv_2-\partial_2\Bv_1$, and
$(\curl\Bv)\times\Bv:=(\curl\Bv)J\Bv$ with
$J=\begin{pmatrix}0&-1\\1&0\end{pmatrix}$.  We write
$\nabla^\perp\psi=(-\partial_2\psi,\partial_1\psi)$ in the plane and
$(e_r,e_\theta)$ for the polar frame.  For a differential form $\alpha$ and a
vector field $X$, $\iota_X\alpha$ is the interior product, $X^\flat$ the
metric-dual $1$-form, and $\dive_\omega$ the divergence with respect to the
form $\omega$.

The variational principle of Theorem~\ref{t:var} is the free boundary
version of the Arnold \cite{Arnold1966}--Ebin--Marsden
\cite{EbinMarsden1970} variational principle for the ideal (incompressible)
fluid equations, extended to ideal MHD; it asserts that a stationary ideal
MHD configuration is the critical point set of the magnetic energy on the
group of volume-preserving diffeomorphisms of a fixed domain.  The idea that magnetohydrostatic
equilibria are extremals of the magnetic energy subject to a constraint on
the field goes back to Woltjer
\cite{Woltjer1958}; the constraint appropriate to ideal MHD is the full
frozen-in law, which is what leads to the formulation on volume-preserving
embeddings rather than to a constraint on the helicity alone.  The
Hamiltonian formulation of ideal MHD \cite{MorrisonGreene1980} places
equilibria at critical points of the energy on Casimir leaves; the frozen-in
orbit is the configuration-space realization of such a leaf.  In the
axisymmetric case the equilibrium equations take the Grad--Shafranov form
\cite{GradRubin1958,Shafranov1966}; the scalar form of the planar equations
used in Section~\ref{ss:scalar} is the same reduction in two dimensions.

The free boundary problem obtained here is the capillary problem of a plasma
bounded by vacuum, and it differs from the classical plasma equilibrium
problems of Berestycki--Brezis \cite{BerestyckiBrezis1980}, Temam
\cite{Temam1975} and Fujii--Hirai \cite{FujiiHirai1983} in that the exterior
vacuum field there is nonzero and surface tension is absent, so that the free
boundary is a flux surface rather than a constant-curvature capillary
surface.  

For stationary
free boundary equilibria, the rigidity and flexibility of configurations with
nested flux surfaces were analyzed in \cite{Constantin2022}; in the planar
simply connected case they prove that nowhere-collinear $v$ and $B$ force the
domain to be a disk and both fields to be circular.  Our choice $v\equiv0$
lies outside that hypothesis: the disk conclusion below is a
\emph{variational minimization} statement --- the minimizer at fixed
stream-function distribution is a disk.  There the fluid may be in motion and the free boundary balances an
external vacuum field, whereas here the fluid is at rest and the boundary
condition is capillary.  The same authors had earlier established Liouville-type rigidity for suitable steady solutions without stagnation points in a two-dimensional periodic channel or in a hollowed-out axisymmetric cylinder, as well as the structural stability of Arnol\'d-stable solutions, \cite{ConstantinDrivasGinsberg2021}.
Three-dimensional stationary free-boundary equilibria supported by an
external surface current are constructed in \cite{EncisoLuquePeralta2025};
the interface condition and the weak, stepped-pressure solutions there
differ from the capillary equilibria considered here.  Steady axisymmetric free-boundary flows in a gravitating hydromagnetic star, with field-aligned and cross-field flows, are studied in \cite{LowMcIntosh2026}; that setting includes gravity and field-aligned motion and has no surface tension.

The \emph{evolution} problem for incompressible ideal MHD with a free surface
and surface tension has been studied in recent years: a priori estimates are obtained in
\cite{HaoLuo2014}, while local well-posedness is established in
\cite{GuLuoZhang2024,HaoYang2025}; the linearized problem and the
failure of well-posedness are analyzed in \cite{HaoLuo2021,HaoLuo2020}, and an
account of these developments is given in \cite{HaoLuo2022}.  For the plasma--vacuum interface
problem, local well-posedness with or without surface tension, and the
vanishing surface tension limit, are established in \cite{LiuXin2025}; the
corresponding results for current-vortex sheets are given in
\cite{LiuXinCVS2025}; moreover \cite{LuoYang2026} proved $H^3$ a
priori estimates in general bounded domains for the surface-tension problem, a
refined blow-up criterion separating topological self-intersection, loss of
boundary regularity, blow-up of the normal velocity and interior blow-up, and
the existence of smooth data whose solutions develop a finite-time
self-intersection of the free boundary.  These are evolution results; the
present paper determines static equilibria and the domains they occupy.

In the absence of surface tension the low-regularity
theory has recently been completed: sharp Hadamard well-posedness at the
regularity $s>\tfrac d2+1$, in arbitrary dimensions and on general (not
necessarily simply connected) domains, together with a matching continuation
criterion, is proved in \cite{IPTT2024}.  Hadamard well-posedness for the
same zero-surface-tension system, on ball-type domains for $5/2<s<4$,
together with nonuniform dependence of the solution map, is obtained in
\cite{HaoYang2026}.  Well-posedness there hinges on
the Taylor sign condition $-\partial_\nu\Pv>0$ on the boundary
\cite{Taylor1950}, a property of the base state: the results of
\cite{IPTT2024} show that it persists on a short time interval for data
satisfying it.  There $\Pv$ is the total pressure, the dynamic boundary
condition in the absence of surface tension is $\Pv=0$ on the free surface,
and the sign condition is imposed on $-\partial_\nu\Pv$.  

All equilibria
constructed below live on disks in the plane; in $\R^3$ the only equilibrium
exhibited is the ball with $\Bv=0$.  Their boundary curvature is constant,
so by \eqref{e:bc} and \eqref{e:bc2d} the
total pressure is constant on $\partial\Om$.  Subtracting that constant from
$\pv$ leaves \eqref{e:eq} unchanged and gives $\Pv=0$ on $\partial\Om$; with
the exterior vacuum at zero pressure, these equilibria are therefore
stationary solutions of the zero-surface-tension system as well.
Proposition~\ref{p:taylor} verifies the Taylor sign condition exactly for
the planar equilibria of Theorem~\ref{t:rearr}; the zero-field ball has
$\Pv$ constant and does not satisfy it.  For the capillary evolution problem
($\kap>0$) the interface
condition is \eqref{e:bc}, and stability of the equilibria is governed by
the second variation of $E$, not by the Taylor sign condition (Remark~\ref{r:hessian}).

Regularity of capillary free
boundaries goes back to the free boundary regularity theory of
Alt--Caffarelli \cite{AltCaffarelli1981}, which applies to minimizers of
area-type functionals with a volume constraint.  For planar Sobolev
homeomorphisms and their regularity we refer to \cite{HenclKoskela2014}.  The
rearrangement method used in Section~\ref{s:2d} is the classical variational
method for steady ideal flows and vortex patches \cite{Burton1989}, in which
the prescribed data are the rearrangement class of the vorticity; here the
rearrangement class is that of the stream function, and the domain is free.

Classical Schwarz
symmetrization acts on a fixed datum in a fixed domain; here the symmetrized
object is the stream function, and the symmetrization is realized by an
explicit area-preserving map of the reference domain, so that the minimizer is
a classical solution of \eqref{e:eq} and \eqref{e:bc2d} on a free domain
rather than an abstract rearranged function.  What is specific to the present
setting is that the level sets of the stream function fixed by the frozen-in
law, not those of the vorticity, carry the topology entering the multi-critical
bounds.
These existence statements are planar; in $\R^3$ the
outcome is the obstruction and the conditional bound described above.

The free boundary problem is the stationarity condition of the following energy:
\begin{equation}\label{e:E}
E(\eta)=\tfrac12\int_{\Om_0}|D\eta\,\Bv_0|^2\,\dd x
+\kap\,\Area\bigl(\eta(\partial\Om_0)\bigr),\qquad \eta\in\mathcal G,
\end{equation}
where $\mathcal G$ is the space of volume-preserving embeddings of $\Om_0$ into $\R^3$,
\begin{equation}\label{e:G}
\mathcal G=\bigl\{\eta\in C^\infty(\overline\Om_0,\R^3):
\eta\ \text{injective on}\ \overline\Om_0,\ \det D\eta=1\bigr\}.
\end{equation}

\medskip

In $\R^3$ the general result is an obstruction: the energy does
not control the full gradient of the embedding (Theorems~\ref{t:cavgen}
and~\ref{t:cav3d}), so the sublevel sets of $E$ are not weakly relatively
compact in $W^{1,2}$; the
axisymmetric reduction gives only a conditional lower bound
(Proposition~\ref{p:axisbound}); and, apart from the ball with $\Bv=0$, no
non-trivial three-dimensional equilibrium is constructed.  In the plane, by contrast, the frozen-in law reduces the problem to a
rearrangement at fixed distribution function, and the sharp bound is attained
by an explicit smooth equilibrium (Theorem~\ref{t:rearr}).

 (i) \emph{The reduction, and its scope.}  The free boundary
condition \eqref{e:bc} follows from stationarity with respect to $\eta$.

\begin{theorem}[Variational principle]\label{t:var}
Let $\eta$ be a volume-preserving embedding of $\Om_0$ and let
$\Bv=\eta_*\Bv_0$.  Then $\eta$ is a critical point of $E$ with respect to
volume-preserving variations if and only if $(\Om,\Bv,\pv)$ satisfies
\eqref{e:eq}--\eqref{e:bc} for some pressure $\pv$.  
\end{theorem}

The equivalence is relative to the reference field: it characterizes the
critical points inside the frozen-in orbit of $\Bv_0$, among fields of the
same field-line topology, not among all divergence-free fields.  The
quantitative form is \eqref{e:dE} of Theorem~\ref{t:var2}: for every
divergence-free $v$ the boundary term carries the normal stress
$\kap H-\tfrac12|\Bv|^2$, which the pressure absorbs into \eqref{e:bc}.

\medskip
(ii) \emph{Energy bounds.}  The
surface term bounds the shape from below, $E(\eta)\ge\kap(36\pi)^{1/3}V^{2/3}>0$
(Proposition~\ref{p:elow}), so minimizing sequences have bounded perimeter.
The magnetic term controls only the component of $D\eta$ along the frozen-in
direction, so the sublevel sets of $E$ are not weakly relatively compact in
$W^{1,2}$.

\begin{enumerate}[(a)]
	\item \emph{Relaxation.}  Over finite-perimeter domains and weak $L^2$
limits of the transported fields the relaxed problem has a minimizer under a
confinement hypothesis, a minimizing sequence lying, after translations, in a
fixed ball (Propositions~\ref{t:relax} and~\ref{t:relax2}); without it the
shape alone admits domains of fixed volume and bounded perimeter whose volume
partly escapes to infinity (Remark~\ref{r:escape}).  In the saddle
configuration of Section~\ref{ss:multi} the gap for the original problem is
proved positive (Propositions~\ref{p:notatt}, \ref{p:strictgap},
\ref{p:qgap}).

\item\emph{Gradient bounds.}  The gradient bound that would supply
compactness is not implied by the energy: along a closed field line one can
shear so that $\int_{\Om_0}|D\eta_k\Bv_0|^2$ stays fixed while
$\int_{\Om_0}|D\eta_k|^2\to\infty$ (Theorems~\ref{t:cavgen}
and~\ref{t:cav3d}).  In the plane this is sharpened by the algebraic identity
$|\cof D\eta|=|D\eta|$ (Lemma~\ref{l:2dmatrix}), so the obstruction is not a
failure of the estimates but of the energy.  Consequently the $W^{1,2}$
closure of $\mathcal G_2$ is strictly larger than $\mathcal G_2$, and weak
lower semicontinuity alone yields a limit in that closure with energy at
most $m_2$ (Lemma~\ref{t:attain2d}).  A minimizer over the enlarged class
$\widehat{\mathcal G}_2$ need not be $C^1$ (Theorem~\ref{t:twistdeg});
this is separate from the Eulerian relaxed class $\overline{\mathcal C}_2^{\mathrm{Lip}}$
of Definition~\ref{d:relaxed2}.
\end{enumerate}

\medskip
(iii)  \emph{The planar method: rearrangement.}  In the plane a
divergence-free field is a gradient perpendicular to its stream function,
$\Bv_0=\nabla^\perp\psi_0$, and the frozen-in law acquires a scalar
consequence: the distribution function $A_0(t)=|\{\psi_0<t\}|$ is invariant.
For $\psi=\psi_0\circ\eta^{-1}$ one then has $\Bv=\nabla^\perp\psi$ with
$\psi$ of the same distribution as $\psi_0$ (Lemma~\ref{l:conj}), so the
minimization is a rearrangement problem at fixed distribution, which the
P\'olya--Szeg\H{o} principle settles (Lemma~\ref{l:ps}).  This gives the sharp
lower bound and, in the single-axis case, its attainment in the smooth class
(Theorem~\ref{t:rearr}): the minimizer is the circular rearrangement
\eqref{e:psistar} of $\psi_0$, an explicit smooth area-preserving
diffeomorphism of $\Om_0$ onto the disk of area $V$ whose image solves
\eqref{e:eq} and \eqref{e:bc2d} classically.

\medskip
For circular fields on a disk the
rearrangement is trivial and the identity is a minimizer
(Theorems~\ref{t:disk} and~\ref{t:circ}); for Morse stream functions the level
topology gives strictly stronger bounds (Theorems~\ref{t:multi}
and~\ref{t:gen}), which at an isolated simple saddle merge are not attained,
the infimum being strictly larger (Propositions~\ref{p:notatt}
and~\ref{p:strictgap}).  In three dimensions the stream function disappears
and with it the method: the energy does not control the gradient, and the
bounds that survive are conditional (Section~\ref{s:3d}).

\section{The configuration space and the energy}\label{s:form}

\subsection{Reference configuration}

Fix a bounded domain $\Om_0\subset\R^3$ with $C^\infty$ boundary (of any
number of connected components), let $V=|\Om_0|$ be its volume, and let
$\nv_0$ be the outward unit normal of $\partial\Om_0$.  Fix also a field
$\Bv_0\in C^\infty(\overline{\Om_0};\R^3)$ with
\begin{equation}
\dive\Bv_0=0\quad\text{in }\Om_0,\qquad \Bv_0\cdot\nv_0=0
\quad\text{on }\partial\Om_0 .
\end{equation}
No connectivity of $\Om_0$ or of $\partial\Om_0$ is assumed anywhere: the
frozen-in law, the energy and the derivation of the free boundary condition
in Section~\ref{s:var} all apply componentwise, a point taken up in
Section~\ref{ss:equiv}.

\subsection{Embeddings and the frozen-in law}

Let $\mathcal G$ be the space of volume-preserving embeddings of $\Om_0$
defined in \eqref{e:G}; the injectivity is
part of the definition, since $\det D\eta=1$ alone gives only a local
diffeomorphism, whereas below we use that $\eta$ is a diffeomorphism of
$\overline\Om_0$ onto $\overline{\Om_\eta}$ (in particular in the change of
variables formula and in $|\Om_\eta|=V$).  This is a configuration
space; the volume-preserving diffeomorphisms of $\R^n$ act on it by left
composition, and the divergence-free variations of Section~\ref{ss:firstvar}
describe its infinitesimal action.  Rigid motions are the subgroup of
Remark~\ref{r:scale}.  For
$\eta\in\mathcal G$ write
\begin{equation}\label{e:Omega}
\Om_\eta=\eta(\Om_0),\qquad
\Bv_\eta(\eta(x))=D\eta(x)\,\Bv_0(x),
\end{equation}
so that $\Bv_\eta$ is the field carried along by the embedding, that is, the
push-forward of $\Bv_0$ as a vector field.  This is the frozen-in law of
ideal MHD.

\begin{lemma}\label{l:div}
For every $\eta\in\mathcal G$, the field $\Bv_\eta$ satisfies
$\dive\Bv_\eta=0$ in $\Om_\eta$ and $\Bv_\eta\cdot\nv=0$ on $\partial\Om_\eta$.
\end{lemma}

\begin{proof}
Let $B=\Bv_\eta$ and $a=(D\eta)^{-1}$, so that $\Bv_0(x)=a(x)\,B(\eta(x))$
by \eqref{e:Omega}; in components $B_0^j=a^j_i\,B^i(\eta)$.  Applying
$\partial_{x^j}$ and summing over $j$,
\[
(\partial_kB^i)(\eta)\,\partial_{x^j}\eta^k\,a^j_i
+B^i(\eta)\,\partial_{x^j}a^j_i
=\dive\Bv_0=0 .
\]
Since $\partial_{x^j}\eta^k\,a^j_i=\delta^k_i$ (as $a=(D\eta)^{-1}$), the
first term is $(\partial_kB^i)(\eta(x))\,\delta^k_i=(\dive B)(\eta(x))$.
The second term vanishes by the Piola identity: since $\det D\eta=1$ one has
$\cof D\eta=(D\eta)^{-T}=a^{T}$, so $(\cof D\eta)^i_j=a^j_i$ and
$\partial_{x^j}(\cof D\eta)^i_j=0$.  Hence $(\dive B)(\eta(x))=0$ for every
$x\in\Om_0$, and since $\eta$ is a diffeomorphism onto $\Om_\eta$ this says
$\dive B=0$ on $\Om_\eta$.

For the boundary, Nanson's formula gives
$\nv_{\eta}\,\dd S_\eta=(\cof D\eta)\,\nv_0\,\dd S_0=(D\eta)^{-T}\nv_0\,\dd S_0$
with the outer normals $\nv_0$ of $\partial\Om_0$ and $\nv_\eta$ of
$\partial\Om_\eta$: the oriented area element of the deformed boundary is
the cofactor matrix applied to that of the reference boundary, with matching
outward orientation because $\det D\eta=1>0$.  The normal pairing is
therefore proportional to
\[
\langle D\eta\Bv_0,(D\eta)^{-T}\nv_0\rangle
=\langle\Bv_0,\nv_0\rangle=0 ,\qedhere
\]
so $\Bv_\eta\cdot\nv_\eta=0$ on $\partial\Om_\eta$.
\end{proof}

\subsection{The energy}

\begin{definition}
The energy of the configuration $\eta\in\mathcal G$ is the functional
\eqref{e:E} of the introduction, written in Eulerian form as
\begin{equation}\label{e:E2}
E(\eta)=\tfrac12\int_{\Om_\eta}|\Bv_\eta|^2\,\dd y
+\kap\,\Area(\partial\Om_\eta).
\end{equation}
Equivalently, in terms of the reference configuration,
\begin{equation}
E(\eta)=\tfrac12\int_{\Om_0}|D\eta(x)\Bv_0(x)|^2\,\dd x
+\kap\int_{\partial\Om_0}\bigl|\cof(D\eta(x))\,\nv_0(x)\bigr|\,\dd S_0(x),
\end{equation}
the second form being the change of variables in \eqref{e:E2} together with
$\det D\eta=1$ and the area formula for the pushed-forward surface measure.
\end{definition}

The first term of \eqref{e:E2} is the magnetic energy and the second is the
surface energy.  The kinetic energy is absent because the state is at rest.

\begin{remark}\label{r:scale}
The energy \eqref{e:E2} is invariant under rigid motions of the ambient
space $\R^n$ ($n=2,3$) composed with $\eta$ on the left: writing
$(R\circ\eta+c)(x):=R\,\eta(x)+c$ for a rotation $R$ and $c\in\R^n$, one has
$E(R\circ\eta+c)=E(\eta)$.  Consequently critical points and minimizers come
in orbits of the Euclidean group.
\end{remark}

\subsection{The scalar form of a planar equilibrium}\label{ss:scalar}

In the plane the equilibrium equations collapse to a single scalar equation,
the form used in the minimization of Section~\ref{s:2d}.

\begin{proposition}[Scalar form]\label{p:scalar}
Let $\Om\subset\R^2$ be a bounded simply connected domain and let
$\Bv\in C^1(\Om;\R^2)$ be divergence free.  Suppose
$(\Om,\Bv,\pv)$ satisfies \eqref{e:eq} with $\pv\in C^1(\Om)$.  Then
$\Bv=\nabla^\perp\psi$ for some $\psi\in C^2(\Om)$.  For every
$x\in\{\nabla\psi\ne0\}$ there is a submersion chart $W_x$ and a $C^1$
function $F_x$ on the interval $\psi(W_x)$ such that
\begin{equation}\label{e:scalar}
\Delta\psi=F_x'(\psi),\qquad \pv=-F_x(\psi)\quad\text{on }W_x .
\end{equation}
If $\psi$ is continuous up to $\partial\Om$, then $\psi=c$ on $\partial\Om$
provided $\partial\Om$ is connected, where $c$ is the constant boundary
value of $\psi$, and the capillary condition
$\pv+\tfrac12|\Bv|^2=\kap\kappa_{\mathrm{curv}}$ reads
\begin{equation}\label{e:scalarbc}
\psi=c\ \ \text{on }\partial\Om,\qquad
\tfrac12|\nabla\psi|^2+\pv=\kap\kappa_{\mathrm{curv}}\ \ \text{on }\partial\Om ,
\end{equation}
with $\kappa_{\mathrm{curv}}$ the scalar curvature of $\partial\Om$ with
respect to the outward normal.
\end{proposition}

\begin{proof}
Simple connectivity and $\dive\Bv=0$ give $\Bv=\nabla^\perp\psi$ with
$\psi\in C^2$, and $\curl\Bv=\Delta\psi$ because
$\curl\nabla^\perp\psi=\Delta\psi$ in the plane.  By
$J\nabla^\perp\psi=-\nabla\psi$ and $\det J=1$,
\[
(\curl\Bv)\times\Bv=\Delta\psi\,J\nabla^\perp\psi=-\Delta\psi\,\nabla\psi ,
\]
so \eqref{e:eq} reads $\nabla\pv=-\Delta\psi\,\nabla\psi$; hence $\nabla\pv$ is
parallel to $\nabla\psi$.  Fix $x$ with $\nabla\psi(x)\ne0$.  By the
submersion theorem, there is a neighborhood $W_x$ on which $\psi$ is one
coordinate and its coordinate slices are connected.  Since $\pv$ has zero
derivative along each such slice, it is a function of $\psi$ on $W_x$;
writing $\pv=-F_x(\psi)$ and differentiating gives
$\Delta\psi=F_x'(\psi)$, which proves \eqref{e:scalar}.

Finally, if $\psi$ is
continuous up to $\partial\Om$, then $\Bv\cdot\nv=0$ makes $\nabla\psi$
normal to $\partial\Om$, so $\psi$ is constant on each component of
$\partial\Om$; with $\partial\Om$ connected this gives $\psi=c$.  Where the
capillary condition $\pv+\tfrac12|\Bv|^2=\kap\kappa_{\mathrm{curv}}$ is
imposed, $|\Bv|=|\nabla\psi|$ turns it into \eqref{e:scalarbc}.
\end{proof}

\section{The variational principle}\label{s:var}

\subsection{The first variation}\label{ss:firstvar}

Let $\eta\in\mathcal G$ and let $v\in C^\infty(\overline\Om;\R^n)$ be divergence
free, where $\Om=\Om_\eta$.  Extend $v$ across a collar of $\partial\Om$ by
extending the closed $(n-1)$-form
$\iota_v(\dd x_1\wedge\cdots\wedge\dd x_n)$ to a closed form on the collar
and dualizing; this keeps $\dive v=0$ on the larger set.  Work in collar
coordinates
$(z,\rho)\in\partial\Om\times(-\varepsilon,\varepsilon)$, with $\rho<0$
inside $\Om$, and fix a nowhere-vanishing $(n-1)$-form $\omega$ on
$\partial\Om$, pulled back to the collar independently of $\rho$.  Every
$(n-1)$-form on the collar decomposes
as $\beta=f\omega+\dd\rho\wedge\gamma$, with $f$ a function and $\gamma$ an
$(n-2)$-form on $\partial\Om$ depending on $\rho$; since
$\dd_{\partial\Om}\omega=0$,
\[
\dd\beta=\dd\rho\wedge\bigl(\partial_\rho f\,\omega
-\dd_{\partial\Om}\gamma\bigr),
\]
so $\beta$ is closed exactly when
$\partial_\rho f\,\omega=\dd_{\partial\Om}\gamma$.  The interior field fixes
the restrictions of $f$ and $\gamma$ to $\rho\le0$ together with all their
$\rho$-derivatives at $\rho=0$; these satisfy
$\partial_\rho^{\,k+1}f\,\omega
=\dd_{\partial\Om}\partial_\rho^{\,k}\gamma$ for every $k\ge0$.  Choose a
smooth extension of $\gamma$ to $\rho>0$ whose jets at $\rho=0$ agree with
those of the interior data (Borel's theorem) and solve
$\partial_\rho f\,\omega=\dd_{\partial\Om}\gamma$ on the whole collar with
the interior trace of $f$.  The jets of $f$ then agree with the interior ones
as well, and $\beta$ is a smooth closed $(n-1)$-form on the collar extending
the interior form.  Dualizing by
$\iota_v(\dd x_1\wedge\cdots\wedge\dd x_n)=\beta$ extends $v$ as a
divergence-free field, because
$\dd\bigl(\iota_v(\dd x_1\wedge\cdots\wedge\dd x_n)\bigr)
=(\dive v)\,\dd x_1\wedge\cdots\wedge\dd x_n$; the extension is defined on
the whole neighborhood $\Om\cup(\partial\Om\times[0,\varepsilon))$ of
$\overline\Om$.  

Denote by $\phi_\varepsilon$ the flow of the extended field,
defined
for $|\varepsilon|\le\varepsilon_0$ on a neighborhood of $\overline\Om$;
since $\dive v=0$ we have $\det D\phi_\varepsilon=1$, so
$\eta_\varepsilon:=\phi_\varepsilon\circ\eta\in\mathcal G$ and
$\eta_\varepsilon$ is a volume-preserving variation of $\eta$ whose velocity
at $\varepsilon=0$ is $v\circ\eta$.  Conversely, if $\eta_\varepsilon$ is a
smooth curve in $\mathcal G$, then
$v:=(\dd\eta_\varepsilon/\dd\varepsilon)|_{\varepsilon=0}\circ\eta^{-1}$ is
a smooth field on $\overline\Om$ with $\dive v=0$, because
$0=\frac{\dd}{\dd\varepsilon}\det D\eta_\varepsilon
=\det D\eta_\varepsilon\,\dive v\circ\eta_\varepsilon$.  Hence
\begin{equation}\label{e:tangent}
T_\eta\mathcal G=\bigl\{v\in C^\infty(\overline\Om;\R^n):\ \dive v=0\bigr\}.
\end{equation}
For $v\in T_\eta\mathcal G$ write
\begin{equation}
\delta E(\eta)[v]:=\frac{\dd}{\dd\varepsilon}\Big|_{\varepsilon=0}
E(\eta_\varepsilon).
\end{equation}

\begin{theorem}[The first variation of the energy]\label{t:var2}
Let $\eta\in\mathcal G$, $\Om=\Om_\eta$, $\Bv=\Bv_\eta$.  Then
\begin{equation}\label{e:dE}
\delta E(\eta)[v]
=\int_\Om v\cdot\bigl(-(\curl\Bv)\times\Bv\bigr)\,\dd y
+\oint_{\partial\Om}\bigl(\kap H-\tfrac12|\Bv|^2\bigr)\,v\cdot\nv\,\dd S
\end{equation}
for every $v\in T_\eta\mathcal G$.
\end{theorem}

\begin{proof}
We compute the two terms of \eqref{e:E2} separately.

\medskip\noindent\textit{The magnetic term.}
Let $B_\varepsilon=\Bv_{\eta_\varepsilon}$.  From \eqref{e:Omega},
$B_\varepsilon(\eta_\varepsilon(x))=D\eta_\varepsilon(x)\Bv_0(x)$.  Since
$\dd\eta_\varepsilon/\dd\varepsilon|_{\varepsilon=0}=v\circ\eta$ we have
$\dd(D\eta_\varepsilon)/\dd\varepsilon|_{\varepsilon=0}=D(v\circ\eta)=(Dv)(\eta)\,D\eta$,
so
\[
\frac{\dd}{\dd\varepsilon}\Big|_{\varepsilon=0}
\left(\frac12\int_{\Om_0}|D\eta_\varepsilon\Bv_0|^2\dd x\right)
=\int_{\Om_0}\Bv(\eta)\cdot(Dv)(\eta)\Bv(\eta)\,\dd x .
\]
Changing variables $y=\eta(x)$, which is volume preserving, gives
\begin{equation}
\frac{\dd}{\dd\varepsilon}\Big|_{\varepsilon=0}
\left(\frac12\int_{\Om_\varepsilon}|B_\varepsilon|^2\dd y\right)
=\int_{\Om}\Bv\cdot\bigl((\Bv\cdot\nabla)v\bigr)\,\dd y .
\end{equation}

We now rewrite the right-hand side.  Integrating by parts and using
$\dive\Bv=0$ (Lemma~\ref{l:div}),
\[
\int_\Om\Bv\cdot((\Bv\cdot\nabla)v)\dd y
=\int_\Om\Bv_i\Bv_j\,\partial_jv_i\,\dd y
=-\int_\Om v\cdot((\Bv\cdot\nabla)\Bv)\,\dd y
+\oint_{\partial\Om}(\Bv\cdot v)(\Bv\cdot\nv)\dd S,
\]
and the boundary term vanishes by $\Bv\cdot\nv=0$.  By the identity
$(\Bv\cdot\nabla)\Bv=(\curl\Bv)\times\Bv+\nabla(\tfrac12|\Bv|^2)$,
\[
-\int_\Om v\cdot((\Bv\cdot\nabla)\Bv)\dd y
=-\int_\Om v\cdot\bigl((\curl\Bv)\times\Bv\bigr)\dd y
-\int_\Om v\cdot\nabla\bigl(\tfrac12|\Bv|^2\bigr)\dd y .
\]
The last integral is
$-\oint_{\partial\Om}\tfrac12|\Bv|^2\,v\cdot\nv\,\dd S
+\int_\Om\tfrac12|\Bv|^2\,\dive v\,\dd y$, and $\dive v=0$.  Combining,
\begin{equation}\label{e:identity}
\int_{\Om}\Bv\cdot\bigl((\Bv\cdot\nabla)v\bigr)\dd y
=-\int_\Om v\cdot\bigl((\curl\Bv)\times\Bv\bigr)\dd y
-\tfrac12\oint_{\partial\Om}|\Bv|^2\,v\cdot\nv\,\dd S .
\end{equation}

\medskip\noindent\textit{The surface term.}
The first variation of area under a boundary displacement $v\cdot\nv$ is
$\dd\Area=\oint_{\partial\Om}H\,v\cdot\nv\,\dd S$, so
\begin{equation}\label{e:dEsurf}
\frac{\dd}{\dd\varepsilon}\Big|_{\varepsilon=0}\kap\,\Area(\partial\Om_\varepsilon)
=\kap\oint_{\partial\Om}H\,v\cdot\nv\,\dd S .
\end{equation}

Adding \eqref{e:identity} and \eqref{e:dEsurf} gives \eqref{e:dE}.
\end{proof}

\subsection{The equivalence}\label{ss:equiv}

The passage from \eqref{e:dE} to the free boundary condition is a statement
about normal traces of solenoidal fields, and it holds with any number of
boundary components.

\begin{lemma}[Normal traces]\label{l:trace}
Let $\Om$ be a bounded domain with $C^\infty$ boundary and connected
components $\Om^1,\dots,\Om^N$, each with any number of boundary components.
Then the normal traces $v\cdot\nv$ of fields $v\in T_\eta\mathcal G$ are
exactly the smooth functions $g$ on $\partial\Om$ satisfying
\[
\oint_{\partial\Om^j}g\,\dd S=0,\qquad j=1,\dots,N .
\]
\end{lemma}

\begin{proof}
Necessity is the divergence theorem applied on each $\Om^j$.  For
sufficiency, solve $\Delta\varphi=0$ in each $\Om^j$ with
$\partial_\nu\varphi=g$ on $\partial\Om^j$; the Neumann problem is solvable
on $\Om^j$ because of the corresponding compatibility condition, and
$v=\nabla\varphi$ is divergence free with $v\cdot\nv=g$.  The extension of
$v$ across a collar used to generate volume-preserving variations is the
divergence-free extension of Section~\ref{ss:firstvar}.
\end{proof}

\begin{proof}[Proof of Theorem~\ref{t:var}]
Let $\eta\in\mathcal G$, $\Om=\Om_\eta$ and $\Bv=\Bv_\eta$.  By
\eqref{e:tangent}, $\eta$ is a critical point of $E$ if and only if
$\delta E(\eta)[v]=0$ for every $v\in T_\eta\mathcal G$.

Suppose first that this holds and write $w=(\curl\Bv)\times\Bv$.  Taking $v$
compactly supported in $\Om$, the boundary term is absent and
$\int_\Om v\cdot w\,\dd y=0$ for all such $v$.  Taking
$v=\curl(\varphi e_i)$ with $\varphi\in C_c^\infty(\Om)$ in $\R^3$, and
$v=\nabla^\perp\varphi$ in the plane, shows that $w$ is a closed
$1$-form, that is, $\curl w=0$ in $\Om$.

Solve $\Delta\phi=\dive w$ in $\Om$ with $\partial_\nu\phi=w\cdot\nv$ on
$\partial\Om$; the compatibility condition on each connected component of
$\Om$ is the divergence theorem, and elliptic regularity gives
$\phi\in C^\infty(\overline\Om)$.  Put $h=w-\nabla\phi$, so that
$\dive h=\curl h=0$ and $h\cdot\nv=0$ on $\partial\Om$; by \eqref{e:tangent},
$h\in T_\eta\mathcal G$.  Since $\dive h=0$ and $h\cdot\nv=0$,
$\int_\Om h\cdot\nabla\phi\,\dd y=0$, and the boundary term of \eqref{e:dE}
vanishes for $v=h$; therefore
\[
0=\delta E(\eta)[h]=-\int_\Om h\cdot w\,\dd y=-\int_\Om|h|^2\,\dd y .
\]
Thus $h=0$ and $w=\nabla\phi=:\nabla\pv$ with $\pv\in C^\infty(\Om)$, which
is \eqref{e:eq}.  Substituting this into
\eqref{e:dE} and using $\dive v=0$,
\[
-\int_\Om v\cdot\nabla\pv\,\dd y
=\int_\Om\pv\,\dive v\,\dd y-\oint_{\partial\Om}\pv\,v\cdot\nv\,\dd S
=-\oint_{\partial\Om}\pv\,v\cdot\nv\,\dd S ,
\]
so that
\[
0=\delta E(\eta)[v]
=\oint_{\partial\Om}\bigl(\kap H-\tfrac12|\Bv|^2-\pv\bigr)v\cdot\nv\,\dd S
\qquad\text{for every }v\in T_\eta\mathcal G .
\]
By Lemma~\ref{l:trace} the factor $\kap H-\tfrac12|\Bv|^2-\pv$ is
$L^2(\partial\Om^j)$-orthogonal to every zero-mean function on the boundary
of each connected component $\Om^j$, hence equals a constant $C_j$ there.

The equation $\nabla\pv=(\curl\Bv)\times\Bv$ determines $\pv$ up to one
additive constant on each $\Om^j$, so replacing $\pv$ by $\pv+C_j$ on
$\Om^j$ leaves $\nabla\pv$ unchanged and gives
$\pv+\tfrac12|\Bv|^2=\kap H$ on $\partial\Om$.  The converse is immediate
from \eqref{e:dE}.  
\end{proof}

\begin{example}[Recovering the ball]
If $\Bv_0=0$ then $E(\eta)=\kap\Area(\partial\Om_\eta)$, and the isoperimetric
inequality shows that the infimum is at least
$\kap\,(36\pi)^{1/3}V^{2/3}$.  It is attained by the embeddings onto balls of
volume $V$, which exist whenever $\Om_0$ is diffeomorphic to a ball (for
instance when $\Om_0$ is itself a ball); the corresponding equilibria are
balls with $\pv=\kap H$ constant.  This is the classical capillary drop.
\end{example}

\begin{example}[Two dimensions]\label{ex:2d}
The construction is not special to three dimensions.  In the plane take
$\Om_0$ to be the unit disk and $\Bv_0=\lam x^\perp=\lam(-x_2,x_1)$, which is
divergence free and tangent to the circle.  The disk with the same field and
with $\pv=\kap+\tfrac12\lam^2-\lam^2|x|^2$ satisfies
$(\curl\Bv)\times\Bv=\nabla\pv$ in $\Om_0$ and
$\pv+\tfrac12|\Bv|^2=\kap$ on $|x|=1$, that is, the two-dimensional boundary
condition \eqref{e:bc2d} of Section~\ref{s:2d}, with $\kappa_{\mathrm{curv}}=1$.
By Theorem~\ref{t:var2d} there the identity embedding is therefore a critical
point of $E$.
\end{example}

\section{The minimization problem}\label{s:min}

\subsection{A lower bound for the energy}

\begin{proposition}\label{p:elow}
The energy $E$ is bounded below on $\mathcal G$:
\begin{equation}\label{e:elow}
E(\eta)\ge\kap\,\Area(\partial\Om_\eta)
\ge\kap\,(36\pi)^{1/3}V^{2/3}>0 .
\end{equation}
\end{proposition}

\begin{proof}
The first inequality is immediate from \eqref{e:E2}.  For the second, note
that $\eta$ is volume preserving, so $|\Om_\eta|=V$, and apply the
isoperimetric inequality:
$\Area(\partial\Om_\eta)\ge(36\pi)^{1/3}V^{2/3}$.
\end{proof}

Thus $m:=\inf_{\mathcal G}E\ge\kap(36\pi)^{1/3}V^{2/3}$.
The shape is controlled conditionally: \eqref{e:elow} bounds the perimeter
along minimizing sequences, and bounded perimeter leaves the shape compact
modulo translation (Lemma~\ref{l:compact}), but without the fixed-ball
hypothesis of that lemma volume can still escape to infinity
(Remark~\ref{r:escape}).

\subsection{The relaxed problem}

Let $\mathcal G$ be as in \eqref{e:G} and let
$\eta_n\in\mathcal G$ be a minimizing sequence, $E(\eta_n)\to m$.  Write
$\Om_n=\Om_{\eta_n}$ and $\Bv_n=\Bv_{\eta_n}$, extended by zero outside
$\Om_n$.

\begin{lemma}\label{l:compact}
Work in $\R^d$ with $d=2$ or $d=3$.  Let $\Om_n$ be sets with
$|\Om_n|=V$ and $\Per(\Om_n)\le C$, and let
$\Bv_n\in L^2(\R^d;\R^d)$ satisfy $\int_{\R^d}|\Bv_n|^2\le C$,
$\dive\Bv_n=0$ and $\Bv_n=0$ almost everywhere on $\R^d\setminus\Om_n$.
Assume in addition
that, after a translation, the sets are contained in a fixed ball,
$\Om_n\subset B(0,R_0)$ for all $n$.  Then there is a subsequence, again
indexed by $n$, a set $\Om_\infty$ of finite perimeter with
$|\Om_\infty|=V$, and a field
$\Bv_\infty\in L^2(\R^d;\R^d)$, such that
$\Om_n\to\Om_\infty$ in $L^1$, $\Bv_n\rightharpoonup\Bv_\infty$ weakly in
$L^2$, and
\[
\Per(\Om_\infty)\le\liminf_n\Per(\Om_n),\qquad
\int_{\R^d}|\Bv_\infty|^2\le\liminf_n\int_{\R^d}|\Bv_n|^2 .
\]
Also $\dive\Bv_\infty=0$ in the sense of distributions, and
$\Bv_\infty=0$ almost everywhere on $\R^d\setminus\Om_\infty$.
\end{lemma}

\begin{proof}
The perimeters are bounded by hypothesis, $\Per(\Om_n)\le C$, and
$|\Om_n|=V$.  Since the $\Om_n$ lie in
the fixed ball $B(0,R_0)$, the compactness theorem for sets of finite
perimeter in a bounded domain gives a subsequence converging in $L^1_{\rm
loc}(\R^d)$, hence in $L^1(\R^d)$ because all the sets are supported in
$B(0,R_0)$, to a set $\Om_\infty$ with $|\Om_\infty|=V$, and the perimeter is
lower semicontinuous, so $\Per(\Om_\infty)\le\liminf\Per(\Om_n)$.

For the fields $\int_{\R^d}|\Bv_n|^2\le C$ by hypothesis, so passing to a
further subsequence of the one just chosen,
$\Bv_n\rightharpoonup\Bv_\infty$ weakly in $L^2$, and the norm is lower
semicontinuous.  Divergence-free is a closed linear condition, preserved
under weak convergence.  For the support condition pair $\Bv_n$ with $\phi\,\chi_{\R^d\setminus\Om_\infty}$ for
$\phi\in C_c^\infty(\R^d)$.  Since $\Bv_n=\Bv_n\chi_{\Om_n}$ and
$\chi_{\Om_n\cap(\R^d\setminus\Om_\infty)}\le|\chi_{\Om_n}-\chi_{\Om_\infty}|$ pointwise,
\[
\Bigl|\int_{\R^d\setminus\Om_\infty}\Bv_n\cdot\phi\Bigr|
=\Bigl|\int_{\R^d}\Bv_n\cdot\phi\,\chi_{\Om_n\cap(\R^d\setminus\Om_\infty)}\Bigr|
\le\|\phi\|_\infty\,\|\Bv_n\|_{L^2}\,
|\Om_n\triangle\Om_\infty|^{1/2}\longrightarrow0 ,
\]
because $|\Om_n\triangle\Om_\infty|\to0$ and $\|\Bv_n\|_{L^2}\le\sqrt C$.  As
$\Bv_n\rightharpoonup\Bv_\infty$ weakly in $L^2$ and
$\phi\,\chi_{\R^d\setminus\Om_\infty}\in L^2$, passing to the limit gives
$\int_{\R^d\setminus\Om_\infty}\Bv_\infty\cdot\phi=0$ for every
$\phi\in C_c^\infty(\R^d)$, hence $\Bv_\infty=0$ almost everywhere on
$\R^d\setminus\Om_\infty$.
\end{proof}

\begin{definition}\label{d:relaxed}
Let $X$ be the space of pairs $(\Om,B)$ where $\Om\subset\R^3$ is a set of
finite perimeter and volume $V$, taken as an almost-everywhere equivalence
class: only $\chi_\Om$, $|\Om|$ and $\Per(\Om)$ enter, so the choice of
representative is immaterial and no topological closure of $\Om$ is
involved.  Let
$B\in L^2(\R^3;\R^3)$ satisfy $\dive B=0$ in distributions
and $B=0$ almost everywhere on $\R^3\setminus\Om$.  Give $X$ the product
topology of $L^1$ convergence of characteristic functions and weak $L^2$
convergence of fields.  Along bounded sequences the support condition is closed
in that topology: if $\chi_{\Om_n}\to\chi_\Om$ in $L^1$,
$B_n\rightharpoonup B$ in $L^2$ with $\sup_n\|B_n\|_{L^2}<\infty$ and
$B_n=0$ almost everywhere off $\Om_n$, then for every
$\phi\in C_c^\infty(\R^3)$
\[
\Bigl|\int_{\R^3\setminus\Om}B_n\cdot\phi\Bigr|
\le\|\phi\|_\infty\,\|B_n\|_{L^2}\,|\Om_n\triangle\Om|^{1/2}\longrightarrow0 ,
\]
so, passing to the weak limit, $\int_{\R^3\setminus\Om}B\cdot\phi=0$, whence
$B=0$ almost everywhere off $\Om$.  Let
\[
\mathcal C_0=\{(\Om_\eta,\Bv_\eta):\eta\in\mathcal G\},\qquad
\overline{\mathcal C}=\operatorname{cl}_X(\mathcal C_0).
\]
Thus $\overline{\mathcal C}$ is the relative topological closure of the
fixed-reference frozen-in orbit, not the class of all divergence-free
fields.  The relaxed energy is
\[
E(\Om,B)=\frac12\int_{\R^3}|B|^2\,\dd y+\kap\Per(\Om).
\]
\end{definition}

\begin{proposition}[Compactness for the relaxed problem]\label{t:relax}
Let $m'=\inf_{\overline{\mathcal C}}E$, and assume that the relaxed problem
admits a minimizing sequence which, after translations, is contained in a
fixed ball.  Then the relaxed energy attains $m'$ on
$\overline{\mathcal C}$, and
\[
m'\ \le\ \inf_{\mathcal G}E=m .
\]
\end{proposition}

\begin{proof}
Let $(\Om_n,\Bv_n)$ be a minimizing sequence of the relaxed problem with
$E(\Om_n,\Bv_n)\to m'$.  Post-composing the underlying embedding with a
translation shows that $\mathcal C_0$, hence $\overline{\mathcal C}$, is
invariant under translations of the ambient space, and so is $E$
(Remark~\ref{r:scale}); the
hypothesis therefore provides a minimizing sequence with the same energies
whose domains all lie in a fixed ball.  Relabeling it $(\Om_n,\Bv_n)$, and
noting that both terms of $E(\Om_n,\Bv_n)$ are nonnegative and
$E(\Om_n,\Bv_n)\to m'$, one has $\Per(\Om_n)\le(m'+1)/\kap$ and
$\int_{\R^3}|\Bv_n|^2\le2(m'+1)$ for all large $n$; the remaining hypotheses
of Lemma~\ref{l:compact} hold by Definition~\ref{d:relaxed} and the
fixed-ball assumption.  Hence
Lemma~\ref{l:compact} gives a subsequence converging to some
$(\Om_\infty,\Bv_\infty)\in X$ in the sense of Definition~\ref{d:relaxed},
with
\[
E(\Om_\infty,\Bv_\infty)
\le\liminf_n\Bigl(\frac12\int|\Bv_n|^2+\kap\Per(\Om_n)\Bigr)=m'.
\]
Since the sequence lies in $\overline{\mathcal C}$, which is closed in $X$ by
definition, its $X$-limit $(\Om_\infty,\Bv_\infty)$ lies in
$\overline{\mathcal C}$.  Hence
$E(\Om_\infty,\Bv_\infty)\ge m'$ as well, so the minimum is attained and
equals $m'$.  Finally, since $\mathcal C_0\subseteq\overline{\mathcal C}$ we
have $m'\le m$.
\end{proof}

\begin{remark}[Escape to infinity]\label{r:escape}
The boundedness hypothesis in Lemma~\ref{l:compact} and
Proposition~\ref{t:relax} is not automatic for the shape alone.  A domain of
fixed volume and bounded perimeter can carry a fixed fraction of its volume
to infinity, by growing a long thin neck that ends in a lobe of fixed
volume: a neck of length $\ell$ and radius $\varepsilon=C/\ell$ has
perimeter $2\pi\varepsilon\ell=2\pi C$ and
volume $\pi\varepsilon^2\ell=\pi C^2/\ell\to0$, so the neck costs a
bounded amount of perimeter and a vanishing amount of volume, while the lobe
carries a fixed fraction of $V$ to distance of order $\ell$.  Hence a sequence
of domains of fixed volume and bounded perimeter need not converge in
$L^1(\R^3)$.  Whether the same escape can be realized by a frozen-in orbit
is decided by the magnetic term: a neck of radius
$\varepsilon$ carrying a flux $\Phi\ne0$ has magnetic energy at least
$\tfrac12\ell\Phi^2/(\pi\varepsilon^2)=\Phi^2\ell^3/(2\pi C^2)$, by the
Cauchy--Schwarz inequality on each cross-section and integration along the neck, so
the construction is available only when the escaping lobe carries away a
vanishing amount of flux.

In the two-dimensional existence results of
Section~\ref{s:2d}
(Theorems~\ref{t:disk}, \ref{t:circ} and~\ref{t:rearr}) the issue does not
arise: they are obtained by direct minimization and do not use the relaxed
problem.
\end{remark}

\section{The two-dimensional problem}\label{s:2d}

\subsection{Formulation, the boundary term, and the matrix identity}

Let $\Om_0\subset\R^2$ be a bounded domain with $C^\infty$ boundary, let
$|\Om_0|=V$, and let $\Bv_0\in C^\infty(\overline{\Om_0};\R^2)$ satisfy
$\dive\Bv_0=0$ and $\Bv_0\cdot\nv_0=0$ on $\partial\Om_0$.  The space
$\mathcal G_2$ of area-preserving embeddings, the push-forward
$\Bv_\eta(\eta(x))=D\eta(x)\Bv_0(x)$ and the energy
\begin{equation}\label{e:E2d}
E(\eta)=\frac12\int_{\Om_0}|D\eta(x)\Bv_0(x)|^2\,\dd x
+\kap\,\Per(\eta(\Om_0))
\end{equation}
are defined as in Sections~\ref{s:form}--\ref{s:min}, with $\Per$ now the
length of the boundary curve.  Only the shape term changes: it is a length,
and its first variation involves the scalar curvature.  Nothing else in
Section~\ref{s:var} depends on the dimension, so the variational principle carries
over with these two substitutions.

For the discussion of non-smooth minimizers it is convenient to also allow
the wider class
\begin{equation}\label{e:G2lip}
\widehat{\mathcal G}_2=\bigl\{\eta\in W^{1,\infty}(\Om_0;\R^2):\
\eta\ \text{a homeomorphism of }\overline{\Om_0}\text{ onto }
\overline{\eta(\Om_0)},\ \det D\eta=1\ \text{a.e.}\bigr\} ,
\end{equation}
which contains $\mathcal G_2$.  The lower bounds below are proved for all
$\eta\in\widehat{\mathcal G}_2$, hence in particular for
$\eta\in\mathcal G_2$; the twists of Theorem~\ref{t:twistdeg} belong to this
enlarged class when $\alpha\in W^{1,\infty}$ and $\alpha\notin C^1$.

The relaxed class $\overline{\mathcal C}$ of Definition~\ref{d:relaxed} is a
three-dimensional object and will not be used for planar domains.  Instead we define the planar relaxed class as follows.

\begin{definition}[Planar relaxed class]\label{d:relaxed2}
Set
\[
\begin{aligned}
X_2=\Bigl\{(\Om,B):\ &\Om\subset\R^2,\ |\Om|=V,\ \Per(\Om)<\infty,\\
&B\in L^2(\R^2;\R^2),\ \dive B=0,\ B=0\ \text{a.e.\ on }\R^2\setminus\Om\Bigr\},
\end{aligned}
\]
endowed with the product topology of $L^1$ convergence of characteristic
functions and weak $L^2$ convergence of fields; here $\Om$ is an
almost-everywhere class as in Definition~\ref{d:relaxed}.  Put
\[
\mathcal C_{0,2}^{\mathrm{Lip}}
=\{(\Om_\eta,\Bv_\eta):\eta\in\widehat{\mathcal G}_2\},\qquad
\overline{\mathcal C}_2^{\mathrm{Lip}}
=\operatorname{cl}_{X_2}\bigl(\mathcal C_{0,2}^{\mathrm{Lip}}\bigr).
\]
Here $\widehat{\mathcal G}_2$ is the class of \eqref{e:G2lip}.
On $\mathcal C_{0,2}^{\mathrm{Lip}}$ the Eulerian and Lagrangian energies
agree, so define
\[
E_2(\Om,B)=\tfrac12\int_{\R^2}|B|^2\,\dd y+\kap\Per(\Om),
\qquad
m_{2,\mathrm{Lip}}'=\inf_{\overline{\mathcal C}_2^{\mathrm{Lip}}}E_2 .
\]
Then $\mathcal C_{0,2}^{\mathrm{Lip}}\subseteq
\overline{\mathcal C}_2^{\mathrm{Lip}}$ gives directly
\[
m_{2,\mathrm{Lip}}'\ \le\ \inf_{\widehat{\mathcal G}_2}E .
\]
The closure $\overline{\mathcal C}_2^{\mathrm{Lip}}$ is taken in $X_2$,
whereas $\overline{\mathcal C}$ of Definition~\ref{d:relaxed} is taken in $X$.
\end{definition}

\begin{center}
\begin{tabular}{lll}
\hline
class & defined in & role \\
\hline
$\overline{\mathcal C}$ & \ref{d:relaxed} & three-dimensional relaxed problem \\
$\widehat{\mathcal G}_2$ & \eqref{e:G2lip} & enlarged class, $\mathcal G_2\subset\widehat{\mathcal G}_2$ \\
$\overline{\mathcal C}_2^{\mathrm{Lip}}$ & \ref{d:relaxed2} & planar relaxed problem \\
\hline
\end{tabular}
\end{center}

\begin{proposition}[Compactness for the planar relaxed problem]\label{t:relax2}
Assume that the planar relaxed problem admits a minimizing sequence
which, after translations, lies in a fixed ball.  Then $E_2$ attains
$m_{2,\mathrm{Lip}}'$ on $\overline{\mathcal C}_2^{\mathrm{Lip}}$.
\end{proposition}

\begin{proof}
Lemma~\ref{l:compact} is stated for $d=2,3$ and applies here with $d=2$.
Translations act as homeomorphisms of $X_2$: they preserve $L^1$
convergence of characteristic functions and weak $L^2$ convergence of
fields.  The orbit $\mathcal C_{0,2}^{\mathrm{Lip}}$ and $E_2$ are translation
invariant (Remark~\ref{r:scale}), hence so is the closure
$\overline{\mathcal C}_2^{\mathrm{Lip}}$; the hypothesis therefore provides a
minimizing sequence
$(\Om_n,\Bv_n)\in\overline{\mathcal C}_2^{\mathrm{Lip}}$ with
$E_2(\Om_n,\Bv_n)\to m_{2,\mathrm{Lip}}'$ and $\Om_n\subset B(0,R_0)$ for
all $n$.  Both terms of $E_2$ are nonnegative, hence
$\Per(\Om_n)\le(m_{2,\mathrm{Lip}}'+1)/\kap$ and
$\int_{\R^2}|\Bv_n|^2\le2(m_{2,\mathrm{Lip}}'+1)$ for all large $n$; the
remaining hypotheses of Lemma~\ref{l:compact} hold by
Definition~\ref{d:relaxed2}.  A subsequence therefore converges in $X_2$ to
some $(\Om_\infty,\Bv_\infty)\in X_2$, and lower semicontinuity of both terms
gives
\[
E_2(\Om_\infty,\Bv_\infty)\le\liminf_n E_2(\Om_n,\Bv_n)=m_{2,\mathrm{Lip}}' .
\]
Since $\overline{\mathcal C}_2^{\mathrm{Lip}}$ is closed in $X_2$, the limit
lies in that class, so $E_2(\Om_\infty,\Bv_\infty)\ge m_{2,\mathrm{Lip}}'$;
equality holds and the infimum is attained.
\end{proof}

\begin{theorem}[Two-dimensional variational principle]\label{t:var2d}
Let $\eta\in\mathcal G_2$ and $\Bv=\Bv_\eta$.  Then, with $\curl$ the scalar
curl $\partial_1\Bv_2-\partial_2\Bv_1$ and
$(\curl\Bv)\times\Bv:=(\curl\Bv)J\Bv$,
\begin{equation}\label{e:dE2d}
\delta E(\eta)[v]
=\int_\Om v\cdot\bigl(-(\curl\Bv)\times\Bv\bigr)\,\dd y
+\oint_{\partial\Om}\bigl(\kap\kappa_{\mathrm{curv}}
-\tfrac12|\Bv|^2\bigr)\,v\cdot\nv\,\dd s
\end{equation}
for every divergence-free $v$, where $\kappa_{\mathrm{curv}}$ is the scalar
curvature of the planar curve $\partial\Om$ with respect to the outward
normal $\nv$.  Consequently $\delta E(\eta)=0$ for all divergence-free $v$
if and only if $(\Om,\Bv,\pv)$ satisfies \eqref{e:eq} with the free boundary
condition
\begin{equation}\label{e:bc2d}
\pv+\tfrac12|\Bv|^2=\kap\kappa_{\mathrm{curv}}\quad\text{on }\partial\Om .
\end{equation}
\end{theorem}

\begin{proof}
The identity \eqref{e:identity} of the proof of Theorem~\ref{t:var2} holds
verbatim, since its derivation uses only $\dive\Bv=0$, $\dive v=0$,
$\Bv\cdot\nv=0$ and integration by parts; in the plane $\curl\Bv$ is a scalar
and $(\curl\Bv)\times\Bv$ is the vector field $(\curl\Bv)J\Bv$.  The first
variation of length under a normal displacement $v\cdot\nv$ is
$\delta\Per=\oint_{\partial\Om}\kappa_{\mathrm{curv}}\,v\cdot\nv\,\dd s$,
with $\kappa_{\mathrm{curv}}=1/r$ on a circle of radius $r$.  Adding the two
contributions gives \eqref{e:dE2d}, and the argument of Section~\ref{ss:equiv},
including the harmonic-field step and Lemma~\ref{l:trace}, gives the
equivalence.
\end{proof}

For area-preserving maps the two-by-two differential has a property that
sharpens the discussion of compactness below: its cofactor has the same norm,
so that a gradient bound would also control the inverse.

\begin{lemma}\label{l:2dmatrix}
Let $M\in\R^{2\times2}$ with $\det M=1$.  Then $\cof M=M^{-T}$ and
\begin{equation}
|\cof M|^2=|M|^2,
\end{equation}
where $|M|=(\operatorname{tr}(M^TM))^{1/2}$.  Consequently, for every
$\eta\in\mathcal G_2$,
\[
\int_{\Om_0}|\cof D\eta|^2\,\dd x=\int_{\Om_0}|D\eta|^2\,\dd x,
\]
and the inverse Jacobian satisfies
$\int_{\Om_0}|(D\eta)^{-1}|^2\,\dd x=\int_{\Om_0}|D\eta|^2\,\dd x$.
\end{lemma}

\begin{proof}
Write $M=\begin{pmatrix}a&b\\c&d\end{pmatrix}$ with $\det M=ad-bc=1$.  Its
cofactor matrix is
\[
\cof M=(\det M)M^{-T}=M^{-T}=\begin{pmatrix}d&-c\\-b&a\end{pmatrix},
\]
so $|\cof M|^2=d^2+c^2+b^2+a^2=|M|^2$.  For $\eta\in\mathcal G_2$ one has
$\det D\eta=1$, hence $\cof D\eta=(D\eta)^{-T}$ and, transposing,
$(D\eta)^{-1}=(\cof D\eta)^{T}$ a.e.; since the Frobenius norm is invariant
under transposition, $|(D\eta)^{-1}|=|\cof D\eta|$.  The integral identities
follow.
\end{proof}

\subsection{Weak lower semicontinuity}\label{ss:attain2d}

\begin{lemma}[Weak lower semicontinuity in the plane]\label{t:attain2d}
Let $\eta_n\in\mathcal G_2$ satisfy $E(\eta_n)\to m_2:=\inf_{\mathcal G_2}E$
and suppose that, after translations and a subsequence, $\eta_n$ converge
weakly in $W^{1,2}$ and uniformly on $\overline{\Om_0}$ to an area-preserving
homeomorphism $\eta$ of $\Om_0$ onto $\Om_\infty$, with
$\Om_n\to\Om_\infty$ in $L^1$.  Then
$D\eta_n\Bv_0\rightharpoonup D\eta\Bv_0$ in $L^2$ and
\begin{equation}\label{e:lsc}
E_{W^{1,2}}(\eta):=\frac12\int_{\Om_0}|D\eta\Bv_0|^2\,\dd x
+\kap\Per(\Om_\infty)\ \le\ m_2 .
\end{equation}
The limit therefore supplies a competitor in the weak closure whose energy
is at most $m_2$; no equality of infima is asserted.
\end{lemma}

Here $E_{W^{1,2}}$ is defined only for the limit $\eta$ of the statement,
whose image has finite perimeter and whose push-forward field lies in
$L^2$; a general map in the weak sequential closure need not have a
measurable image domain or a push-forward field.

\begin{proof}
By the assumed weak $W^{1,2}$ convergence, $D\eta_n\rightharpoonup D\eta$ in
$L^2$.  Because $\Bv_0\in L^\infty$, it follows that
$D\eta_n\Bv_0\rightharpoonup D\eta\Bv_0$ weakly in $L^2(\Om_0)$.
Convexity of the magnetic energy gives its weak lower semicontinuity, and by
the assumed $L^1$ convergence of the image domains the perimeter is lower
semicontinuous \cite{AFP2000}.  Therefore
$E_{W^{1,2}}(\eta)\le\liminf_n E(\eta_n)=m_2$, which is \eqref{e:lsc}.
\end{proof}

\subsection{The magnetic energy does not control the full gradient}%
\label{ss:gap2d}

The frozen-in law controls $D\eta$ only in the direction of $\Bv_0$.
Proposition~\ref{p:gap2d} gives an explicit shear with bounded magnetic energy
and unbounded gradient.

\begin{proposition}[The magnetic bound does not control the gradient]%
\label{p:gap2d}
There exist reference configurations $(\Om_0,\Bv_0)$ and
area-preserving embeddings $\eta_n\in\mathcal G_2$ for which
\[
\sup_n\int_{\Om_0}|D\eta_n\Bv_0|^2\,\dd x<\infty
\qquad\text{while}\qquad
\int_{\Om_0}|D\eta_n|^2\,\dd x\to\infty .
\]
\end{proposition}

\begin{proof}
Take $\Om_0=\{x:r_1<|x|<r_2\}$ an annulus and
$\Bv_0(x)=e_\theta=(-x_2/r,x_1/r)$ the unit angular field; it is
divergence free and tangent to both boundary circles, and the number of
boundary components plays no role in the computation.  For a fixed smooth
$f:[r_1,r_2]\to\R$ with $f(r_1)=f(r_2)=0$ define the twist map
$\eta_\alpha(r,\theta)=(r,\theta+\alpha f(r))$, extended to Cartesian
coordinates.  Its pull-back of the area form is
\[
\eta_\alpha^*(r\,\dd r\wedge\dd\theta)
=r\,\dd r\wedge\dd(\theta+\alpha f(r))
=r\,\dd r\wedge\dd\theta ,
\]
because $\dd r\wedge\dd r=0$; thus $\det D\eta_\alpha=1$ and
$\eta_\alpha\in\mathcal G_2$.

In the orthonormal frame $(e_r,e_\theta)$ the differential is
$D\eta_\alpha e_\theta=e_\theta$ and
$D\eta_\alpha e_r=e_r+r\alpha f'(r)e_\theta$; hence
$|D\eta_\alpha\Bv_0|=|D\eta_\alpha e_\theta|=1$ and
\[
\int_{\Om_0}|D\eta_\alpha\Bv_0|^2\,\dd x
=\int_{\Om_0}1\,\dd x
=\pi(r_2^2-r_1^2),
\]
independent of $\alpha$.  On the other hand,
$|D\eta_\alpha|^2
=|D\eta_\alpha e_r|^2+|D\eta_\alpha e_\theta|^2
=1+r^2\alpha^2f'(r)^2+1$, so
\[
\int_{\Om_0}|D\eta_\alpha|^2\,\dd x
=2\pi(r_2^2-r_1^2)
+2\pi\alpha^2\int_{r_1}^{r_2}r^3f'(r)^2\,\dd r\to\infty
\]
as $\alpha\to\infty$ whenever $f'\not\equiv0$.  Since $\eta_\alpha$ preserves
$r$ it maps $\Om_0$ onto itself, so the perimeter is the constant
$\Per(\Om_0)$; the energy \eqref{e:E2d} therefore stays bounded while the
full gradient diverges.  Taking $\alpha=n\to\infty$ gives the sequence
$\eta_n$ of the statement.
\end{proof}

\begin{example}[Unbounded gradients are not excluded by the energy]\label{ex:cav}
The disk with $\Bv_0=\lam x^\perp$ of Example~\ref{ex:2d} is no exception to
Proposition~\ref{p:gap2d}: there $\Bv_0=\lam r e_\theta$ is angular
everywhere, so the rescaled twists $\eta_n(r,\theta)=(r,\theta+n\beta(r))$
apply on the disk as well, with magnetic energy independent of $n$ and
gradient diverging as $n\to\infty$ whenever $\beta'\not\equiv0$.  A twist is
smooth at the origin exactly when
$\beta(r)=\beta_0+h(r^2)$ with $h\in C^\infty([0,1])$, the condition of
Theorem~\ref{t:disk}; such twists attain equality in \eqref{e:magmin}.
More generally, Theorem~\ref{t:cavgen} gives the same
failure of the gradient bound whenever a closed regular field line is
available, while for the circular disk field attainability follows directly
from Theorem~\ref{t:disk}, because the identity disk realizes the global
minimum of $E$.
\end{example}

\begin{theorem}[Energy does not control transverse gradients]\label{t:cavgen}
Let $\Om_0\subset\R^2$ be a bounded domain with smooth boundary and let
$\Bv_0\in C^\infty(\overline{\Om_0},\R^2)$ be divergence free, tangent to
$\partial\Om_0$, and not identically zero.  Assume that $\Bv_0$ possesses a
closed field line in $\Om_0$.  Then there is a sequence $\eta_k\in\mathcal G_2$
of area-preserving diffeomorphisms of $\overline{\Om_0}$ onto itself such that
\begin{equation}\label{e:cavgen}
\int_{\Om_0}|D\eta_k\Bv_0|^2\,\dd x=\int_{\Om_0}|\Bv_0|^2\,\dd x
\qquad\text{for every }k,
\end{equation}
while $\int_{\Om_0}|D\eta_k|^2\,\dd x\to\infty$; in particular
$E(\eta_k)=E(\mathrm{id})$ for every $k$.
\end{theorem}

\begin{proof}
\emph{1. Flow coordinates on a tube.}  Let $\gamma$ be a closed field line of
$\Bv_0$, $\Phi_s$ the flow of $\Bv_0$, $T$ the period of $\gamma$, and
$\Sigma=\{y(n):|n|<\varepsilon_1\}$ a smooth segment through $x_0\in\gamma$
transverse to $\Bv_0$.  Since $\dive\Bv_0=0$, the $1$-form
$-\Bv_{0,2}\dd x_1+\Bv_{0,1}\dd x_2$ is closed and its integral along $\gamma$
vanishes; a tubular neighborhood of $\gamma$ is an annulus whose fundamental
group is generated by $\gamma$, so that form is $-\dd\psi$ there and
$\Bv_0=\nabla^\perp\psi$.  Then $\psi$ is constant along field lines, while
$\psi|_\Sigma$ is strictly monotone, so a field line leaving $y(n)$ returns to
$\Sigma$ only at $y(n)$.  Write the return condition as $F(s,n)=0$ for the
flow $\Phi_s(y(n))$ meeting $\Sigma$; transversality of $\Bv_0$ to $\Sigma$
gives $\partial_sF(T,0)\ne0$, so the implicit function theorem yields a unique
smooth return time $T(n)$ for $|n|$ small, with $T(0)=T$.  Every nearby field
line therefore closes up at $y(n)$ as a smooth embedded closed curve, with period $T(n)$.
Because $\dive\Bv_0=0$ the flow preserves area, so
\[
(s,n)\ \longmapsto\ \Phi_s\bigl(y(n)\bigr),\qquad
s\in\R\ (\text{modulo }T(n)),\quad |n|<\varepsilon_0 ,
\]
is a $C^\infty$ diffeomorphism onto an annular neighborhood $U$ of $\gamma$
foliated by the closed field lines $n=\mathrm{const}$; in these coordinates $\Bv_0=\partial_s$.  The area
element is $\mathcal J(n)\,\dd s\,\dd n$: independence of $s$ follows because
the flow preserves area and commutes with $s$-translations, and
$\mathcal J>0$ on $(-\varepsilon_0,\varepsilon_0)$.

\emph{2. A shear with growing transverse derivative.}  Fix
$\chi\in C_c^\infty((-\varepsilon_0,\varepsilon_0))$ with $0\le\chi\le1$ and
$\chi=1$ on $(-\varepsilon_0/2,\varepsilon_0/2)$, put
$\beta_k(n)=A\sin(kn)\chi(n)$ with $A>0$ fixed, and set
$\tilde\eta_k(s,n)=(s+\beta_k(n),n)$ on $U$.  Then
\begin{equation}\label{e:shearmat}
D\tilde\eta_k=\begin{pmatrix}1 & \beta_k'(n)\\[2pt] 0 & 1\end{pmatrix},
\qquad\text{so}\qquad \det D\tilde\eta_k=1 ,
\end{equation}
and the first column is $(1,0)$, that is,
$D\tilde\eta_k\partial_s=\partial_s=\Bv_0$.

\emph{3. Descent to $\Om_0$.}  Two facts show that $\tilde\eta_k$ descends to
a map of $\Om_0$.  First, $\tilde\eta_k$ respects the identification
$s\sim s+T(n)$, because $\beta_k$ depends on $n$ alone; thus it is a
diffeomorphism of the annulus $U$ onto itself, with inverse
$(s,n)\mapsto(s-\beta_k(n),n)$.  Second, $\chi$ vanishes near
$|n|=\varepsilon_0$, so $\tilde\eta_k$ equals the identity on a collar of
$\partial U$.  Hence $\eta_k:=\tilde\eta_k$ on $U$ and
$\eta_k:=\mathrm{id}$ on $\Om_0\setminus U$ is a smooth diffeomorphism of
$\overline{\Om_0}$ onto itself.  It preserves area: the leaves are mapped to
themselves and, on each leaf, $\eta_k$ is a translation in $s$ modulo $T(n)$,
which preserves the leaf measure $\mathcal J(n)\,\dd s$.  Finally
$D\eta_k\Bv_0=\Bv_0$ on $U$ by \eqref{e:shearmat} and trivially off $U$, which
is \eqref{e:cavgen} and gives $E(\eta_k)=E(\mathrm{id})$.

\emph{4. The gradient blows up.}  On the relatively compact set $U$ the frame
$(\partial_s,\partial_n)$ is comparable to an orthonormal one: there are
$0<a\le b$, depending on $U$ only, with
$a|\xi|\le|\xi_1\partial_s+\xi_2\partial_n|\le b|\xi|$.  Since
$D\tilde\eta_k\partial_n=\beta_k'(n)\partial_s+\partial_n$, this gives
$|D\eta_k|\ge a|\beta_k'(n)|/b$ on $U$, and integrating against the area
element,
\begin{equation}
\int_{\Om_0}|D\eta_k|^2\,\dd x
\ \ge\ c\int\beta_k'(n)^2\,\dd n,\qquad
c=\frac{a^2}{b^2}\min_{|n|\le\varepsilon_0}\bigl(T(n)\mathcal J(n)\bigr)>0 .
\end{equation}
Now $\beta_k'(n)=Ak\cos(kn)\chi(n)+A\sin(kn)\chi'(n)$, so
\[
\int\beta_k'(n)^2\,\dd n
=A^2k^2\int\chi^2\cos^2(kn)\,\dd n+O(1)
\ \sim\ \tfrac12A^2k^2\int\chi(n)^2\,\dd n\to\infty ,
\]
the cross term being $O(1)$ after one integration by parts.  Hence
$\int_{\Om_0}|D\eta_k|^2\,\dd x\to\infty$, as required.
\end{proof}

The equality analysis in the two direct-minimization theorems below is
common to both and is therefore isolated first.

\begin{lemma}[Rigidity in the equality case]\label{l:twistrig}
Let $R>0$, $\Om_0=\{x\in\R^2:|x|<R\}$, and $\eta\in\mathcal G_2$.  Assume
that for every $r\in(0,R)$ and every $\theta$
\begin{equation}\label{e:twistcond}
\bigl|D\eta(re^{i\theta})\,e_\theta\bigr|=1 ,
\end{equation}
and that $\eta(\{|x|<r\})$ is a disk of radius $r$.  Then, up to a rigid
motion and the choice of the origin of $\theta$,
\begin{equation}\label{e:twistform}
\eta(re^{i\theta})=c(r)+r\,e^{i(\theta+\alpha(r))}
\end{equation}
for a $C^1$ function $c$ and a function $\alpha$; and $\eta$ is smooth at the
origin if and only if $\alpha(r)=\alpha_0+h(r^2)$ with $\alpha_0\in\R$ and
$h\in C^\infty([0,R^2])$.
\end{lemma}

\begin{proof}
By \eqref{e:twistcond} the velocity $r\,D\eta(re^{i\theta})e_\theta$ of the
curve $\theta\mapsto\eta(re^{i\theta})$ has constant magnitude $r$.  That
curve is the boundary of the disk $\eta(\{|x|<r\})$, hence a circle of radius
$r$; so it is traversed at constant speed $r$ and winds around it once.
Since $\eta$ is orientation preserving (its Jacobian is $1$) and maps
$\{|x|<r\}$ onto the open disk, the curve is traversed counterclockwise, and
\eqref{e:twistform} follows.  The Jacobian of \eqref{e:twistform} in polar
coordinates is
\[
\det D\eta
=\frac{1}{r}\det\Bigl(\partial_r\eta,\ \partial_\theta\eta\Bigr)
=\frac1r\det\bigl(c'(r)+e^{i(\theta+\alpha(r))}+\cdots,\ r\,i\,
e^{i(\theta+\alpha(r))}\bigr)
=1+c'(r)\cdot e^{i(\theta+\alpha(r))},
\]
where $\cdots$ denotes the terms from differentiating $e^{i(\theta+\alpha(r))}$
with respect to $r$.  Since $\eta$ is area preserving this equals $1$ for
every $r$ and every $\theta$, so $c'(r)\cdot e^{i(\theta+\alpha(r))}=0$ for
all $\theta$; as $\theta$ varies the unit vector $e^{i(\theta+\alpha(r))}$
sweeps the whole unit circle, so $c'(r)=0$.  Thus $c$ is a constant
vector and $\eta$ is a twist
$(r,\theta)\mapsto(r,\theta+\alpha(r))$ up to a rigid motion.

It remains to identify the twists that are smooth at the origin.  For such a
twist $\eta(x)=u(|x|)x+v(|x|)Jx$ with $u=\cos\alpha$ and $v=\sin\alpha$;
here $\eta(x)\cdot x=u(|x|)|x|^2$ and $\eta(x)\cdot Jx=v(|x|)|x|^2$ are
smooth radial functions vanishing to second order at $0$, so $u$ and $v$ are
smooth at $0$.  A smooth radial function on the plane is a smooth function
of $r^2$, its Taylor expansion at $0$ containing only even powers of $r$;
conversely $\eta(x)=\tilde u(|x|^2)x+\tilde v(|x|^2)Jx$ is smooth.  Since
$u+iv$ takes values in $S^1$ and $[0,R]$ is simply connected, this is
equivalent to $\alpha(r)=\alpha_0+h(r^2)$ with $\alpha_0\in\R$ and
$h\in C^\infty([0,R^2])$.
\end{proof}

\begin{theorem}[Rigid rotation on the disk minimizes the magnetic
energy]\label{t:disk}
Let $\Om_0=\Dm$, $\Bv_0=\lam x^\perp$, and assume $\lam\ne0$.
Then for every $\eta\in\mathcal G_2$,
\begin{equation}\label{e:magmin}
\tfrac12\int_{\Dm}|D\eta(x)(\lam x^\perp)|^2\,\dd x
\ge \tfrac{\pi\lam^2}{4}
= \tfrac12\int_{\Dm}|\lam x^\perp|^2\,\dd x ,
\end{equation}
with equality if and only if $\eta$ is a twist
$(r,\theta)\mapsto(r,\theta+\alpha(r))$ up to a rigid motion, where
$\alpha(r)=\alpha_0+h(r^2)$ with $h\in C^\infty([0,1])$; this is the
condition that the twist be smooth at the origin.  Consequently
the identity disk, and more generally every such twist, is a global minimizer of
$E$ over $\mathcal G_2$: the infimum $m_2$ is attained by a genuine
area-preserving embedding, and a complete two-dimensional free-boundary
equilibrium exists.
\end{theorem}

\begin{proof}
Write $\Bv_0(x)=\lam Jx$ with $J=\begin{pmatrix}0&-1\\1&0\end{pmatrix}$, so
that $|D\eta\Bv_0|^2=\lam^2|D\eta Jx|^2$.  We identify $\R^2$ with $\C$ when
evaluating $J$ on $e^{i\theta}$, so that $Je^{i\theta}=i\,e^{i\theta}$; thus
$J(re^{i\theta})=ri\,e^{i\theta}$.  For fixed $r\in(0,1)$ let
$L(r)$ be the length of the image curve $\eta(\{x:|x|=r\})$.

Parametrizing this circle by $\theta$ gives the velocity
$\dd\eta(re^{i\theta})/\dd\theta=D\eta(re^{i\theta})(ri\,e^{i\theta})$;
since $ri\,e^{i\theta}=J(re^{i\theta})=Jx$ has magnitude $r$, its magnitude is
$|D\eta Jx|$,
and therefore
\[
L(r)=\int_0^{2\pi}|D\eta Jx|\,\dd\theta .
\]
By the Cauchy--Schwarz inequality on the circle,
\[
\int_0^{2\pi}|D\eta Jx|^2\,\dd\theta
\ge \frac{L(r)^2}{2\pi} .
\]
Integrating in $r$ against the area element $r\,\dd r$ yields
\[
\int_{\Dm}|D\eta Jx|^2\,\dd x
=\int_0^1 r\int_0^{2\pi}|D\eta Jx|^2\,\dd\theta\,\dd r
\ge \frac1{2\pi}\int_0^1 r\,L(r)^2\,\dd r .
\]
Now define $\varrho(y)=|\eta^{-1}(y)|$ on $\Om=\eta(\Dm)$, so that the
level set $\{\varrho=r\}$ is exactly $\eta(\{|x|=r\})$, of length $L(r)$;
and $\{\varrho<r\}=\eta(\{|x|<r\})$ is a topological disk of area $\pi r^2$.

Its perimeter is $L(r)$, and the isoperimetric inequality gives
$L(r)\ge 2\pi r$.  Hence
\[
\int_0^1 r\,L(r)^2\,\dd r
\ge \int_0^1 r\,(2\pi r)^2\,\dd r
=4\pi^2\int_0^1 r^3\,\dd r=\pi^2 ,
\]
and consequently
\[
\int_{\Dm}|D\eta Jx|^2\,\dd x\ge \frac{\pi^2}{2\pi}=\frac{\pi}{2},
\qquad
\frac12\int_{\Dm}|D\eta\Bv_0|^2\,\dd x\ge \frac{\pi\lam^2}{4}.
\]
Equality in the
Cauchy--Schwarz inequality on the circle $\{|x|=r\}$ requires
$|D\eta(re^{i\theta})(ri\,e^{i\theta})|$ to be independent of $\theta$, say
equal to $a(r)$; then $L(r)=2\pi a(r)$.  Equality in the isoperimetric
inequality requires the image curve, of perimeter $L(r)$ and enclosing area
$\pi r^2$, to be a circle of radius $r$, that is $L(r)=2\pi r$; comparing the
two gives $a(r)=r$, whence $|D\eta(re^{i\theta})e_\theta|=1$ and the image of
$\{|x|<r\}$ is a disk of radius $r$.  These conclusions hold for almost every
$r$; since $(r,\theta)\mapsto D\eta(re^{i\theta})e_\theta$ and
$(r,\theta)\mapsto\eta(re^{i\theta})$ are continuous, they hold for every
$r\in(0,1)$.  Lemma~\ref{l:twistrig} then shows that
$\eta$ is a twist $(r,\theta)\mapsto(r,\theta+\alpha(r))$ up to a rigid
motion, smooth at the origin exactly when $\alpha(r)=\alpha_0+h(r^2)$ with
$\alpha_0\in\R$ and $h\in C^\infty([0,1])$.

Finally,
$E(\eta)=\tfrac12\int_{\Dm}|D\eta\Bv_0|^2\,\dd x
+\kap\Per(\eta(\Dm))$.  The magnetic term is at least
$\pi\lam^2/4$ by \eqref{e:magmin}, and $\Per(\eta(\Dm))\ge 2\pi$ by
the isoperimetric inequality because $|\eta(\Dm)|=\pi$.  Thus
$E(\eta)\ge \pi\lam^2/4+2\pi\kap$ for every $\eta\in\mathcal G_2$, while the
identity (or any twist) achieves equality.  Hence the identity disk is a
global minimizer and the infimum is attained by a genuine area-preserving
embedding.
\end{proof}

\subsection{Regularity of a minimizer}\label{ss:reg}

The compactness hypothesis for a minimizing sequence $(\eta_n)$ is the
uniform gradient bound
\begin{equation}\label{e:gradbd}
\sup_n\int_{\Om_0}|D\eta_n(x)|^2\,\dd x<\infty .
\end{equation}

\begin{theorem}[The twist degeneracy on the disk]\label{t:twistdeg}
Let $\Om_0=\Dm$ and $\Bv_0=\lam x^\perp$.  For every
$\alpha\in W^{1,\infty}(0,1)$ the map
$\eta_\alpha(r,\theta)=(r,\theta+\alpha(r))$ is an area-preserving
homeomorphism of $\Dm$ onto itself belonging to the enlarged class
$\widehat{\mathcal G}_2$ of \eqref{e:G2lip}; it
satisfies
$\eta_\alpha(\Dm)=\Dm$
and achieves equality in \eqref{e:magmin}, so it
attains the infimum of $E$ over $\widehat{\mathcal G}_2$, while the same
value is attained over $\mathcal G_2$ by the identity.  If $\alpha'$ is
discontinuous at some $r_0\in(0,1)$ then
$\eta_\alpha$ is a minimizer whose gradient is bounded yet which is not of
class $C^1$.
\end{theorem}

\begin{proof}
Throughout write $x=re^{i\theta}$ with $r\in[0,1]$ and $\theta\in\R$.
For $r>0$ the map $\eta_\alpha$ is given by
$\eta_\alpha(re^{i\theta})=re^{i(\theta+\alpha(r))}$, and we set
$\eta_\alpha(0)=0$; this defines a single-valued map of $\Dm$ because
$\theta\mapsto\theta+\alpha(r)$ preserves the identification
$\theta\sim\theta+2\pi$ and the value at $r=0$ is independent of $\theta$.

\medskip
\noindent\textit{Step 1: $\eta_\alpha$ is a homeomorphism belonging to
$\widehat{\mathcal G}_2$.}
Continuity follows from the continuity of $r\mapsto\alpha(r)$: on $\Dm$ one
has
\begin{align*}
|\eta_\alpha(re^{i\theta})-\eta_\alpha(r'e^{i\theta'})|
&\le |re^{i(\theta+\alpha(r))}-r'e^{i(\theta+\alpha(r))}|+|r'e^{i(\theta+\alpha(r))}-r'e^{i(\theta'+\alpha(r))}|\\
  &\; +|r'e^{i(\theta'+\alpha(r))}-r'e^{i(\theta'+\alpha(r'))}|\\
&\le |r-r'|+r'\,|\theta-\theta'|+r'\,|\alpha(r)-\alpha(r')| ,
\end{align*}
which tends to $0$ as $(r',\theta')\to(r,\theta)$ since $\alpha$ is
Lipschitz.  The map $\beta(r,\theta):=(r,\theta-\alpha(r))$, equivalently
$\beta(re^{i\theta})=re^{i(\theta-\alpha(r))}$, is defined in the same way and
is likewise continuous; directly from the definition,
$\eta_\alpha\circ\beta=\beta\circ\eta_\alpha=\mathrm{id}$ on $\Dm$.  Hence
$\eta_\alpha$ is a homeomorphism of $\Dm$ onto itself with inverse $\beta$.
Moreover $\alpha\in W^{1,\infty}$ makes $\eta_\alpha\in W^{1,\infty}$, and the
computation
\[
\det D\eta_\alpha
=\tfrac1r\det\bigl(\partial_r\eta_\alpha,\partial_\theta\eta_\alpha\bigr)
=\tfrac1r\det\bigl(e^{i(\theta+\alpha)}
+r\alpha'(r)\,ie^{i(\theta+\alpha)},\ ri\,e^{i(\theta+\alpha)}\bigr)
=1
\]
shows $\det D\eta_\alpha=1$ a.e.; thus
$\eta_\alpha\in\widehat{\mathcal G}_2$, the class defined in
\eqref{e:G2lip}.

\medskip
\noindent\textit{Step 2: the magnetic inequality \eqref{e:magmin} holds on
$\widehat{\mathcal G}_2$.}
The proof of Theorem~\ref{t:disk} uses only two ingredients, neither of which
needs $\eta$ to be smooth.  First, for $\eta\in\widehat{\mathcal G}_2$ the set
$\eta(\{|x|<r\})$ is a topological disk whose area is $\pi r^2$, because
$\eta$ is a homeomorphism with $\det D\eta=1$ a.e. and $|\{|x|<r\}|=\pi r^2$;
hence its boundary $\eta(\{|x|=r\})$, of length $L(r)$, satisfies
$L(r)\ge2\pi r$ by the isoperimetric inequality for rectifiable curves.
Second, $\eta\in W^{1,\infty}$ makes $\theta\mapsto\eta(re^{i\theta})$
absolutely continuous for a.e.\ $r$, so that
$L(r)=\int_0^{2\pi}|D\eta(re^{i\theta})(ri\,e^{i\theta})|\,\dd\theta$ and the
Cauchy--Schwarz inequality
$\int_0^{2\pi}|D\eta Jx|^2\,\dd\theta\ge L(r)^2/(2\pi)$ holds for a.e.\ $r$.
Integrating as in the proof of Theorem~\ref{t:disk} and using
$|\Bv_0(x)|=|\lam|\,r$ then gives
\[
\frac12\int_{\Dm}|D\eta\Bv_0|^2\,\dd x\ \ge\ \frac{\pi\lam^2}{4}
\qquad\text{for every }\eta\in\widehat{\mathcal G}_2 ,
\]
i.e.\ \eqref{e:magmin} with $\mathcal G_2$ replaced by
$\widehat{\mathcal G}_2$.

\medskip
\noindent\textit{Step 3: $\eta_\alpha$ saturates
\eqref{e:magmin}.}
At a point $re^{i\theta}$ with $r>0$, differentiating
$\eta_\alpha(re^{i\theta})=re^{i(\theta+\alpha(r))}$ in the orthonormal frame
$(e_r,e_\theta)$ gives, with $e_r,e_\theta$ evaluated at the image angle
$\theta+\alpha(r)$,
\[
D\eta_\alpha\,e_r=e_r+r\alpha'(r)\,e_\theta,\qquad
D\eta_\alpha\,e_\theta=e_\theta .
\]
Hence $|D\eta_\alpha|^2=2+r^2\alpha'^2$, and since
$\Bv_0(x)=\lam x^\perp=\lam r\,e_\theta$ (with $e_\theta$ now at the source
angle $\theta$) we get
\[
D\eta_\alpha\Bv_0=\lam r\,D\eta_\alpha e_\theta
=\lam r\,e_\theta ,\qquad
|D\eta_\alpha\Bv_0|=|\lam|\,r=|\Bv_0| ,
\]
the last equality being independent of $\theta$.  Therefore
$\tfrac12\int_{\Dm}|D\eta_\alpha\Bv_0|^2\,\dd x
=\tfrac12\int_{\Dm}|\Bv_0|^2\,\dd x=\pi\lam^2/4$, so $\eta_\alpha$ attains
equality in \eqref{e:magmin} over $\widehat{\mathcal G}_2$.  Since
$\eta_\alpha(\Dm)=\Dm$ (Step~1),
$\Per(\eta_\alpha(\Dm))=\Per(\Dm)=2\pi$.  Consequently
$E(\eta_\alpha)=\pi\lam^2/4+2\pi\kap$, which by Step~2 is the least value of
$E$ on $\widehat{\mathcal G}_2$; hence $\eta_\alpha$ is a global minimizer of
$E$ over $\widehat{\mathcal G}_2$.  The identity (equivalently, any twist)
attains the same value over the smaller class $\mathcal G_2$.

\medskip
\noindent\textit{Step 4: failure of smoothness.}
Assume now $\alpha'\in L^\infty$ is discontinuous at some $r_0\in(0,1)$.  On
one hand $r^2\alpha'(r)^2\le\|\alpha'\|_\infty^2$, so
$|D\eta_\alpha|^2=2+r^2\alpha'^2$ is bounded and the gradient bound
\eqref{e:gradbd} holds.  On the other hand
$D\eta_\alpha e_r=e_r+r\alpha'(r)e_\theta$ is discontinuous along the circle
$\{|x|=r_0\}$, so $D\eta_\alpha$ has no continuous representative there and
$\eta_\alpha\notin C^1$; the map
$\alpha(r)=c\min\{r,r_0\}$ with $c\ne0$ realizes this.  Thus a global
minimizer on which \eqref{e:gradbd} holds need not be of class $C^1$.
\end{proof}

\subsection{Direct minimization for circular fields and its sharpness}%
\label{ss:circ}

Theorem~\ref{t:disk} is not an accident of the unit disk: the same direct
minimization works for every \emph{circular} reference field on a disk,
namely one whose field lines are the concentric circles.  Such fields are
exactly the radial--angular fields $\Bv_0(x)=\rho(|x|)e_\theta$ on a disk.
On an annulus the level sets are again circles, but the isoperimetric step
fails and the identity is no longer a minimizer (Remark~\ref{r:annulus}), so
the disk hypothesis cannot be dropped.

\begin{theorem}[Circular fields on a disk minimize directly]\label{t:circ}
Let $\Om_0=\{x:|x|<R\}$ be the disk of radius $R>0$ and let
$\Bv_0(x)=\rho(|x|)e_\theta$, where $\rho(r)=r q(r^2)$ for some
$q\in C^\infty([0,R^2])$ with $q\ne0$ on $[0,R^2]$; then $\rho$ has a
constant sign on $(0,R]$ and $\Bv_0$ extends smoothly across the axis, being
tangent to $\partial\Om_0$.  Then for every $\eta\in\mathcal G_2$,
\begin{equation}\label{e:circmin}
\frac12\int_{\Om_0}|D\eta\Bv_0|^2\,\dd x
\ge \frac12\int_{\Om_0}|\Bv_0|^2\,\dd x ,
\end{equation}
with equality if and only if $\eta$ is a twist
$(r,\theta)\mapsto(r,\theta+\alpha(r))$ up to a rigid motion, with
$\alpha(r)=\alpha_0+h(r^2)$ for some $h\in C^\infty([0,R^2])$; this is the
smoothness condition at the origin characterized in
Lemma~\ref{l:twistrig}.  Consequently
the identity (and more generally every twist) is a global minimizer of $E$
over $\mathcal G_2$: the infimum is attained by a genuine area-preserving
embedding, and a complete two-dimensional free-boundary equilibrium exists.
\end{theorem}

Theorem~\ref{t:disk} covers the case $\rho(r)=\lam r$ on the unit disk with
$\lam\ne0$; the sign of $\lam$ is immaterial, since the energy is quadratic
in $\Bv_0$ and $(\curl(-\Bv))\times(-\Bv)=(\curl\Bv)\times\Bv$, so that
$\Bv_0\mapsto-\Bv_0$ leaves $E$, its minimizers and their equality case
unchanged.

\begin{proof}
Write $x=re^{i\theta}$ with $r\in(0,R)$ and $\theta\in\R$, and let
$(e_r,e_\theta)$ be the polar frame at $x$, so that
$\Bv_0(x)=\rho(r)e_\theta$ and $\dd x=r\,\dd r\,\dd\theta$.

\medskip
\textit{Step 1: reduction to a one-variable inequality.}
Since $\Bv_0(x)=\rho(r)e_\theta$ we have
$|D\eta(x)\Bv_0(x)|^2=\rho(r)^2\,|D\eta(x)e_\theta|^2$, hence
\begin{equation}\label{e:circsplit}
\int_{\Om_0}|D\eta\Bv_0|^2\,\dd x
=\int_0^R\rho(r)^2\Bigl(\int_0^{2\pi}|D\eta(re^{i\theta})e_\theta|^2
\,\dd\theta\Bigr)r\,\dd r .
\end{equation}
For $r\in(0,R)$ set $g(r,\theta):=|D\eta(re^{i\theta})e_\theta|$ and let
$L(r)$ be the length of the image curve $\eta(\{|x|=r\})$.  Parametrizing
this curve by $\theta$, its velocity is
$\partial_\theta\eta(re^{i\theta})=r\,D\eta(re^{i\theta})e_\theta$, of
magnitude $r\,g(r,\theta)$; therefore
\begin{equation}\label{e:circL}
L(r)=\int_0^{2\pi}r\,g(r,\theta)\,\dd\theta
=r\int_0^{2\pi}g(r,\theta)\,\dd\theta .
\end{equation}

\medskip
\textit{Step 2: the magnetic energy is bounded below.}
The Cauchy--Schwarz inequality on the circle gives, in view of \eqref{e:circL},
\[
\int_0^{2\pi}g(r,\theta)^2\,\dd\theta
\ \ge\ \frac1{2\pi}\Bigl(\int_0^{2\pi}g(r,\theta)\,\dd\theta\Bigr)^2
=\frac1{2\pi}\frac{L(r)^2}{r^2}=\frac{L(r)^2}{2\pi r^2}.
\]
Substituting into \eqref{e:circsplit},
\begin{equation}\label{e:circLB}
\int_{\Om_0}|D\eta\Bv_0|^2\,\dd x
\ \ge\ \frac1{2\pi}\int_0^R\rho(r)^2\,\frac{L(r)^2}{r}\,\dd r .
\end{equation}
On the other hand $\eta\in\mathcal G_2$ is area preserving, so
$\eta(\{|x|<r\})$ is a topological disk of area $\pi r^2$ and its boundary is
the curve $\eta(\{|x|=r\})$, of length $L(r)$.  The isoperimetric inequality
therefore gives $L(r)\ge2\pi r$, the radius of the area-minimizing circle
being $r$.  Inserting this into \eqref{e:circLB},
\begin{equation}\label{e:circfin}
\int_{\Om_0}|D\eta\Bv_0|^2\,\dd x
\ \ge\ \frac1{2\pi}\int_0^R\rho(r)^2\,\frac{4\pi^2r^2}{r}\,\dd r
=2\pi\int_0^R\rho(r)^2\,r\,\dd r
=\int_{\Om_0}|\Bv_0|^2\,\dd x ,
\end{equation}
the last equality because $|\Bv_0(x)|^2=\rho(r)^2$ and
$\dd x=r\,\dd r\,\dd\theta$ give
$\int_{\Om_0}|\Bv_0|^2\,\dd x
=2\pi\int_0^R\rho(r)^2\,r\,\dd r$.  This is \eqref{e:circmin}.

\medskip
\textit{Step 3: equality case.}
Equality in \eqref{e:circfin} forces equality in Step 2 at almost every $r$:
equality in the Cauchy--Schwarz inequality makes $g(r,\cdot)$ constant, say
$g(r,\theta)=g(r)$, so that \eqref{e:circL} gives $L(r)=2\pi r\,g(r)$, and
equality in the isoperimetric inequality makes the image curve a circle of
radius $r$, so $L(r)=2\pi r$ and $g(r)=1$.  Since $g$ and $L$ are
continuous in $r$, both equalities hold for every $r$, and then
$\eta(\{|x|<r\})$ is a disk of radius $r$; Lemma~\ref{l:twistrig} gives that $\eta$ is a twist up to a rigid
motion, smooth at the origin exactly when $\alpha(r)=\alpha_0+h(r^2)$ with
$h\in C^\infty([0,R^2])$.  Conversely, for every such twist
$|D\eta(re^{i\theta})e_\theta|=1$, so equality holds in \eqref{e:circmin},
which completes the characterization.

\medskip
\textit{Step 4: minimality of $E$.}
The final step is again the same as in Theorem~\ref{t:disk}.  For
$\eta\in\mathcal G_2$,
$E(\eta)=\tfrac12\int_{\Om_0}|D\eta\Bv_0|^2\,\dd x+\kap\,\Per(\eta(\Om_0))$;
the magnetic term is at least $\tfrac12\int_{\Om_0}|\Bv_0|^2\,\dd x$ by
\eqref{e:circmin}, while $\Per(\eta(\Om_0))\ge2\pi R$ by the isoperimetric
inequality, since $\eta$ is area preserving and $|\eta(\Om_0)|=\pi R^2$ (in
Theorem~\ref{t:disk} this reads $\Per(\eta(\Dm))\ge2\pi$, the case $R=1$).
Hence $E(\eta)\ge\tfrac12\int_{\Om_0}|\Bv_0|^2\,\dd x+2\pi\kap R=E(\mathrm{id})$
for every $\eta\in\mathcal G_2$, with equality for $\eta=\mathrm{id}$.  Thus
the identity is a global minimizer over $\mathcal G_2$, and by Step 3 the
minimizers are exactly the twists up to rigid motions; the infimum is
attained by a genuine area-preserving embedding.
\end{proof}

\begin{remark}[The disk hypothesis is necessary]\label{r:annulus}
For an annulus the identity is not a minimizer.  Let
$\Om_0=\{R_1<|x|<R_2\}$ with $0<R_1<R_2$, let $\Bv_0=\rho(r)e_\theta$ with
$\rho\not\equiv0$, take $-R_1^2<C<0$, and set
$\eta_C(r,\theta)=(\sqrt{r^2+C},\theta)$.  Here $f(r)=\sqrt{r^2+C}$ is smooth
and positive on $[R_1,R_2]$ and satisfies $f(r)<r$ together with $ff'=r$;
the latter gives $\det D\eta_C=\frac1r ff'=1$, so $\eta_C$ is an
area-preserving diffeomorphism of $\Om_0$ onto the annulus
$\{\sqrt{R_1^2+C}<|x|<\sqrt{R_2^2+C}\}$, hence an element of $\mathcal G_2$.
Since $f(r)<r$ we have $|D\eta_Ce_\theta|=f(r)/r<1$ wherever $\rho\ne0$, so
the magnetic energy of $\eta_C$ is strictly smaller than that of the
identity; the perimeter of the image is likewise smaller, namely
$2\pi\bigl(f(R_1)+f(R_2)\bigr)<2\pi(R_1+R_2)$, because $f(r)=\sqrt{r^2+C}<r$
for $r>0$ when $C<0$.  Thus $E(\eta_C)<E(\mathrm{id})$,
so the identity is not a minimizer on an annulus; the direct argument of
Theorem~\ref{t:circ} is specific to the disk, where the enclosed area is
determined by the level.
\end{remark}

\subsection{Existence for a general single-axis reference field}\label{ss:rearr}

For a single-axis stream function, circular rearrangement gives a sharp
lower bound and an explicit minimizer.  We prove the comparison, show that
the rearranging map is a smooth area-preserving diffeomorphism, and obtain
the existence of a smooth equilibrium in the original configuration class
$\mathcal G_2$ for every reference field satisfying {\rm (H1)}.

The method is the classical rearrangement method for steady ideal flows and
vortex patches \cite{Burton1989}.  The frozen-in law fixes the distribution
function of the stream function, and the Schwarz rearrangement minimizes the
Dirichlet energy among functions with that distribution.  Here this gives the
circular profile \eqref{e:psistar}; the free boundary also minimizes the
perimeter at fixed area.  The rearranging map is shown below to be a smooth
area-preserving diffeomorphism, hence a member of the original class
$\mathcal G_2$; the comparison is therefore carried out in $\mathcal G_2$
itself, and the infimum is attained there.

\emph{The enlarged admissible class.}  The lower bounds below hold for every
$\eta$ in the enlarged class $\widehat{\mathcal G}_2$ of \eqref{e:G2lip},
while the minimizing map constructed below lies in $\mathcal G_2$.

\emph{Standing assumptions.}  Let $\Om_0\subset\R^2$ be bounded and simply
connected with $C^\infty$ boundary and area $V$, and let
$\Bv_0\in C^\infty(\overline\Om_0;\R^2)$ satisfy $\dive\Bv_0=0$,
$\Bv_0\cdot\nv_0=0$ on $\partial\Om_0$, and $\Bv_0\not\equiv0$.  Since
$\Om_0$ is simply connected there is a stream function
$\psi_0\in C^\infty(\overline\Om_0)$, unique up to an additive constant, with
$\Bv_0=\nabla^\perp\psi_0$, and tangency is equivalent to $\psi_0$ being
constant on $\partial\Om_0$.  The problem is invariant under
$\Bv_0\mapsto-\Bv_0$: the energy \eqref{e:E2d} is quadratic in $\Bv_0$ and
$(\curl(-\Bv))\times(-\Bv)=(\curl\Bv)\times\Bv$; by this invariance a single
nondegenerate extremum may be taken to be a minimum.  We assume
\begin{enumerate}
\item[(H1)] $\psi_0$ has exactly one critical point $p_0$ in
$\overline\Om_0$; it lies in $\Om_0$ and is a nondegenerate minimum.
\end{enumerate}
This is the natural two-dimensional form of the hypothesis that the field has
a single magnetic axis: all field lines are closed curves nested about $p_0$.

Write $m=\min_{\overline{\Om_0}}\psi_0$ and let $M=\max_{\overline{\Om_0}}\psi_0$
be the maximum of $\psi_0$ over $\overline{\Om_0}$.  Under {\rm (H1)} the only
critical point is the minimum $p_0$, so $m=\psi_0(p_0)$ and the maximum is
attained on $\partial\Om_0$, where $\psi_0=M>m$.  We fix the additive
constant in $\psi_0$ by this normalization, so that the boundary value is $M$
throughout.  For $t\in(m,M)$ the level set
$\{\psi_0=t\}$ contains no critical point, hence is a smooth compact
$1$-manifold in $\Om_0$, i.e.\ a finite disjoint union of circles.  By the
Morse lemma $\{\psi_0<t\}$ is a small disk about $p_0$ for $t$ slightly above
$m$, and flowing along the smooth vector field
$\nabla\psi_0/|\nabla\psi_0|^2$ gives a diffeomorphism
$\{\psi_0\le t_1\}\to\{\psi_0\le t_2\}$ for all $m<t_1<t_2<M$.  Hence
\begin{equation}\label{e:subdisc}
\text{for every }t\in(m,M),\quad
\{\psi_0<t\}\ \text{is a topological disk and}\ \{\psi_0=t\}
\ \text{is a single circle.}
\end{equation}
Define the \emph{area function} of the reference field,
\begin{equation}
A_0(t)=\bigl|\{x\in\Om_0:\psi_0(x)<t\}\bigr|,\qquad t\in[m,M],
\end{equation}
so that $A_0(m)=0$ and $A_0(M)=V$; by the coarea formula
\begin{equation}\label{e:A0prime}
A_0'(t)=\oint_{\{\psi_0=t\}}\frac{1}{|\nabla\psi_0|}\,\dd\mathcal H^1>0 ,
\end{equation}
so $A_0$ is a smooth, strictly increasing bijection of $[m,M]$ onto $[0,V]$.
(The coarea formula \eqref{e:A0prime} is applied at levels $t\in(m,M)$, where
$\nabla\psi_0\ne0$; the smoothness at the endpoint $t=m$, where
$\nabla\psi_0$ vanishes at $p_0$, together with $A_0'(m)>0$, is established
separately in Lemma~\ref{l:morseiso}.)

\begin{lemma}[Conjugation and distribution invariance]\label{l:conj}
Let $\eta\in\widehat{\mathcal G}_2$, $\Om=\eta(\Om_0)$, $\Bv=\eta_*\Bv_0$ and
$\psi=\psi_0\circ\eta^{-1}$.  Then
\begin{enumerate}[\rm (i)]
\item $\Bv=\nabla^\perp\psi$ a.e.\ in $\Om$;
\item $|\{y\in\Om:\psi(y)<t\}|=A_0(t)$ for every $t\in[m,M]$;
\item $\int_\Om|\Bv|^2\,\dd y=\int_{\Om_0}|D\eta\,\Bv_0|^2\,\dd x$.
\end{enumerate}
\end{lemma}

\begin{proof}
(iii) is the change of variables $y=\eta(x)$: since $\det D\eta=1$ a.e.\ we
have $\dd y=\dd x$, and $|\Bv(\eta(x))|^2=|D\eta(x)\Bv_0(x)|^2$.

For (i) we first note that a $2\times2$ matrix $A$ with $\det A=1$
satisfies $AJ=JA^{-T}$.  Writing $A=\begin{pmatrix}a&b\\c&d\end{pmatrix}$
with $ad-bc=1$ one has $A^{-1}=\begin{pmatrix}d&-b\\-c&a\end{pmatrix}$, hence
$A^{-T}=\begin{pmatrix}d&-c\\-b&a\end{pmatrix}$, and therefore
\[
JA^{-T}=
\begin{pmatrix}0&-1\\1&0\end{pmatrix}
\begin{pmatrix}d&-c\\-b&a\end{pmatrix}
=\begin{pmatrix}b&-a\\d&-c\end{pmatrix}
=AJ .
\]
Consequently, at $y=\eta(x)$,
\[
\Bv(y)=D\eta(x)\Bv_0(x)=D\eta(x)J\nabla\psi_0(x)
=J\bigl(D\eta(x)\bigr)^{-T}\nabla\psi_0(x) ,
\]
whereas
$\nabla\psi(y)=\nabla(\psi_0\circ\eta^{-1})(y)
=\bigl(D\eta^{-1}(y)\bigr)^T\nabla\psi_0(x)
=\bigl(D\eta(x)\bigr)^{-T}\nabla\psi_0(x)$.

A homeomorphism $\eta\in W^{1,\infty}$ with $J_\eta=\det D\eta=1$ a.e.\ is a
planar mapping of finite distortion, with
$K=|D\eta|^2/J_\eta\le\|D\eta\|_\infty^2$.  By the inverse theorem for
planar Sobolev homeomorphisms with positive Jacobian
\cite{HenclKoskela2014}, $\eta^{-1}\in W^{1,1}_{\mathrm{loc}}$ and
$D\eta^{-1}(\eta(x))=(D\eta(x))^{-1}$ for a.e.\ $x$.  In the plane
$|(D\eta)^{-1}|=|\cof D\eta|=|D\eta|$, and the area formula for the
Lipschitz homeomorphism $\eta$ gives, for every $p<\infty$,
\[
\int_\Om|D\eta^{-1}|^p\,\dd y
=\int_{\Om_0}|(D\eta)^{-1}|^p\,\dd x
\le\|D\eta\|_\infty^p|\Om_0| .
\]
Letting $p\to\infty$ in the last equation gives
$\|D\eta^{-1}\|_{L^\infty(\Om)}\le\|D\eta\|_\infty$, so
$\eta^{-1}\in W^{1,\infty}(\Om)$ and the Sobolev chain rule
applies to $\psi_0\circ\eta^{-1}$.  Hence $\Bv=J\nabla\psi=\nabla^\perp\psi$
a.e.\ and $\psi\in W^{1,\infty}(\Om)\cap C^0(\overline\Om)$, as claimed; in
particular $\psi$ is locally Lipschitz.

For (ii), $\{y\in\Om:\psi(y)<t\}=\eta(\{x\in\Om_0:\psi_0(x)<t\})$, and
$|\eta(E)|=|E|$ for every measurable $E\subset\Om_0$ because $\eta$ is a
homeomorphism with $\det D\eta=1$ a.e.
\end{proof}

Thus, the magnetic energy of any competitor is the Dirichlet energy of a
function $\psi$ on $\Om$ whose distribution function is \emph{fixed}, equal
to that of $\psi_0$; the only freedom left is the shape of the level sets.
The next lemma states that among all such $\psi$ the one with circular level
sets has the least Dirichlet energy.  This is the P\'olya--Szeg\H{o}
principle for the Schwarz symmetrization
\cite{PolyaSzego1951,Talenti1976,BrothersZiemer1988}; we recall the short
proof, which also exhibits the mechanism.  For every competitor
$\eta\in\widehat{\mathcal G}_2$ the function $\psi=\psi_0\circ\eta^{-1}$
lies in $W^{1,\infty}(\Om)$ by Lemma~\ref{l:conj}(i).  Since $\eta$ is bi-Lipschitz, the
chain rule gives $\nabla\psi=(D\eta^{-1})^{\mathsf T}(\nabla\psi_0\circ\eta^{-1})$
a.e., so $\nabla\psi=0$ a.e.\ only on $\eta(\operatorname{Crit}\psi_0)$; under {\rm (H1)} or {\rm (H2)}, and
in general under the assumptions of Section~\ref{ss:gen} (all interior critical
points of $\psi_0$ nondegenerate and $\nabla\psi_0\ne0$ on $\partial\Om_0$), that set is
finite and hence Lebesgue-null.  Thus $|\{\nabla\psi=0\}|=0$ and $\nabla\psi\ne0$
$\mathcal H^1$-a.e.\ on $\{\psi=t\}$ for a.e.\ $t$, so the coarea formula and the level-set
estimates below apply to $\psi$ for a.e.\ level $t$ throughout
Section~\ref{ss:multi} and Section~\ref{ss:gen}.

\begin{lemma}[Level-set lower bound]\label{l:ps}
Let $\Om\subset\R^2$ be a bounded domain and let $\psi\in W^{1,\infty}(\Om)\cap
C^0(\overline\Om)$ satisfy $\psi=M$ on $\partial\Om$, $\min\psi=m$,
$\psi\le M$ in $\overline\Om$, $|\{\nabla\psi=0\}|=0$, and
$|\{y\in\Om:\psi(y)<t\}|=A_0(t)$ for every $t\in[m,M]$.  Then
\begin{equation}\label{e:ps}
\int_\Om|\nabla\psi|^2\,\dd y\ \ge\
\int_m^M\frac{4\pi A_0(t)}{A_0'(t)}\,\dd t .
\end{equation}
\end{lemma}

\begin{proof}
Since $\psi\in W^{1,\infty}(\Om)$, the coarea formula for Sobolev functions
\cite{EvansGariepy2025} gives both
\begin{equation}\label{e:coarea}
\int_\Om|\nabla\psi|^2\,\dd y
=\int_\R\Bigl(\int_{\{\psi=t\}}|\nabla\psi|\,\dd\mathcal H^1\Bigr)\dd t
\end{equation}
and, since $|\{\nabla\psi=0\}|=0$, on differentiating $|\{\psi<t\}|=A_0(t)$,
\begin{equation}\label{e:A0p}
A_0'(t)=\int_{\{\psi=t\}}\frac{1}{|\nabla\psi|}\,\dd\mathcal H^1
\qquad\text{for a.e. }t\in(m,M) .
\end{equation}
Fix such a $t$ and let $\ell(t)=\mathcal H^1(\{\psi=t\})$.  By
the Cauchy--Schwarz inequality on the level set,
\[
\Bigl(\int_{\{\psi=t\}}|\nabla\psi|\,\dd\mathcal H^1\Bigr)
\Bigl(\int_{\{\psi=t\}}\frac{1}{|\nabla\psi|}\,\dd\mathcal H^1\Bigr)
\ \ge\ \ell(t)^2 ,
\]
so by \eqref{e:A0p}
\begin{equation}\label{e:csleaves}
\int_{\{\psi=t\}}|\nabla\psi|\,\dd\mathcal H^1\ \ge\
\frac{\ell(t)^2}{A_0'(t)} .
\end{equation}
Since $\psi=M$ on $\partial\Om$, $\psi\le M$ in $\overline\Om$, and $t<M$, the
sublevel set
$F_t=\{y\in\Om:\psi(y)<t\}$ is compactly contained in $\Om$, so
$\partial^*F_t\subseteq\{\psi=t\}$ and therefore
$\Per(F_t)=\mathcal H^1(\partial^*F_t)\le\ell(t)$, the perimeter being taken
in the sense of sets of finite perimeter \cite{AFP2000}.  The isoperimetric
inequality in the plane, $\Per(F)\ge2\sqrt{\pi|F|}$, gives
\[
\ell(t)\ \ge\ \Per(F_t)\ \ge\ 2\sqrt{\pi A_0(t)} .
\]
Substituting into \eqref{e:csleaves} and integrating over $t\in(m,M)$ by
\eqref{e:coarea} yields \eqref{e:ps}; since the integrand is nonnegative,
restricting the $t$-integral to $(m,M)$ gives a lower bound.

Suppose in addition that $\psi\in C^2(\Om)\cap C^0(\overline\Om)$ and that
every $t\in(m,M)$ is a regular value of $\psi$, as holds for the smooth
reference datum $\psi_0$ under {\rm (H1)}.  If equality holds in \eqref{e:ps},
then equality holds in \eqref{e:csleaves} and in the isoperimetric step for
a.e.\ $t$, hence for every $t\in(m,M)$ because the integrands are continuous
on regular levels: on $\{\psi=t\}$ the gradient has constant magnitude
$a(t)>0$, and $F_t$ is a disk $B(c(t),r(t))$ with
$r(t)=\sqrt{A_0(t)/\pi}$ and
$\nabla\psi=a(t)(\cos\theta,\sin\theta)$ on its boundary.  The disks are
nested, since $F_{t_1}\subset F_{t_2}$ for $t_1<t_2$ gives
\[
|c(t_2)-c(t_1)|\ \le\ r(t_2)-r(t_1).
\]
Because $\psi\in C^2$ with $\nabla\psi\ne0$ on $\{m<\psi<M\}$, the coarea
formula gives
\[
\frac{\dd}{\dd t}\int_{F_t}y\,\dd y
=\int_{\{\psi=t\}}\frac{y}{|\nabla\psi|}\,\dd\mathcal H^1 ,
\]
whose right-hand side is continuous in $t$.  Since $r(t)^2=A_0(t)/\pi$ is
$C^1$ with $r>0$ on $(m,M)$, both $r$ and the center of mass
$c(t)=(\pi r(t)^2)^{-1}\int_{F_t}y\,\dd y$ of the disk $F_t$ are $C^1$ in
$t$, with $2r(t)r'(t)=A_0'(t)/\pi>0$.  Differentiating
$\psi\bigl(c(t)+r(t)\,e^{i\theta}\bigr)=t$ gives
$a(t)\bigl(c'(t)\cdot e^{i\theta}+r'(t)\bigr)=1$ for every
$\theta$, hence $c'(t)=0$: the disks are concentric and $\psi$ is radial.
\end{proof}

In dimension two the area-preserving Morse lemma --- the \emph{Morse--Darboux}
\emph{lemma} --- states that a nondegenerate critical point of a function admits
area-preserving coordinates in which the function is a function of $|w|^2/2$,
and the normalizing map can be taken $C^\infty$ across the critical point
\cite{ColinDeVerdiereVey1979,Eliasson1990}.

\begin{lemma}[Area-preserving Morse coordinates]\label{l:morseiso}
Assume {\rm (H1)}, let $p_0$ be the critical point and $m=\psi_0(p_0)$.
There exist neighborhoods $U$ of $p_0$ and $\mathcal V$ of
$0\in\R^2$, a smooth diffeomorphism $\chi:U\to\mathcal V$ with
$\det D\chi\equiv1$, and a smooth strictly increasing $\Phi$ with $\Phi(0)=m$
and $\Phi'(0)>0$ such that, writing $I=|w|^2/2$ for $w\in\mathcal V$,
\[
\psi_0\circ\chi^{-1}(w)=\Phi(I).
\]
Thus $\chi$ sends each level curve $\{\psi_0=t\}\cap U$ onto the circle
$\{|w|=\sqrt{2I(t)}\}$, where $I(t)=\Phi^{-1}(t)$.  Parametrize a level curve
$\{\psi_0=t\}$ by its \emph{normalized area parameter} $\zeta\in[0,1)$: fixing
a smooth reference point $\sigma_t(0)$ on it, let $\sigma_t(\zeta)$ be the
point for which the arc from $\sigma_t(0)$ to $\sigma_t(\zeta)$ along
$\{\psi_0=t\}$ satisfies
\[
\int_{\text{arc}}\frac{\dd\mathcal H^1}{|\nabla\psi_0|}=\zeta\,A_0'(t) .
\]
\end{lemma}

If the reference point is chosen on $\chi^{-1}(\{w_2=0,\ w_1>0\})$, then
$\zeta=(\arg w)/(2\pi)$ $(\operatorname{mod}1)$.

\begin{proof}
\emph{1. The area-preserving normal form.}  The existence of $\chi$ and
$\Phi$ is the area-preserving Morse lemma for smooth symplectic
surfaces \cite{ColinDeVerdiereVey1979,Eliasson1990}: near a nondegenerate critical point
of a smooth function on a surface there are $C^\infty$ area-preserving
coordinates in which the function is a function of $|w|^2/2$, the normalizing
map being smooth across the critical point.  (The analytic precursor
\cite{Vey1977} treats holomorphic germs; it is the smooth-category result of
\cite{ColinDeVerdiereVey1979} that applies here.)  The
mechanism is the Birkhoff normal form of a one-degree-of-freedom Hamiltonian
system, where no small divisors arise.  We thus have $U$, $\mathcal V$, $\chi$ and $\Phi$ as in the
statement, with
\[
\psi_0\circ\chi^{-1}(w)=\Phi(I)\qquad(w\in\mathcal V).
\]
Since $\Phi$ is strictly increasing, $\chi$ sends each level curve
$\{\psi_0=t\}\cap U$ onto the circle $\{|w|=\sqrt{2I(t)}\}$, where
$I(t)=\Phi^{-1}(t)$.

\emph{2. The flow angle is the polar angle.}  Let $X=\nabla^\perp\psi_0$.
It is divergence free, so its flow is area preserving, and $\psi_0$ is
constant along it, so $X$ is tangent to the level curves of $\psi_0$; near
the nondegenerate minimum $p_0$ these are smooth embedded circles shrinking
to $p_0$, and $X$ is periodic on each of them.  Because $\chi$ is area
preserving and carries $\psi_0$ to $\Phi(|w|^2/2)$, the push-forward
$\chi_*X$ is the Hamiltonian field of $\Phi(|w|^2/2)$ with respect to the
standard area form, namely $\Phi'(I)\,(-w_2,w_1)$; its flow is the rigid
rotation about $0$ with angular speed $\Phi'(I)$.  Hence the normalized
flow-time angle $\theta$ of $X$, measured from a reference point on
$\chi^{-1}(\{w_2=0,\ w_1>0\})$, is exactly the polar angle $\arg w$ of the
image point.

\emph{3. Smoothness of the area function at $m$.}  By the ordinary Morse
lemma choose coordinates $\xi$ near $p_0$ in which
$\psi_0=m+\frac12\xi^TQ\xi$ \emph{exactly}, with $Q>0$.  Put $z=Q^{1/2}\xi$,
so that $\psi_0=m+\frac12|z|^2$; in these coordinates the Euclidean area
element $\dd x_1\,\dd x_2$ becomes $a(z)\,\dd z_1\,\dd z_2$, where
$a(z)=\bigl|\det(\partial x/\partial z)\bigr|$ is the Jacobian density of the
chart $z\mapsto x$, smooth and strictly positive.  For $t$ near $m$ the
sublevel set is the disk $\{|z|<r(t)\}$ with $r(t)=\sqrt{2(t-m)}$, so
\[
A_0(t)=\int_{|z|<r(t)}a(z)\,\dd z
=\int_0^{r(t)}s\,a_{\rm av}(s)\,\dd s ,
\qquad
a_{\rm av}(s)=\int_0^{2\pi}a(s\cos\alpha,s\sin\alpha)\,\dd\alpha .
\]
The function $a_{\rm av}$ is smooth and even, hence a smooth function of
$s^2$; writing $a_{\rm av}(s)=\widetilde a(s^2)$ and substituting $u=s^2$
gives
\[
A_0(t)=\frac12\int_0^{2(t-m)}\widetilde a(u)\,\dd u ,
\]
which is smooth in $t$ near $t=m$, with $A_0(m)=0$ and
$A_0'(m)=a_{\rm av}(0)=2\pi\,a(0)>0$.  Since $\chi$ is area preserving, the
sublevel set $\{\psi_0<t\}$ and its image $\{|w|<\sqrt{2I(t)}\}$ have equal
area, so $A_0(t)=2\pi I(t)$, that is $I=A_0(\psi_0)/(2\pi)$; this is
consistent with part~1 and gives $\Phi'(0)=2\pi/A_0'(m)>0$.

\emph{4. The normalized area parameter.}  The pair $(I,\theta)$ is a
system of action-angle coordinates, so the area element is
$\dd I\,\dd\theta$; write the metric as
$E\,\dd I^2+2F\,\dd I\dd\theta+G\,\dd\theta^2$.  Then $EG-F^2=1$, so the
inverse metric satisfies $g^{II}=G/(EG-F^2)=G$.  As $\psi_0=\Phi(I)$ depends
only on $I$,
\[
|\nabla\psi_0|^2=g^{II}\Phi'(I)^2=G\,\Phi'(I)^2,\qquad
\dd\mathcal H^1=\sqrt G\,\dd\theta\ \ \text{on a level curve},
\]
so $\dd\mathcal H^1/|\nabla\psi_0|=\dd\theta/\Phi'(I)$.  Integrating
over a full level curve gives $A_0'(t)=2\pi/\Phi'(I)$.  Comparing with the
defining relation
$\int_{\text{arc}}\dd\mathcal H^1/|\nabla\psi_0|=\zeta\,A_0'(t)$ for the arc
from $\sigma_t(0)$ to $\sigma_t(\zeta)$, and writing $\Delta\theta$ for the
change of $\theta$ along that arc,
\[
\int_{\text{arc}}\frac{\dd\mathcal H^1}{|\nabla\psi_0|}
=\frac{\Delta\theta}{\Phi'(I)}
=\zeta\,\frac{2\pi}{\Phi'(I)},
\]
so $\zeta=\Delta\theta/(2\pi)$.  By part~2, $\Delta\theta$ is the change of
$\arg w$ along the arc, and the reference point lies on
$\{w_2=0,\ w_1>0\}$; hence $\zeta=(\arg w)/(2\pi)$ $(\operatorname{mod}1)$.
\end{proof}

\begin{theorem}[Sharp rearrangement bound and attainment]\label{t:rearr}
Assume {\rm (H1)}, let $R=\sqrt{V/\pi}$, let $\mathbb D_R$ be the disk of
radius $R$ centered at the origin, and define the \emph{circular
rearrangement} of the stream function by
\begin{equation}\label{e:psistar}
\varrho(t)=\sqrt{A_0(t)/\pi}\ \ (t\in[m,M]),\qquad
\psi_*(y)=\varrho^{-1}(|y|)\ \ (y\in\overline{\mathbb D_R}) ,
\end{equation}
where $\varrho^{-1}$ denotes the inverse of the increasing bijection
$\varrho:[m,M]\to[0,R]$; put $\Bv_*=\nabla^\perp\psi_*$.  Then:
\begin{enumerate}[\rm (i)]
\item For every $\eta\in\widehat{\mathcal G}_2$,
\begin{equation}\label{e:Emin}
E(\eta)\ \ge\ \underbrace{\tfrac12\int_m^M\frac{4\pi A_0(t)}{A_0'(t)}
\,\dd t}_{=:\,M_*}\ +\ \kap\,2\sqrt{\pi V} ,
\end{equation}
and equality is attained by an explicit smooth $\eta_*\in\mathcal G_2$ with
$\eta_*(\Om_0)=\mathbb D_R$ and $(\eta_*)_*\Bv_0=\Bv_*$.  Since
$\mathcal G_2\subset\widehat{\mathcal G}_2$, the infimum over the original
smooth class $\mathcal G_2$ equals the right-hand side of \eqref{e:Emin} and
is attained.
\item $\psi_*$ is radial and smooth, $\Bv_*=\rho(|y|)\,e_\theta$ with
$\rho=\psi_*'$, and the shape term is the isoperimetric minimum because
$\eta_*(\Om_0)=\mathbb D_R$ is a disk.
\item $(\mathbb D_R,\Bv_*,\pv_*)$ is a classical solution of
\eqref{e:eq} and \eqref{e:bc2d}, where $\pv_*$ is the radial function fixed by
\begin{equation}\label{e:pstar}
\pv_*'(r)=-\Bigl(\rho(r)\rho'(r)+\frac{\rho(r)^2}{r}\Bigr),\qquad
\pv_*(R)=\frac{\kap}{R}-\tfrac12\rho(R)^2 .
\end{equation}
\end{enumerate}
\end{theorem}

\begin{proof}
\emph{1. Construction of $\eta_*$.}  For $t\in(m,M)$ let
$C_0(t)=\{\psi_0=t\}$, a single circle by \eqref{e:subdisc}, and recall from
\eqref{e:A0prime} that
$\oint_{C_0(t)}|\nabla\psi_0|^{-1}\dd\mathcal H^1=A_0'(t)$.  We need a
reference point $\sigma_t(0)\in C_0(t)$ depending smoothly on $t$ which,
near $p_0$, lies on the curve $\chi^{-1}(\{w_2=0,\ w_1>0\})$ of
Lemma~\ref{l:morseiso}.  Let $U$ and $\chi$ be as in that lemma and let
$\lambda\in C^\infty(U)$ be $1$ near $p_0$ and $0$ near $\partial U$.  In
the $\chi$ coordinates the radial field $g(|w|^2)w$, with $g>0$ smooth, is
transverse to every circle $|w|=\sqrt{2\Phi^{-1}(t)}$ and its trajectories are
the rays;
denote by $Z_{\rm rad}$ the field on $U$ that equals $g(|w|^2)w$ in these
coordinates.  The field
$Z_{\rm grad}=\nabla\psi_0/|\nabla\psi_0|^2$ is transverse to the same
foliation on all of $\Om_0\setminus\{p_0\}$.  Transversality, written
$\dd\psi_0(Z)>0$, is preserved under convex combination, so the field
$Z:=\lambda Z_{\rm rad}+(1-\lambda)Z_{\rm grad}$ on $U$ and
$Z:=Z_{\rm grad}$ on $\Om_0\setminus U$ is smooth on $\overline\Om_0$ and
satisfies $\dd\psi_0(Z)>0$ on $\Om_0\setminus\{p_0\}$.  Choose
$p_1\in\chi^{-1}(\{w_2=0,\ w_1>0\})$ so close to $p_0$ that
$\{\psi_0\le\psi_0(p_1)\}\subset\{\lambda=1\}$, and let $\sigma_t(0)$ be
the point of $C_0(t)$ on the trajectory of $Z$ through $p_1$, followed
forward for $t>\psi_0(p_1)$ and backward for $t<\psi_0(p_1)$.  Along the
trajectory $\psi_0$ increases strictly in forward time and decreases
strictly in backward time, and neither end can converge to an intermediate
level, since $\dd\psi_0(Z)$ is bounded below on every compact set on which
$\psi_0$ is bounded away from $m$; hence the trajectory runs from $p_0$ to
$\partial\Om_0$ and meets each level circle exactly once, so $\sigma_t(0)$
is defined on all of $(m,M)$.  It lies on the positive
$w_1$-axis for $t\le\psi_0(p_1)$, since there the trajectory is confined to
$\{\lambda=1\}$, where $Z=Z_{\rm rad}$; and $t\mapsto\sigma_t(0)$ is
smooth on $(m,M)$, because for $t\le\psi_0(p_1)$ it is
$\chi^{-1}(\sqrt{2\Phi^{-1}(t)},0)$, smooth as $\Phi'(I)>0$, while for
$t\ge\psi_0(p_1)$ smooth dependence of the flow applies, the crossing time
of the level $t$ being smooth in $t$ because $\dd\psi_0(Z)>0$.

Parametrize
$C_0(t)$ by the normalized area parameter $\zeta\in[0,1)$ of
Lemma~\ref{l:morseiso}.
Set
\begin{equation}\label{e:etastar}
\eta_*\bigl(\sigma_t(\zeta)\bigr)=\varrho(t)\,(\cos2\pi\zeta,\ \sin2\pi\zeta),
\qquad \eta_*(p_0)=0 .
\end{equation}

\emph{2. $\eta_*$ is a homeomorphism.}  It maps $C_0(t)$ bijectively onto
the circle of radius $\varrho(t)$ centered at the origin.  The level curves
$\{C_0(t)\}_{t\in(m,M)}$ foliate $\Om_0\setminus\{p_0\}$ and the circles
foliate $\overline{\mathbb D_R}\setminus\{0\}$; since $\varrho(t)\downarrow0$
as $t\downarrow m$ and $\varrho(t)\uparrow R$ as $t\uparrow M$, the definition
\eqref{e:etastar} extends continuously by $\eta_*(p_0)=0$ and by sending
$\partial\Om_0=\{\psi_0=M\}$ onto the circle of radius $R$.  Hence $\eta_*$
is a continuous bijection of the compact set $\overline\Om_0$ onto
$\overline{\mathbb D_R}$, therefore a homeomorphism.

\emph{3. $\eta_*$ is area-preserving.}  Compute the area element in the
coordinates $(t,\zeta)$.  On the source side, the strip between $C_0(t)$ and
$C_0(t+dt)$ has, along the arc $[\zeta,\zeta+\dd\zeta]$, area
$\frac{\dd t}{|\nabla\psi_0|}\,\dd\mathcal H^1$; by the definition of $\zeta$
the measure $|\nabla\psi_0|^{-1}\dd\mathcal H^1$ of that arc equals
$A_0'(t)\,\dd\zeta$.  Hence the area element of the source, in the
coordinates $(t,\zeta)$, is $A_0'(t)\,\dd t\,\dd\zeta$.  On the target side,
writing $y=\varrho(t)(\cos2\pi\zeta,\sin2\pi\zeta)$ and using
$\varrho^2=A_0/\pi$, the area element is
\[
s\,\dd s\,\dd\theta
=\varrho(t)\varrho'(t)\,\dd t\cdot2\pi\,\dd\zeta
=\pi\,(\varrho^2)'(t)\,\dd t\,\dd\zeta
=A_0'(t)\,\dd t\,\dd\zeta ,
\]
where $s=|y|$.  The two coincide, so $|\det D\eta_*|=1$ a.e.; and $\eta_*$ is
orientation preserving, because in the coordinates $(t,\zeta)$ used above on
the source and in the polar coordinates $(s,\theta)$ of the target it reads
$(t,\zeta)\mapsto(\varrho(t),2\pi\zeta)$, with $\varrho$ increasing and both
coordinate systems positively oriented.  Hence
$\det D\eta_*=1$ a.e.\ and $\eta_*$ is area-preserving.

\emph{4. Regularity: $\eta_*$ is smooth and belongs to $\mathcal G_2$.}

(a) \emph{Off the critical point.}  On
$\overline\Om_0\setminus\{p_0\}$ the vector field $X=\nabla^\perp\psi_0$ is
smooth and nonvanishing, and $|X|=|\nabla\psi_0|$, so that
$\dd\tau=\dd\mathcal H^1/|\nabla\psi_0|$ is the flow time of $X$.  Thus
$\zeta=\tau/A_0'(t)$, where $\tau$ is the time needed to flow from the
reference point to the given point and $A_0'(t)$ is the period, smooth by
\eqref{e:A0prime}.  Let $\Xi(t,s)$ be the point reached by following the flow
of $X$ from $\sigma_t(0)$ for time $s$.  Then $\Xi$ is smooth in both
variables and its differential is everywhere nonsingular: the image of
$\partial_t$ is transverse to the level curves, since the section
$t\mapsto\sigma_t(0)$ moves across them, while the image of $\partial_s$ is
$X$, which is tangent to them.  Moreover $\psi_0\circ\Xi(t,s)=t$, and the flow
of $X$ on $C_0(t)$ has minimal period $A_0'(t)$, so
$\Xi(t,s)=\Xi(t',s')$ exactly when $t=t'$ and
$s-s'\in A_0'(t)\mathbb Z$.  The map $(t,s)\mapsto(t,s+A_0'(t))$ is a
smooth free diffeomorphism whose iterates act properly, since $A_0'$ is
bounded below on every compact subinterval of $(m,M)$; hence $\Xi$ is a
smooth covering map from $(m,M)\times\R$ onto $\Om_0\setminus\{p_0\}$, with
this map as deck transformation.  On the cover $(\psi_0,\tau)$ are smooth
coordinates and $\zeta=\tau/A_0'(t)$ is smooth.  Under the deck
transformation $\zeta$ increases by $1$, so the pair
$(\cos2\pi\zeta,\sin2\pi\zeta)$ is invariant and descends to a smooth
$\R^2$-valued function on $\Om_0\setminus\{p_0\}$.  Since
$\nabla\psi_0\ne0$ on $\partial\Om_0$ by {\rm (H1)}, the level curves converge
smoothly to $\partial\Om_0$ as $t\uparrow M$, so the section $\sigma_t(0)$ and
the period $A_0'(t)$ extend smoothly to $t=M$ and the descent is smooth on
$\overline\Om_0\setminus\{p_0\}$.  Since $t=\psi_0$ and
$\varrho$ are smooth with $\varrho>0$ there, the map
$\eta_*=\varrho(t)(\cos2\pi\zeta,\sin2\pi\zeta)$ is smooth on
$\overline\Om_0\setminus\{p_0\}$.

(b) \emph{At the critical point.}  With the reference point chosen as in
step~1, let $U$ and $\chi:U\to\mathcal V$ be as in Lemma~\ref{l:morseiso} and
shrink $U$, if necessary, to a disk about $p_0$ so small that the sublevel set
$\{\psi_0<t\}$ is contained in $U$ whenever $t=\psi_0(x)$ with $x\in U$ and,
in addition, $\psi_0<\psi_0(p_1)$ throughout $U$ (both are possible because
$\psi_0(p_0)=m<\psi_0(p_1)$); by step~1 the reference point $\sigma_t(0)$ then
lies on $\chi^{-1}(\{w_2=0,\ w_1>0\})$ at every level met in $U$.

For $x\in U\setminus\{p_0\}$ put $w=\chi(x)$ and $t=\psi_0(x)$.  The map
$\chi$ sends the sublevel set $\{\psi_0<t\}\subset U$ onto the disk
$\{w'\in\mathcal V:|w'|<|w|\}$: indeed
$\psi_0\circ\chi^{-1}(w')=\Phi(|w'|^2/2)$ is increasing in $|w'|$ and equals
$\Phi(|w|^2/2)=\psi_0(x)=t$ on $|w'|=|w|$.  Since $\chi$ is area preserving
this disk has the same area as the sublevel set, namely $A_0(t)$; as it is a
Euclidean disk of radius $|w|$, we get $A_0(t)=\pi|w|^2$ and therefore
$\varrho(t)=\sqrt{A_0(t)/\pi}=|w|$.  By the final
assertion of Lemma~\ref{l:morseiso}, with the above choice of reference
point $\zeta=(\arg w)/(2\pi)$.  Consequently
\[
\eta_*(x)=\varrho(t)\,(\cos2\pi\zeta,\ \sin2\pi\zeta)
=|w|\,\bigl(\cos(\arg w),\ \sin(\arg w)\bigr)=w=\chi(x)
\]
for $x\in U\setminus\{p_0\}$; by continuity the identity $\eta_*=\chi$ holds
on all of $U$, both sides sending $p_0$ to $0$.  Since $\chi$ is a smooth
diffeomorphism, $\eta_*$ is smooth at $p_0$.

(c) \emph{Conclusion.}  Thus $\eta_*\in C^\infty(\overline\Om_0;\R^2)$.
Its Jacobian determinant is $1$ off $p_0$ by step 3 and equals
$\det D\chi=1$ on $U$, so $\det D\eta_*\equiv1$.  A homeomorphism (step 2)
that is a local diffeomorphism is a diffeomorphism onto its image; hence
$\eta_*$ is a smooth diffeomorphism of $\overline\Om_0$ onto
$\overline{\mathbb D_R}$ with $\det D\eta_*=1$, i.e.\
$\eta_*\in\mathcal G_2\subset\widehat{\mathcal G}_2$.

(d) \emph{The profile is smooth.}  By the smoothness of $A_0$ at $m$
established in Lemma~\ref{l:morseiso}, with $A_0'(m)>0$, and since
$\nabla\psi_0\ne0$ on $\partial\Om_0$ gives $A_0'(M)>0$, here $A_0$ is
smooth up to $M$ as well, since the level curves $\{\psi_0=t\}$ converge
smoothly to $\partial\Om_0$ as $t\uparrow M$, inversion of
$\varrho^2=A_0(t)/\pi$ shows that $\psi_*$ is a smooth function of $|y|^2$
near $0$ and smooth up to $|y|=R$.  Hence
$\Bv_*=\nabla^\perp\psi_*$ is smooth on $\overline{\mathbb D_R}$.

\emph{5. $\eta_*$ attains the bound.}  By \eqref{e:etastar} and the
definition of $\psi_*$ one has $\psi_*\circ\eta_*=\psi_0$ on $\Om_0$: both
sides equal $t$ on $C_0(t)$ and equal $m$ at $p_0$.  Hence
Lemma~\ref{l:conj}(i) gives
$(\eta_*)_*\Bv_0=\nabla^\perp\psi_*=\Bv_*$.  The sublevel sets of $\psi_*$
are concentric disks, so the isoperimetric step in the proof of
Lemma~\ref{l:ps} is an equality; and $|\nabla\psi_*|$ is constant on each
level circle, so the Cauchy--Schwarz step is an equality as well.  Hence
\[
\int_{\mathbb D_R}|\nabla\psi_*|^2\,\dd y
=\int_m^M\frac{4\pi A_0(t)}{A_0'(t)}\,\dd t ,
\]
and by Lemma~\ref{l:conj}(iii) the magnetic energy of $\eta_*$ is
$\tfrac12\int_{\mathbb D_R}|\nabla\psi_*|^2=M_*$.  Since
$\eta_*(\Om_0)=\mathbb D_R$ has perimeter $2\pi R=2\sqrt{\pi V}$, we get
$E(\eta_*)=M_*+\kap\,2\sqrt{\pi V}$, which is the right-hand side of
\eqref{e:Emin}.

\emph{6. The bound holds for every $\eta$.}  Let
$\eta\in\widehat{\mathcal G}_2$ and $\psi=\psi_0\circ\eta^{-1}$.  By
Lemma~\ref{l:conj}(i), $\psi\in W^{1,\infty}(\Om)\cap C^0(\overline\Om)$
and $\Bv=\nabla^\perp\psi$, so $\psi$ satisfies the hypotheses of
Lemma~\ref{l:ps}.  By
Lemma~\ref{l:conj}(ii) the distribution function of $\psi$ is $A_0$, so
Lemma~\ref{l:ps} and Lemma~\ref{l:conj}(iii) give
$\int_{\Om_0}|D\eta\Bv_0|^2=\int_\Om|\Bv|^2
=\int_\Om|\nabla\psi|^2\ge\int_m^M\frac{4\pi A_0}{A_0'}\dd t$.  For the
shape term, $|\eta(\Om_0)|=V$ by area preservation, so the isoperimetric
inequality gives $\Per(\eta(\Om_0))\ge2\sqrt{\pi V}$.  Adding the two gives
\eqref{e:Emin}.

\emph{7. The minimizer is a classical equilibrium.}  Write
$\Bv_*=\psi_*'(|y|)\,e_\theta=\rho(|y|)e_\theta$.  Since $\psi_*$ is smooth
and radial with $\psi_*'(0)=0$, the function $\rho$ is smooth and odd to
leading order, so $\Bv_*$ is a smooth field on $\overline{\mathbb D_R}$; it
is divergence-free by construction and tangential because $\psi_*$ is
constant on $|y|=R$.  For the force balance use the identity
$(\curl\Bv)\times\Bv=(\Bv\cdot\nabla)\Bv-\nabla(\tfrac12|\Bv|^2)$.  With
$\Bv=\rho e_\theta$ one has
$\nabla(\tfrac12|\Bv|^2)=\rho\rho'\,e_r$ and, since
$\partial_\theta e_\theta=-e_r$,
\[
(\Bv\cdot\nabla)\Bv=\frac{\rho}{r}\partial_\theta(\rho\,e_\theta)
=-\frac{\rho^2}{r}\,e_r .
\]
Hence $(\curl\Bv_*)\times\Bv_*=-(\rho\rho'+\rho^2/r)e_r=\nabla\pv_*$ by
\eqref{e:pstar}.  Near $r=0$ we have $\rho(r)=\rho_1r+O(r^3)$, so
$\rho\rho'+\rho^2/r=O(r)$ and $\pv_*$ is a smooth function of $r^2$, hence a
smooth radial function on the disk.  Finally, on $|y|=R$ the curvature is
$\kappa_{\mathrm{curv}}=1/R$ and
$\pv_*(R)+\tfrac12\rho(R)^2=\kap/R=\kap\,\kappa_{\mathrm{curv}}$ by the
normalization in \eqref{e:pstar}, which is exactly \eqref{e:bc2d}.  Thus
$(\mathbb D_R,\Bv_*,\pv_*)$ satisfies \eqref{e:eq} and \eqref{e:bc2d}
classically.
\end{proof}

\begin{corollary}[The rearrangement is the normal form near the axis]%
\label{c:axis}
In the notation of Lemma~\ref{l:morseiso} and Theorem~\ref{t:rearr}, and with
the reference point chosen as in step~1 of the proof there, one has
$\eta_*=\chi$ on $U$.
\end{corollary}

\begin{proof}
This is step~4(b) of the proof of Theorem~\ref{t:rearr}: on
$U\setminus\{p_0\}$ one has $\varrho(t)=|w|$ and $\zeta=(\arg w)/(2\pi)$,
hence $\eta_*=\chi$, and both sides send $p_0$ to $0$.
\end{proof}

Consequently the smoothness of the rearrangement across the magnetic axis is
exactly the smoothness of the area-preserving normal form, and it does not
depend on the construction of the reference section $\sigma_t(0)$.

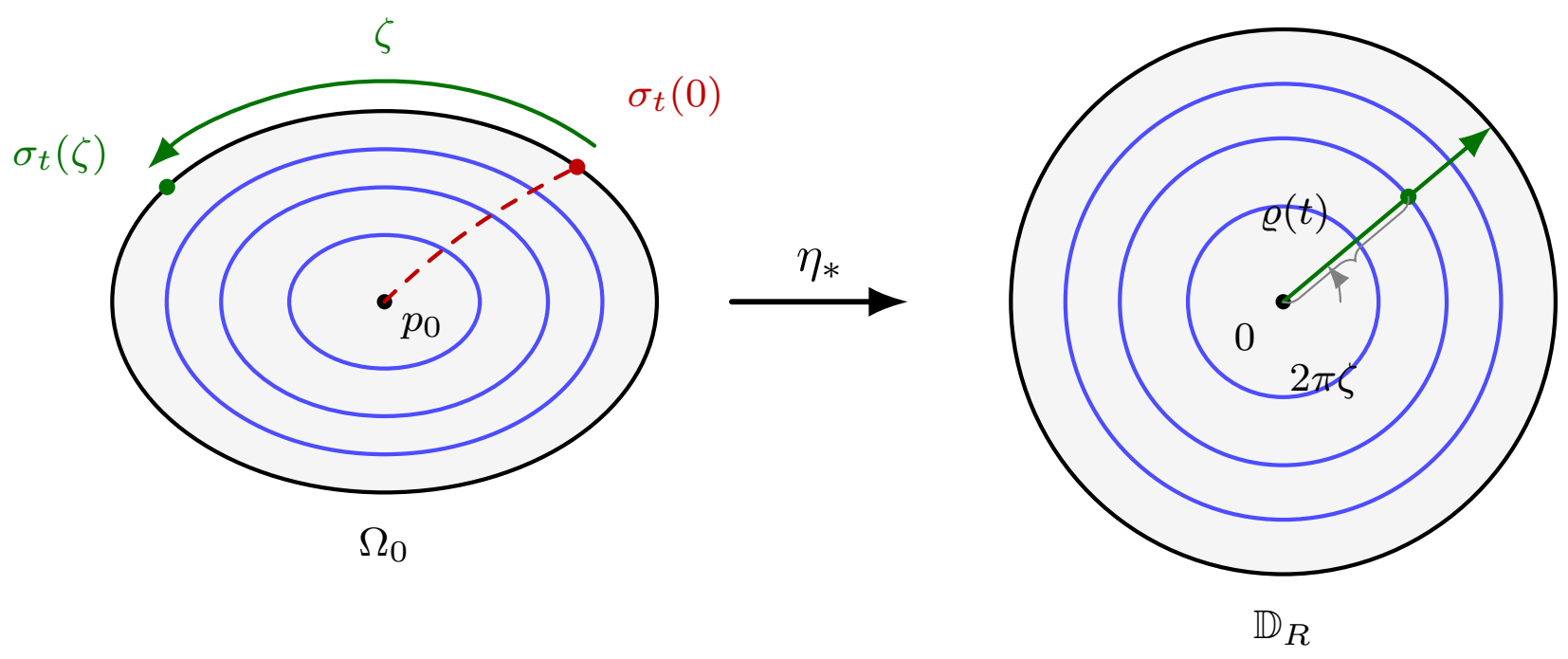
\begin{figure}[!ht]
\centering
\begin{tikzpicture}[scale=0.95,>=Latex,line cap=round]
  \begin{scope}
    \draw[fill=black!4,thick] (0,0) ellipse (2.00 and 1.40);
    \foreach \s in {0.35,0.60,0.80}{%
      \draw[blue!70,thick] (0,0) ellipse ({2.00*\s} and {1.40*\s});}
    \fill (0,0) circle (1.6pt);
    \node[anchor=north west,inner sep=1pt] at (0.08,-0.04) {\scriptsize $p_0$};
    \draw[red!75!black,dashed,thick] (0,0) .. controls (0.55,0.55) and (1.05,0.80) .. ({2.00*cos(45)},{1.40*sin(45)});
    \fill[red!75!black] ({2.00*cos(45)},{1.40*sin(45)}) circle (1.7pt);
    \node[red!75!black,anchor=south west,inner sep=1pt] at (1.74,1.32) {\scriptsize $\sigma_t(0)$};
    \fill[green!45!black] ({2.00*cos(143)},{1.40*sin(143)}) circle (1.7pt);
    \node[green!45!black,anchor=east,inner sep=1pt] at (-1.98,1.08) {\scriptsize $\sigma_t(\zeta)$};
    \draw[green!45!black,thick,->] ({2.18*cos(45)},{1.62*sin(45)}) arc (45:143:2.18 and 1.62);
    \node[green!45!black,anchor=south,inner sep=1pt] at (0,1.78) {\scriptsize $\zeta$};
    \node[anchor=north] at (0,-1.52) {\scriptsize $\Om_0$};
  \end{scope}
  \draw[->,line width=1.1pt] (2.55,0)--(3.85,0) node[midway,above]{\small$\eta_*$};
  \begin{scope}[shift={(6.60,0)}]
    \draw[fill=black!4,thick] (0,0) circle (2.00);
    \foreach \s in {0.35,0.60,0.80}{%
      \draw[blue!70,thick] (0,0) circle ({2.00*\s});}
    \fill (0,0) circle (1.6pt);
    \node[anchor=north east,inner sep=1pt] at (-0.16,-0.12) {\scriptsize $0$};
    \draw[green!45!black,thick,->] (0,0)--(40:2.00);
    \fill[green!45!black] (40:1.20) circle (1.7pt);
    \draw[decorate,decoration={brace,amplitude=3pt,mirror},gray]
      (0,0) -- (40:1.20) node[midway,above left=4pt,inner sep=0pt,black]
      {\scriptsize $\varrho(t)$};
    \draw[->,gray] (0.42,0) arc (0:40:0.42);
    \node[black,anchor=north,inner sep=1pt] at (0.30,-0.40) {\scriptsize $2\pi\zeta$};
    \node[anchor=north] at (0,-2.12) {\scriptsize $\mathbb D_R$};
  \end{scope}
\end{tikzpicture}
\caption{The circular rearrangement $\eta_*$ of Theorem~\ref{t:rearr}.  The
level curve $C_0(t)=\{\psi_0=t\}$ of the reference stream function (left) is
sent onto the circle of radius $\varrho(t)=\sqrt{A_0(t)/\pi}$ (right); the
normalized area parameter $\zeta$, the fraction of
$\int|\nabla\psi_0|^{-1}\dd\mathcal H^1$ along the level, is preserved, so
that $\zeta$ becomes the polar angle divided by $2\pi$.  The radius is fixed
by the enclosed area, so the map is area preserving, and near $p_0$ it
coincides with the area-preserving normal form of Lemma~\ref{l:morseiso}
(Corollary~\ref{c:axis}).  The angular positions marked in the two panels
are schematic only and are not to be read as numerically equal.}
\label{f:rearr}
\end{figure}

The map is pictured in Figure~\ref{f:rearr}: level curves are sent onto
concentric circles, the area parameter $\zeta$ is preserved, and the radius
$\varrho(t)$ is fixed by the enclosed area.

\begin{proposition}[Taylor sign condition]\label{p:taylor}
In the situation of Theorem~\ref{t:rearr}, let $\Pv_*=\pv_*+\tfrac12|\Bv_*|^2$
be the total pressure of the constructed equilibrium, $\rho=\psi_*'$ and
$R=\sqrt{V/\pi}$.  Then
\begin{equation}\label{e:taylor}
-\partial_\nu\Pv_*\big|_{\partial\mathbb D_R}=-\,\Pv_*'(R)
=\frac{\rho(R)^2}{R},
\end{equation}
so the Taylor sign condition $-\partial_\nu\Pv_*>0$ holds if and only if
$\rho(R)\ne0$, that is, if and only if $|\Bv_*|=\rho(R)$ does not vanish on the
free boundary, and its strength is then exactly $\rho(R)^2/R$.  Under {\rm (H1)}
this alternative is never degenerate: $\rho(R)=2\pi R/A_0'(M)>0$.
\end{proposition}

\begin{proof}
By \eqref{e:pstar},
\[
\Pv_*'(r)=\pv_*'(r)+\rho(r)\rho'(r)
=-\Bigl(\rho\rho'+\frac{\rho^2}{r}\Bigr)+\rho\rho'=-\frac{\rho(r)^2}{r} .
\]
Since $\Pv_*$ is radial and $\nv=e_r$ on $\partial\mathbb D_R$, this gives
\eqref{e:taylor}.  Next, $\psi_*(\varrho(t))=t$ and
$\varrho'(t)=A_0'(t)/(2\pi\varrho(t))$ give
$\rho(\varrho(t))=1/\varrho'(t)=2\pi\varrho(t)/A_0'(t)$; at $t=M$, where
$\varrho(M)=R$, this is $\rho(R)=2\pi R/A_0'(M)$.  Under {\rm (H1)} one has
$\nabla\psi_0\ne0$ on $\partial\Om_0$, so $A_0'(M)<\infty$ and $\rho(R)>0$.
\end{proof}

\begin{remark}[The second variation is flat along twists]\label{r:hessian}
By Theorem~\ref{t:disk} every smooth twist $\eta_\alpha$ of the disk is a
global minimizer; here $\alpha(r)=\alpha_0+h(r^2)$ with
$h\in C^\infty([0,1])$, the condition for smoothness at the origin, and such
$\alpha$ are called admissible.  Hence, for every admissible $\alpha$ and
every $s\in\R$,
$E(\eta_{s\alpha})=E(\mathrm{id})$, so the second derivative of
$s\mapsto E(\eta_{s\alpha})$ at $s=0$ vanishes:
\[
\delta^2E(\mathrm{id})[\alpha,\alpha]=0\qquad\text{for every admissible }\alpha.
\]
The Hessian therefore vanishes on the whole twist direction, so the
degeneracy is a flat direction of the second variation.
\end{remark}

\begin{example}[A non-circular, non-rotationally symmetric datum]%
\label{ex:ellipse}
Take $a>b>0$, $\Om_0=\{x:x_1^2/a^2+x_2^2/b^2<1\}$ and
$\psi_0=x_1^2/a^2+x_2^2/b^2-1$, so that $\Bv_0=\nabla^\perp\psi_0$
vanishes only at the origin and $\psi_0=0$ on $\partial\Om_0$: hypothesis
{\rm (H1)} holds, with $m=-1$, $M=0$ and $V=\pi ab$.  The sublevel sets are
ellipses with $A_0(t)=\pi ab(1+t)$, $A_0'=\pi ab$, so
\[
M_*=\tfrac12\int_{-1}^0\frac{4\pi A_0(t)}{A_0'(t)}\,\dd t
=\tfrac12\int_{-1}^0 4\pi(1+t)\,\dd t=\pi ,
\]
while a direct computation gives
\[
\tfrac12\int_{\Om_0}|\Bv_0|^2\,\dd x
=\frac{\pi(a^2+b^2)}{2ab}=\pi+\frac{\pi(a-b)^2}{2ab}\ \ge\ \pi ,
\]
with equality only for $a=b$.  Since $\Per(\Om_0)>2\pi\sqrt{ab}=2\sqrt{\pi V}$
for $a\ne b$, the identity is strictly worse than the rearrangement:
\[
E(\mathrm{id})-E(\eta_*)=\frac{\pi(a-b)^2}{2ab}
+\kap\bigl(\Per(\Om_0)-2\pi\sqrt{ab}\bigr)>0\qquad(a\ne b).
\]
The rearranged configuration is the rigid rotation $\Bv_*=\frac{2r}{ab}
e_\theta$ with $r=|y|$, on the disk of radius $\sqrt{ab}$, with
$\Delta\psi_*=4/(ab)$, a constant; the whole computation is carried out in
Appendix~\ref{app:ellipse}.  This is the simplest instance in which
Theorem~\ref{t:rearr} is verified on data that are neither circular nor
rotationally symmetric.
\end{example}

\begin{remark}[Fields with several critical points]\label{r:multicrit}
If $\psi_0$ has more than one critical point, the lower bound \eqref{e:ps}
remains valid: its proof uses only the Cauchy--Schwarz inequality, the isoperimetric
inequality, and the boundary condition $\psi=M$ on $\partial\Om$ (satisfied
whenever the maximum of $\psi_0$ is attained on $\partial\Om_0$, in particular
under {\rm (H2)} below), none of which depends on the number of critical
points.  It is no longer sharp: when $\{\psi_0<t\}$
has $k\ge2$ components of areas $a_1,\dots,a_k$ with $\sum a_i=A_0(t)$, the
image sublevel set also has $k$ components (a homeomorphism preserves the
number of components), so the sharp isoperimetric bound at that level is
$\ell(t)\ge\sum_i2\sqrt{\pi a_i}$, which is strictly larger than
$2\sqrt{\pi A_0(t)}$ by strict concavity of $\sqrt{\ \cdot\ }$.  The
level-wise bound is obtained in Section~\ref{ss:multi} by tracking the component
structure exactly (the number of components and their areas are homeomorphism
invariants) and applying the isoperimetric inequality to each component
separately (Theorem~\ref{t:multi}); its nonattainment at an isolated simple
merge is Proposition~\ref{p:notatt}.
\end{remark}

\subsection{The multi-critical case: a sharper bound}\label{ss:multi}

As anticipated in Remark~\ref{r:multicrit}, in this subsection we keep the
standing assumptions of Section~\ref{ss:rearr}, in particular the simple
connectivity of $\Om_0$, which is used below in
Lemma~\ref{l:invar}; we drop {\rm (H1)} and assume
\begin{enumerate}
\item[(H2)] all interior critical points of $\psi_0$ are nondegenerate,
$\psi_0$ has no interior local maximum, and the boundary value $c_0=M$ is
regular in the sense that $\nabla\psi_0\ne0$ on $\partial\Om_0$.
\end{enumerate}
The interior critical points are then finite in number and the usual Morse
coordinates are available at them.  The boundary restriction of $\psi_0$ is
constant and therefore everywhere degenerate as a function on $\partial\Om_0$;
the condition $\nabla\psi_0\ne0$ on $\partial\Om_0$ is what keeps the boundary
value out of the critical-value set, so that logarithmic behavior
$A_0'(t)\sim\log(1/|t-s|)$ occurs only at interior saddle values and the
levels used below lie in $(m,M)$.
A nondegenerate interior critical point is either \emph{elliptic} (the
Hessian of $\psi_0$ is definite, so the point is a magnetic axis about which
the field lines close) or \emph{hyperbolic}, a saddle.

Interior local maxima of $\psi_0$ are elliptic critical points of the opposite orientation (magnetic axes with reversed circulation);
they are excluded here only so that every sublevel component is simply
connected (the general bound of Section~\ref{ss:gen} does not need this
restriction).  Under {\rm (H2)} the maximum $M$ is attained on $\partial\Om_0$,
so the boundary value is $c_0=M$; consequently $\psi=M$ on $\partial\Om$ for
every $\eta\in\widehat{\mathcal G}_2$, and below we may freely use
$\psi=M$ on $\partial\Om$ with $M=\max_{\overline\Om}\psi$.  Let
$t\in(m,M)$ be a regular
value of $\psi_0$, that is, one with $\nabla\psi_0\neq0$ on
$\{\psi_0=t\}$.  Then $\{\psi_0<t\}$ has finitely many connected components;
write $k(t)$ for their number and
\[
a_1(t),\dots,a_{k(t)}(t),\qquad \sum_{i=1}^{k(t)}a_i(t)=A_0(t) ,
\]
for their areas.  Note that $k(t)\ge2$ for some $t$ precisely when $\psi_0$
has a saddle: a new component is born at each minimum, two components merge
at each saddle, and $\{\psi_0<t\}$ is connected for $t$ close to $M$.

\begin{lemma}[The component structure is invariant]\label{l:invar}
Under {\rm (H2)}, for every regular value $t$ each component of
$\{\psi_0<t\}$ is a topological disk.  In addition, for every
$\eta\in\widehat{\mathcal G}_2$ and $\psi=\psi_0\circ\eta^{-1}$, the set
$\{\psi<t\}$ has exactly $k(t)$ connected components, with the same mutual
inclusion relations as those of $\{\psi_0<t\}$, and their areas are
$a_1(t),\dots,a_{k(t)}(t)$.
\end{lemma}

\begin{proof}
We first note that $\{\psi_0<t\}$ is compactly contained in $\Om_0$: on
$\partial\Om_0$ one has $\psi_0=M>t$.  Since $t$ is a regular value, each
component $\Gamma$ of $\{\psi_0=t\}$ is a Jordan curve in $\Om_0$, and simple
connectivity of $\Om_0$ puts the disk $D_\Gamma$ it bounds inside $\Om_0$
(an annulus would already provide a counterexample to this step).
We claim that $\psi_0<t$ throughout $D_\Gamma$.  Let
$\mu=\max_{\overline{D_\Gamma}}\psi_0$.  Since $\psi_0=t$ on
$\Gamma=\partial D_\Gamma$, one has $\mu\ge t$; if $\mu>t$ then $\mu$ is
attained at an interior point of $D_\Gamma$, an interior local maximum of
$\psi_0$, which {\rm (H2)} excludes.  Hence $\mu=t$.  If $\psi_0(x_0)\ge t$ at an
interior point $x_0\in D_\Gamma$, then either $\psi_0(x_0)>t$, which
contradicts $\mu=t$ at once, or $\psi_0(x_0)=t$, so that $x_0$ lies on a
level component other than $\Gamma$ and transversality at the regular level
$t$ ($\nabla\psi_0\ne0$) puts points with $\psi_0>t$ arbitrarily close to
$x_0$, still inside $D_\Gamma$; in both cases $\mu>t$.  Therefore
$\psi_0<t$ on $D_\Gamma$.

Next, no level component $\Gamma'$ different from $\Gamma$ is contained in
$D_\Gamma$: by the claim just applied to $\Gamma'$, transversality at
$\Gamma'$ puts $\psi_0>t$ just outside $D_{\Gamma'}$, an open subset of
$D_\Gamma$ on which the same claim gives $\psi_0<t$, a contradiction.  Two
Jordan curves in the plane are either nested or bound disjoint disks, so
nestedness being excluded, the disks $D_\Gamma$ are pairwise disjoint.  The set
$\Om_0\setminus\bigcup_\Gamma\Gamma$ is a finite union of regions $R$ on each
of which $\psi_0\ne t$, so each $R$ is contained either in $\{\psi_0<t\}$ or
in $\{\psi_0>t\}$.  A region $R\subset\{\psi_0<t\}$ cannot accumulate on
$\partial\Om_0$ (where $\psi_0=M>t$), so its boundary is contained in the
union of the level curves.  Let $\Gamma\subseteq\partial R$.  By
transversality $R$ lies, near $\Gamma$, on the side of $\Gamma$ where
$\psi_0<t$, that is, inside $D_\Gamma$; since $R$ is connected and disjoint
from $\Gamma$, it cannot meet both sides of $\Gamma$, so $R\subseteq D_\Gamma$.
By the preceding paragraphs no level component lies inside $D_\Gamma$ and
none strictly contains $D_\Gamma$, so $D_\Gamma$ is itself one of the
regions; as $R$ is a region contained in $D_\Gamma$ and meeting its
interior, $R=D_\Gamma$.  Hence
$\{\psi_0<t\}$ is the disjoint union of the disks $D_\Gamma$.  Therefore
every component of $\{\psi_0<t\}$ is a topological disk.

For the invariance, $\{\psi<t\}=\eta(\{\psi_0<t\})$ and $\eta$ is a
homeomorphism of $\overline\Om_0$ onto $\overline\Om$, so it maps the
connected components of $\{\psi_0<t\}$ bijectively onto those of
$\{\psi<t\}$; and $|\eta(E)|=|E|$ for every measurable $E\subset\Om_0$
because $\det D\eta=1$ a.e., so the areas are preserved.
\end{proof}

\begin{theorem}[Component lower bound]\label{t:multi}
Assume {\rm (H2)}.  Then for every $\eta\in\widehat{\mathcal G}_2$ the
magnetic energy alone satisfies
\begin{equation}\label{e:Ecomp}
\tfrac12\int_{\Om_0}|D\eta\,\Bv_0|^2\,\dd x\ \ge\
\tfrac12\int_m^M\frac{4\pi\bigl(\sum_{i=1}^{k(t)}\sqrt{a_i(t)}\bigr)^2}{A_0'(t)}
\,\dd t\ =:\ M_{\mathrm{comp}} ,
\end{equation}
and consequently
$E(\eta)\ge M_{\mathrm{comp}}+\kap\,2\sqrt{\pi V}$.
For a single component ($k(t)=1$ for a.e.\ $t$) the constant agrees with that
of \eqref{e:Emin}, that is, $M_{\mathrm{comp}}=M_*$, and in general
$M_{\mathrm{comp}}\ge M_*$, with equality if and only if $k(t)=1$ for a.e.\
$t$; in particular $M_{\mathrm{comp}}>M_*$ strictly whenever $\psi_0$ has a
saddle.
\end{theorem}

\begin{proof}
Let $\psi=\psi_0\circ\eta^{-1}$ and let $t$ be a regular value.  By
Lemma~\ref{l:invar} the set $\{\psi<t\}$ has exactly $k(t)$ components
$Q_1(t),\dots,Q_{k(t)}(t)$, each a topological disk, with
$|Q_i(t)|=a_i(t)$.  Since $t<M$ and $\psi=M$ on $\partial\Om$, the set
$\{\psi<t\}$ is compactly contained in $\Om$; hence
$\partial^*\{\psi<t\}\subseteq\{\psi=t\}$ and, the perimeter being additive
over the components,
\[
\ell(t)=\mathcal H^1(\{\psi=t\})\ \ge\
\Per(\{\psi<t\})=\sum_{i=1}^{k(t)}\Per\bigl(Q_i(t)\bigr) .
\]
Each $Q_i(t)$ has area $a_i(t)$, so the isoperimetric inequality gives
$\Per(Q_i(t))\ge2\sqrt{\pi a_i(t)}$ and therefore
\[
\ell(t)\ \ge\ 2\sqrt\pi\sum_{i=1}^{k(t)}\sqrt{a_i(t)} .
\]
Substituting this bound into \eqref{e:csleaves} and integrating over $t$
exactly as in the proof of Lemma~\ref{l:ps} yields \eqref{e:Ecomp}; the
level sets for $t\notin(m,M)$ are empty and contribute nothing.  For the
shape term, $|\eta(\Om_0)|=V$ and the isoperimetric inequality gives
$\Per(\eta(\Om_0))\ge2\sqrt{\pi V}$.

Finally,
\[
4\pi\Bigl(\sum_{i=1}^{k}\sqrt{a_i}\Bigr)^2
=4\pi\sum_{i=1}^{k}a_i+8\pi\sum_{i<j}\sqrt{a_ia_j}
\ \ge\ 4\pi A_0 ,
\]
with equality, at a fixed regular level $t$, if and only if at most one of the
$a_i(t)$ is nonzero.  Hence $M_{\mathrm{comp}}\ge M_*$, with equality if and
only if $k(t)=1$ for a.e.\ $t$, and the inequality is strict as soon as
$k(t)\ge2$ on a set of positive measure, which happens exactly when $\psi_0$
has a saddle.
\end{proof}

\begin{proposition}[The bound is not attained across a saddle]\label{p:notatt}
Assume {\rm (H2)} and let $s$ be a critical value of $\psi_0$ at which the
following happens: there is a unique saddle point $q$ with $\psi_0(q)=s$, and
as $t\uparrow s$ the sublevel set $\{\psi_0<t\}$ has exactly two components
in a neighborhood of $q$, which merge at $t=s$; no other critical point has
value in a neighborhood of $s$.  Write $a_1,a_2>0$ for the limiting areas of
those components as $t\uparrow s$ (and $a_i(t)$ for their areas at level
$t<s$).  Then no $\eta\in\widehat{\mathcal G}_2$
attains equality in \eqref{e:Ecomp}.
\end{proposition}

\begin{proof}
Suppose $\eta$ attains equality in \eqref{e:Ecomp}.  Write
$\psi=\psi_0\circ\eta^{-1}$; by Lemma~\ref{l:invar} its distribution
function is $A_0$.  The bound is the integral over $t$ of nonnegative
pointwise quantities: at each regular level $t$, with $Q_i(t)$ the
components of $\{\psi<t\}$, one has \eqref{e:csleaves} and
\[
\ell(t)\ :=\ \mathcal H^1(\{\psi=t\})\ \ge\
\sum_i\Per\bigl(Q_i(t)\bigr)\ \ge\
2\sqrt\pi\sum_i\sqrt{a_i(t)} .
\]
Since the integrands are nonnegative and
\[
\int_{\Om_0}|D\eta\,\Bv_0|^2\,\dd x
=\int_\Om|\nabla\psi|^2\,\dd y
\ \ge\ \int_m^M\frac{\ell(t)^2}{A_0'(t)}\,\dd t
\ \ge\ 2M_{\mathrm{comp}} ,
\]
equality throughout forces equality in both pointwise inequalities for a.e.\
regular level $t$, so that for a.e.\ $t$ each component $Q_i(t)$ is a disk
(the equality case of the isoperimetric inequality) and $|\nabla\psi|$ is
constant on each component of $\{\psi=t\}$.

Let $t_0<s<t_1$ be levels of $\psi$ chosen in the
full-measure set of the preceding paragraph, with $s$ the saddle value.

Since $\psi\in W^{1,\infty}(\Om)$ by Lemma~\ref{l:conj}(i), the coarea
formula, the Cauchy--Schwarz step, finite-perimeter sublevel boundaries and
the equality conditions above are available on that full-measure set of
levels.  What $\psi$ inherits from $\psi_0$ is its
distribution function and its level topology, namely the number and areas of the
sublevel components and the values at which they merge
(Lemma~\ref{l:invar}), though not the shapes of its level curves.
The two components of
$\{\psi<t_0\}$ that merge at level $s$ are then disks $D_1,D_2$ of areas
$a_1(t_0),a_2(t_0)$, and the component $F_{t_1}$ of $\{\psi<t_1\}$ containing
$D_1\cup D_2$ is a disk as well.  Two
disjoint disks of radii $r_1,r_2$ are contained in a disk of radius $R$ only
if $R\ge r_1+r_2$: the segment joining the centers has length at least
$r_1+r_2$, and extending it by $r_1$ at one end and $r_2$ at the other stays
inside the containing disk, so $2R\ge2(r_1+r_2)$.  Hence
\[
|F_{t_1}|\ \ge\ \bigl(\sqrt{a_1(t_0)}+\sqrt{a_2(t_0)}\bigr)^2 .
\]
But $|F_{t_1}|\to a_1+a_2$ as $t_1\downarrow s$, whereas
$(\sqrt{a_1(t_0)}+\sqrt{a_2(t_0)})^2\to(\sqrt{a_1}+\sqrt{a_2})^2
=a_1+a_2+2\sqrt{a_1a_2}>a_1+a_2$ as $t_0\uparrow s$.  Choosing
$t_0<s<t_1$ close enough to $s$ therefore gives
$|F_{t_1}|<(\sqrt{a_1(t_0)}+\sqrt{a_2(t_0)})^2$, a contradiction.
\end{proof}

The packing step used above is the equality case of a statement that is exact
in every configuration, and the exact form is what a quantified gap needs.

\begin{lemma}[Quantified packing]\label{l:packing}
Let $B(c,R)\supset B(c_1,r_1)\cup B(c_2,r_2)$ with $r_1,r_2>0$ and
$B(c_1,r_1)\cap B(c_2,r_2)=\varnothing$.  Then
\begin{equation}
R-r_1-r_2=\tfrac12\bigl[d_1+d_2+d_3+d_4\bigr],
\end{equation}
where
\[
d_1=R-|c_1-c|-r_1,\qquad d_2=R-|c_2-c|-r_2 ,
\]
\[
d_3=|c_1-c|+|c_2-c|-|c_1-c_2|,\qquad d_4=|c_1-c_2|-r_1-r_2 ,
\]
and $d_1,d_2,d_3,d_4\ge0$.  In particular $R=r_1+r_2$ if and only if the four
slacks vanish, that is, if and only if the two disks are tangent to the outer
circle and to each other, with $c$ on the segment $[c_1,c_2]$.
\end{lemma}

\begin{proof}
The four slacks are nonnegative by containment, the triangle inequality and
disjointness, and their sum equals $2R-2r_1-2r_2$.
\end{proof}

For the isolated two-component merge covered by
Proposition~\ref{p:notatt}, the infimum is in fact strictly larger than the
component bound.  The argument below does not cover simultaneous or interacting
critical events.

\begin{proposition}[A strict gap across a saddle]\label{p:strictgap}
Assume {\rm (H2)} and let $s$ be a critical value of $\psi_0$ satisfying the
same isolated simple-merge condition as in Proposition~\ref{p:notatt}:
there is a unique saddle point of value $s$, exactly two sublevel components
merge there, and no other critical value lies near $s$; let their limiting
areas be $a_1,a_2>0$ as $t\uparrow s$.  Then
\[
\inf_{\eta\in\widehat{\mathcal G}_2}E(\eta)\ >\
M_{\mathrm{comp}}+\kap\,2\sqrt{\pi V}
\]
strictly: the component bound of Theorem~\ref{t:multi} is not sharp at the
isolated simple merge considered here.
\end{proposition}

\begin{proof}
Suppose the contrary and let $\eta_n$ be a minimizing sequence with
$E(\eta_n)\to M_{\mathrm{comp}}+\kap\,2\sqrt{\pi V}$.  Write
$\Om_n=\eta_n(\Om_0)$, $\psi_n=\psi_0\circ\eta_n^{-1}\in
W^{1,\infty}(\Om_n)$ by Lemma~\ref{l:conj}(i), and for a regular
value $t$ let $Q_{1,n}(t),\dots,Q_{k(t),n}(t)$ be the components of
$\{\psi_n<t\}$.

\emph{1. Both terms converge separately.}  By Theorem~\ref{t:multi} the
magnetic part is at least $M_{\mathrm{comp}}$, and by the isoperimetric
inequality $\Per(\Om_n)\ge2\sqrt{\pi V}$.  Since their sum tends to
$M_{\mathrm{comp}}+\kap\,2\sqrt{\pi V}$, both do:
\[
\tfrac12\int_{\Om_0}|D\eta_n\,\Bv_0|^2\,\dd x\to M_{\mathrm{comp}},
\qquad \Per(\Om_n)\to2\sqrt{\pi V}.
\]

\emph{2. Vanishing of the deficit forces roundness a.e.}  In the proof of
Theorem~\ref{t:multi} one has, with $\ell_n(t)=\mathcal H^1(\{\psi_n=t\})$
and $\psi_n=\psi_0\circ\eta_n^{-1}$, the chain of inequalities
\[
\int_{\Om_0}|D\eta_n\,\Bv_0|^2\,\dd x
=\int_{\Om_n}|\nabla\psi_n|^2\,\dd y
\ \ge\ \int_m^M \frac{\ell_n(t)^2}{A_0'(t)}\,\dd t
\ \ge\ \int_m^M\frac{4\pi\bigl(\sum_{i=1}^{k(t)}\sqrt{a_i(t)}\bigr)^2}{A_0'(t)}\,\dd t
=2M_{\mathrm{comp}} ;
\]
the first inequality is \eqref{e:coarea}--\eqref{e:csleaves} applied to
$\psi_n$, with distribution function $A_0$ by
Lemma~\ref{l:invar}; the second is the pair
$\ell_n(t)\ge\sum_i\Per(Q_{i,n}(t))$ and
$\Per(Q_{i,n}(t))\ge2\sqrt{\pi a_i(t)}$ inserted pointwise.  The left and
right ends differ by
$\int_{\Om_0}|D\eta_n\Bv_0|^2\,\dd x-2M_{\mathrm{comp}}$, which tends to $0$
by step~1.  Writing the total deficit as
\[
\int_{\Om_0}|D\eta_n\,\Bv_0|^2\,\dd x
-\int_m^M\frac{\ell_n(t)^2}{A_0'(t)}\,\dd t
+\int_m^M\Bigl[\frac{\ell_n(t)^2}{A_0'(t)}
-\frac{4\pi\bigl(\sum_i\sqrt{a_i(t)}\bigr)^2}{A_0'(t)}\Bigr]\dd t ,
\]
both terms are integrals of nonnegative integrands and the sum tends to $0$;
hence each tends to $0$.  Passing to a subsequence, for a.e.\ $t\in(m,M)$,
\[
\ell_n(t)\ \to\ 2\sqrt\pi\sum_{i=1}^{k(t)}\sqrt{a_i(t)},
\]
and, since $\ell_n(t)\ge\sum_i\Per(Q_{i,n}(t))\ge2\sqrt\pi\sum_i\sqrt{a_i(t)}$,
also $\Per(Q_{i,n}(t))\to2\sqrt{\pi a_i(t)}$ for each $i$ and a.e.\ $t$.

\emph{3. Passage to the limit: everything becomes a ball.}  We use the
sharp quantitative isoperimetric inequality \cite{FuscoMaggiPratelli2008} (in the plane
explicit-constant quantitative isoperimetric inequalities are classical, e.g.\ the
Bonnesen inequality \cite{Bonnesen1924}; we keep the dimension-free form above for
uniformity):
a set $F$ of area $a$ obeys $\Per(F)\ge2\sqrt{\pi a}\bigl(1+c\,\mathcal
A(F)^2\bigr)$, where $\mathcal A(F)=\inf_B|F\Delta B|/|F|$ is the Fraenkel
asymmetry, the infimum being over balls of area $a$.  Consequently, if
$|F_n|=a$ and $\Per(F_n)\to2\sqrt{\pi a}$, then $\mathcal A(F_n)\to0$, i.e.\
there are balls $B_n$ of area $a$ with
$\|\mathbf 1_{F_n}-\mathbf 1_{B_n}\|_{L^1}\to0$.

Apply this to $\Om_n$: for each $n$ translate $\Om_n$ so that its
approximating ball is centered at the origin, giving sets that we again denote
$\Om_n$, so that $\mathbf 1_{\Om_n}\to\mathbf 1_{B(0,R)}$ in $L^1$ with
$R=\sqrt{V/\pi}$.  The sets $Q_{i,n}(t_0)$ are translated by the same amount;
inclusions, disjointness and $L^1$ convergence are preserved by a common
translation, so the argument below applies to the translated sets.  Choose $t_0<s<t_1$, both in the full-measure set of
step~2, so close to $s$ that
\begin{equation}\label{e:window}
|E_{t_1}|\ <\ \bigl(\sqrt{a_1(t_0)}+\sqrt{a_2(t_0)}\bigr)^2 ,
\end{equation}
where $E_{t_1}$ is the component of $\{\psi_0<t_1\}$ containing the two
merging wells.  This is possible because $|E_{t_1}|\to a_1+a_2$ as
$t_1\downarrow s$, while
$(\sqrt{a_1(t_0)}+\sqrt{a_2(t_0)})^2\to a_1+a_2+2\sqrt{a_1a_2}>a_1+a_2$ as
$t_0\uparrow s$.

Because $t_0,t_1$ lie in the full-measure set of step~2, the sublevel
components at these two levels are well defined for all large $n$.
Let $W_n$ be the
component of $\{\psi_n<t_1\}$ containing the two merging wells
$Q_{1,n}(t_0),Q_{2,n}(t_0)$ (translated as above); by step~2 applied at the
fixed level $t_1$, $\Per(W_n)\to2\sqrt{\pi|E_{t_1}|}$, since $|W_n|=|E_{t_1}|$
(by Lemma~\ref{l:invar} the areas of the components of $\{\psi_n<t_1\}$ agree
with those of $\psi_0$).
The quantitative inequality therefore gives
$\|\mathbf 1_{W_n}-\mathbf 1_{B(c_n,R_1)}\|_1\to0$ with
$R_1=\sqrt{|E_{t_1}|/\pi}$; applying step~2 at the fixed level $t_0$
likewise gives
$\|\mathbf 1_{Q_{i,n}(t_0)}-\mathbf 1_{B(d_{i,n},r_i)}\|_1\to0$ with
$r_i=\sqrt{a_i(t_0)/\pi}$.  Since
$Q_{i,n}(t_0)\subseteq W_n\subseteq\Om_n$ and almost all of $\Om_n$ lies in
$B(0,R)$, the centers $c_n$ and $d_{i,n}$ are bounded; passing to a
subsequence with $c_n\to c$ and $d_{i,n}\to d_i$ and using the continuity of
translations in $L^1$,
\[
\mathbf 1_{W_n}\to\mathbf 1_{B(c,R_1)},\qquad
\mathbf 1_{Q_{i,n}(t_0)}\to\mathbf 1_{B(d_i,r_i)}
\quad\text{in }L^1 .
\]

\emph{4. The limit violates the packing constraint.}  From
$Q_{i,n}(t_0)\subseteq W_n$ we get, in the limit,
$B(d_i,r_i)\subseteq B(c,R_1)$ up to a null set; since both are open balls,
the inclusion holds as sets.  From $Q_{1,n}(t_0)\cap Q_{2,n}(t_0)=\varnothing$ and
$\|\mathbf 1_{A_1}\mathbf 1_{A_2}-\mathbf 1_{B_1}\mathbf 1_{B_2}\|_1
\le\|\mathbf 1_{A_1}-\mathbf 1_{B_1}\|_1+\|\mathbf 1_{A_2}-\mathbf 1_{B_2}\|_1$,
applied with $A_i=Q_{i,n}(t_0)$ and $B_i=B(d_i,r_i)$, the $L^1$ convergence
above gives $|B(d_1,r_1)\cap B(d_2,r_2)|=0$, hence
$|d_1-d_2|\ge r_1+r_2$.  Containment in $B(c,R_1)$ forces, measuring the
extent along the direction of $d_2-d_1$,
\[
|d_1-d_2|+r_1+r_2\ \le\ 2R_1 ,
\]
so that $R_1\ge r_1+r_2$, that is
\[
|E_{t_1}|\ \ge\ \bigl(\sqrt{a_1(t_0)}+\sqrt{a_2(t_0)}\bigr)^2 ,
\]
contradicting \eqref{e:window}.
\end{proof}

The proof of Proposition~\ref{p:strictgap} uses the packing constraint only in
its equality case.  Quantifying it needs two elementary geometric estimates,
and then the two propositions that combine them with the analytic steps above.

\begin{lemma}[Cap and lens]\label{l:caplens}
Let $0\le h\le r$ and let $\mathrm{cap}_r(h)$ be the area of the part of a
disk of radius $r$ cut off by a chord at distance $r-h$ from its center.  Then
\begin{equation}\label{e:cap}
\mathrm{cap}_r(h)\ \ge\ \frac{\pi}{2}\,\sqrt r\,h^{3/2}.
\end{equation}
If two disks of radii $r_1,r_2$ have centers at distance $r_1+r_2-s$ with
$0\le s\le\min(r_1,r_2)$, then their intersection $\mathcal L$ satisfies
\begin{equation}\label{e:lens}
|\mathcal L|\ \ge\ \frac{\pi}{2\sqrt2}\,\sqrt{\min(r_1,r_2)}\;s^{3/2}.
\end{equation}
\end{lemma}

\begin{proof}
With $h=rs$ one has $\mathrm{cap}_r(h)=r^2f(s)$, where
$f(s)=\arccos(1-s)-(1-s)\sqrt{2s-s^2}$ and $f'(s)=2\sqrt{s(2-s)}$.  Hence
$F(s):=f(s)-\tfrac{\pi}{2}s^{3/2}$ has $F(0)=F(1)=0$ and
$F'(s)=\sqrt s\,[2\sqrt{2-s}-3\pi/4]$, whose bracket decreases from
$2\sqrt2-3\pi/4>0$ to $2-3\pi/4<0$; so $F$ increases and then decreases, and therefore
$F\ge0$, that is, \eqref{e:cap}.  For two overlapping disks the intersection is
the union of the two caps cut off by the radical axis, of heights $h_1,h_2\ge0$
with $h_1+h_2=s$; by \eqref{e:cap} and the convexity of $t\mapsto t^{3/2}$,
\[
|\mathcal L|\ \ge\ \frac{\pi}{2}\,\sqrt{\min(r_1,r_2)}\,
\bigl(h_1^{3/2}+h_2^{3/2}\bigr)
\ \ge\ \frac{\pi}{2\sqrt2}\,\sqrt{\min(r_1,r_2)}\;s^{3/2}. \qedhere
\]
\end{proof}

\begin{proposition}[Quantitative packing for approximate disks]%
\label{p:qpack}
Let $E\supset E_1\cup E_2$ with $E_1\cap E_2=\varnothing$, and let
$B(c,R)$ and $B(c_i,r_i)$ be disks with
$R^2=|E|/\pi$, $r_i^2=|E_i|/\pi$,
$|E\triangle B(c,R)|\le \varepsilon$, $|E_i\triangle B(c_i,r_i)|\le \varepsilon_i$,
$|c_1-c_2|\ge\max(r_1,r_2)$, and
\[
\varepsilon+\varepsilon_i<\tfrac{\pi}{2}r_i^2\ \ (i=1,2).
\]
Then
\begin{equation}\label{e:qpack}
R\ \ge\ r_1+r_2-\tfrac12\bigl(\kappa_1+\kappa_2+\kappa_4\bigr),\qquad
\kappa_i=\Bigl(\frac{2(\varepsilon+\varepsilon_i)}{\pi\sqrt{r_i}}\Bigr)^{2/3},\quad
\kappa_4=\Bigl(\frac{2\sqrt2\,(\varepsilon_1+\varepsilon_2)}{\pi\sqrt{\min(r_1,r_2)}}\Bigr)^{2/3}.
\end{equation}
\end{proposition}

\begin{proof}
For the three disks define $d_1,\dots,d_4$ by the formulas of
Lemma~\ref{l:packing}; their sum is always $2(R-r_1-r_2)$, and $d_3\ge0$ by the
triangle inequality.  By the smallness assumption, if $d_i<0$ then
$|d_i|\le r_i$: otherwise $|c_i-c|>R$, so the center of $B(c_i,r_i)$ lies
outside $B(c,R)$ and more than half of $B(c_i,r_i)$ lies outside it; hence
$\varepsilon+\varepsilon_i>\pi r_i^2/2$, contradicting the strict
smallness assumption.  If $d_4<0$, the
hypothesis $|c_1-c_2|\ge\max(r_1,r_2)$ gives $|d_4|\le\min(r_1,r_2)$, the range
of the lens estimate of Lemma~\ref{l:caplens}.  Thus the cap and lens estimates
of Lemma~\ref{l:caplens} may be applied below.  If $d_1<0$, the
disk $B(c_1,r_1)$ protrudes beyond $B(c,R)$ by depth $|d_1|$; the part of it
outside $B(c,R)$ contains the cap cut off by the half-plane of $B(c,R)$
orthogonal to $c_1-c$, so by \eqref{e:cap}
\[
\frac{\pi}{2}\sqrt{r_1}\,|d_1|^{3/2}\ \le\ |B(c_1,r_1)\setminus B(c,R)|\ \le\ \varepsilon+\varepsilon_1 ,
\]
that is, $d_1\ge-\kappa_1$; if $d_1\ge0$ the bound is trivial.  Likewise
$d_2\ge-\kappa_2$.  If $d_4<0$ the two disks overlap with depth $|d_4|$; by
$|c_1-c_2|\ge\max(r_1,r_2)$ that depth is at most $\min(r_1,r_2)$, the range of
the lens estimate \eqref{e:lens}, and
$B_1\cap B_2\subset(B_1\setminus E_1)\cup(B_2\setminus E_2)$ because
$E_1\cap E_2=\varnothing$, so $d_4\ge-\kappa_4$; if $d_4\ge0$ this is trivial.
Adding the four bounds gives
$2(R-r_1-r_2)\ge-(\kappa_1+\kappa_2+\kappa_4)$.
\end{proof}

\begin{proposition}[Quantified gap]\label{p:qgap}
In the situation of Proposition~\ref{p:strictgap} there is a constant
$c_*>0$, whose dependence on $\psi_0$ is explicit up to the universal constant
$c_F$ and the data-dependent constant $C_2$ below, such that
\[
\inf_{\widehat{\mathcal G}_2}E\ \ge\
M_{\mathrm{comp}}+\kap\,2\sqrt{\pi V}+c_* .
\]
More precisely, put $\rho_0=2\sqrt{a_1a_2}>0$, where
$a_i=\lim_{t\uparrow s}a_i(t)$, choose $\delta>0$ so small that
$(\sqrt{a_1(t_0)}+\sqrt{a_2(t_0)})^2-|E_{t_1}|\ge\rho_0/2$ for every
$t_0\in(s-\delta,s)$ and $t_1\in(s,s+\delta)$, and that $s$ is the only
critical value in $[s-\delta,s+\delta]$.  Set
\[
I_-=(s-\delta,s-\delta/2),\qquad I_+=(s+\delta/2,s+\delta),\qquad
I=I_-\cup I_+ ,
\]
put $\eta_0=\delta/4$, let $\alpha=\inf_I\min_i\sqrt{a_i}$,
$\beta=\inf_I\sum_i\sqrt{a_i}$ and
$\gamma_\delta=\sup_I A_0'<\infty$, and let
$R_{\max}=\max_{t\in[s,s+\delta]}\sqrt{|E_t|/\pi}$.  Then
\begin{equation}\label{e:qgap}
\inf_{\widehat{\mathcal G}_2}E-\bigl(M_{\mathrm{comp}}+\kap2\sqrt{\pi V}\bigr)
\ \ge\ \frac{C_1\eta_0}{2}\Bigl(\frac{\rho_0}{6\pi C_2R_{\max}^2}\Bigr)^{3},
\qquad C_1=\frac{4\pi c_F\alpha\beta}{\gamma_\delta},
\end{equation}
where $c_F$ is the constant of the quantitative isoperimetric inequality
\cite{FuscoMaggiPratelli2008} and $C_2$ is the explicit constant
\eqref{e:C2}, depending only on $R_{\max}$, $\rho_0$ and the extrema of
$r_i(t)=\sqrt{a_i(t)/\pi}$ over $t\in[s-\delta,s]$.
\end{proposition}

\begin{proof}
Let $E_*=M_{\mathrm{comp}}+\kap2\sqrt{\pi V}$ and put
\[
c_*=\frac{C_1\eta_0}{2}\Bigl(\frac{\rho_0}{6\pi C_2R_{\max}^2}\Bigr)^3 .
\]
Let $r_{\min}$ and $r_{\max}$ be the minimum and maximum of
$r_i(t)=\sqrt{a_i(t)/\pi}$ over the two merging wells and
$t\in[s-\delta,s]$; then $r_{\min}>0$.  For $t>s$ the merged component
satisfies $|E_t|\ge a_1+a_2\ge a_i(t')$ for $t'\le s$, so
$r_{\max}\le R_{\max}$.  Put
\[
\rho_*=\min\Bigl\{\frac{r_{\min}^2}{2(R_{\max}^2+r_{\max}^2)},
\ \frac{r_{\min}^2}{4\sqrt2\,r_{\max}^2}\Bigr\},\qquad
C_2^0=(4\sqrt2)^{2/3},\qquad
K=2^{1/3}r_{\max}r_{\min}^{-1/3}R_{\max}^{-2/3},
\]
and take
\begin{equation}\label{e:C2}
C_2=\max\Bigl\{C_2^0,\ KC_2^0,\
\frac{\rho_0}{6\pi R_{\max}^2(2\rho_*)^{2/3}}\Bigr\}.
\end{equation}
Then
\begin{equation}\label{e:csmall}
\frac12\Bigl(\frac{\rho_0}{6\pi C_2R_{\max}^2}\Bigr)^{3/2}\le\rho_*
\end{equation}
and $C_2\ge KC_2^0$.

Suppose, for contradiction, that some $\eta\in\widehat{\mathcal G}_2$ has
$\delta_E:=E(\eta)-E_*<c_*$.  Write $\Om=\eta(\Om_0)$,
$M_{\mathrm{mag}}=\tfrac12\int_{\Om_0}|D\eta\,\Bv_0|^2$ and
$\psi=\psi_0\circ\eta^{-1}$, so that $M_{\mathrm{mag}}=\tfrac12\int_\Om|\nabla\psi|^2$ by
Lemma~\ref{l:conj}.  Isoperimetry gives $\Per(\Om)\ge2\sqrt{\pi V}$, hence
$\delta_E\ge M_{\mathrm{mag}}-M_{\mathrm{comp}}$.

\emph{1. The magnetic deficit.}  Near a nondegenerate saddle the Morse
coordinates give $\psi_0=\psi_0(p)+x_1^2-x_2^2$ with $p$ the critical point,
so \eqref{e:A0prime} gives $A_0'(\psi_0(p)+\tau)\ge c\log(1/\tau)$ for small
$\tau>0$; the same holds below the critical value.  Thus $A_0'$ is bounded on
$I$ and unbounded on every interval containing the critical value, which is
why $\gamma_\delta=\sup_I A_0'$ is finite.  For
a.e.\ $t$ let $Q_i(t)$ be the components of $\{\psi<t\}$, $a_i(t)=|Q_i(t)|$
and $\mathcal A_i(t)$ their Fraenkel asymmetries, and put
\[
S(t)=2\sqrt\pi\sum_i\sqrt{a_i(t)},\qquad
\varsigma(t)=\sum_i\bigl(\Per(Q_i(t))-2\sqrt{\pi a_i(t)}\bigr)\ \ge\ 0 .
\]
Since $\psi\in W^{1,\infty}(\Om)$ by Lemma~\ref{l:conj}(i), the coarea
formula and the Cauchy--Schwarz inequality give, with
$\ell(t)=\mathcal H^1(\{\psi=t\})\ge S(t)+\varsigma(t)$,
\begin{align*}
M_{\mathrm{mag}}-M_{\mathrm{comp}}
&\ \ge\ \int_m^M\frac{\ell(t)^2-S(t)^2}{2A_0'(t)}\,\dd t\\
&\ \ge\ \int_m^M\frac{S(t)\,\varsigma(t)}{A_0'(t)}\,\dd t ,
\end{align*}
because $\ell^2-S^2\ge2S\varsigma$ and $M_{\mathrm{comp}}=\tfrac12\int S^2/A_0'$.

By the quantitative isoperimetric inequality
$\varsigma(t)\ge2\sqrt\pi c_F\sum_i\sqrt{a_i(t)}\,\mathcal A_i(t)^2$; on $I$,
\[
S(t)\varsigma(t)\ge4\pi c_F\alpha\beta\sum_i\mathcal A_i(t)^2,
\qquad A_0'(t)\le\gamma_\delta ,
\]
and therefore
\[
M_{\mathrm{mag}}-M_{\mathrm{comp}}\ \ge\ C_1\int_I
\sum_i\mathcal A_i(t)^2\,\dd t,\qquad C_1=\frac{4\pi c_F\alpha\beta}{\gamma_\delta}.
\]
Since $\delta_E\ge M_{\mathrm{mag}}-M_{\mathrm{comp}}$ and each of $I_-,I_+$ has length
$\delta/2=2\eta_0$, each contains a level $t$ with
$\sum_i\mathcal A_i(t)^2\le2\delta_E/(C_1\delta)=\delta_E/(2C_1\eta_0)$.

\emph{2. Approximating disks.}  Choose such levels $t_0\in I_-$ and
$t_1\in I_+$.  At $t_1$ the component $W$ of the sublevel set containing
the two merging wells satisfies $|W|=|E_{t_1}|$, and at $t_0$ the two
components that merge at $s$ are $Q_1(t_0),Q_2(t_0)$; moreover
$Q_i(t_0)\subset W$ and $Q_1(t_0)\cap Q_2(t_0)=\varnothing$.  By
\cite{FuscoMaggiPratelli2008} there are disks $B(c,R)$, $B(c_i,r_i)$ with
$R^2=|E_{t_1}|/\pi$, $r_i^2=a_i(t_0)/\pi$ and
\[
|W\triangle B(c,R)|\le\pi R^2\mathcal A_W(t_1),\qquad
|Q_i(t_0)\triangle B(c_i,r_i)|\le\pi r_i^2\mathcal A_i(t_0).
\]

\emph{3. Quantitative packing.}  Put
$\rho=\bigl(\delta_E/(2C_1\eta_0)\bigr)^{1/2}$, so that all three asymmetries
are at most $\rho$.  Since $\delta_E<c_*$ and $\delta_E=2C_1\eta_0\rho^2$,
\[
\rho<\Bigl(\frac{c_*}{2C_1\eta_0}\Bigr)^{1/2}
=\frac12\Bigl(\frac{\rho_0}{6\pi C_2R_{\max}^2}\Bigr)^{3/2}\le\rho_* .
\]
The $L^1$ errors of Proposition~\ref{p:qpack} therefore satisfy
$\varepsilon\le\pi R_{\max}^2\rho$, $\varepsilon_i\le\pi r_{\max}^2\rho$ and
$\varepsilon+\varepsilon_i\le\pi(R_{\max}^2+r_{\max}^2)\rho<\frac\pi2
r_{\min}^2\le\frac\pi2 r_i^2$, which is the strict smallness hypothesis of
Proposition~\ref{p:qpack}.  Moreover
$\kappa_4\le(4\sqrt2\,r_{\max}^2\rho/\sqrt{r_{\min}})^{2/3}<r_{\min}$ by
\eqref{e:csmall}.  Since $\varepsilon+\varepsilon_i\le2\pi
R_{\max}^2\rho$ and $\varepsilon_1+\varepsilon_2\le2\pi r_{\max}^2\rho$,
\eqref{e:qpack} and $r_{\max}\le R_{\max}$ give
\[
\kappa_j\le C_2^0\,\rho^{2/3}R_{\max}^{4/3}r_{\min}^{-1/3}\qquad (j=1,2,4).
\]
The hypothesis $|c_1-c_2|\ge\max(r_1,r_2)$ of
Proposition~\ref{p:qpack} holds as soon as $\kappa_4<\min(r_1,r_2)$: if
$|c_1-c_2|<\max(r_1,r_2)$, the lens $B(c_1,r_1)\cap B(c_2,r_2)$ contains the
lens of depth $\min(r_1,r_2)$, so the endpoint $s=\min(r_1,r_2)$ of
\eqref{e:lens} and monotonicity of the lens area in the depth give
$\varepsilon_1+\varepsilon_2\ge(\pi/(2\sqrt2))\min(r_1,r_2)^2$, which by the
definition of $\kappa_4$ gives $\kappa_4\ge\min(r_1,r_2)\ge r_{\min}$.  In
the relevant case Proposition~\ref{p:qpack} gives
\[
R\ \ge\ r_1+r_2-\tfrac12(\kappa_1+\kappa_2+\kappa_4),\qquad
\kappa_j\ \le\ C_2^0\,\rho^{2/3}R_{\max}^{4/3}r_{\min}^{-1/3},
\]
with $r_{\min}$ and $r_{\max}$ as defined at the start of the proof.

\emph{4. The deficit.}  By the choice of $\delta$,
$\tfrac12\rho_0\le(\sqrt{a_1(t_0)}+\sqrt{a_2(t_0)})^2-|E_{t_1}|
=\pi\bigl((r_1+r_2)^2-R^2\bigr)$, hence, since $R\le r_1+r_2\le2r_{\max}$,
\[
r_1+r_2-R\ =\ \frac{(r_1+r_2)^2-R^2}{r_1+r_2+R}\ \ge\
\frac{\rho_0}{2\pi(r_1+r_2+R)}\ \ge\ \frac{\rho_0}{8\pi r_{\max}} .
\]
Comparing with Step 3,
\[
\tfrac32C_2^0\rho^{2/3}R_{\max}^{4/3}r_{\min}^{-1/3}\ \ge\
\frac{\rho_0}{8\pi r_{\max}},
\]
that is, $\rho^2\ge(\rho_0/(12\pi C_2^0))^3\,r_{\min}r_{\max}^{-3}R_{\max}^{-4}$.
Since $\delta_E=2C_1\eta_0\rho^2$, and
$K^3=2r_{\max}^3/(r_{\min}R_{\max}^2)$,
\[
\delta_E\ \ge\ \frac{C_1\eta_0\rho_0^3 r_{\min}}
{864\pi^3(C_2^0)^3r_{\max}^3R_{\max}^4}
=\frac{C_1\eta_0}{2}\Bigl(\frac{\rho_0}{6\pi KC_2^0R_{\max}^2}\Bigr)^3
\ \ge\ \frac{C_1\eta_0}{2}\Bigl(\frac{\rho_0}{6\pi C_2R_{\max}^2}\Bigr)^3=c_* ,
\]
contradicting $\delta_E<c_*$.  As $\eta$ was arbitrary, \eqref{e:qgap} follows.
\end{proof}

\smallskip
\noindent\emph{The constants in \eqref{e:qgap}.}  The constants of the proof
are summarized below; apart from the universal constant $c_F$ all of them are
determined by $\psi_0$.

\begin{samepage}
\begin{center}
\small
\renewcommand{\arraystretch}{1.35}
\begin{tabular}{@{}>{\raggedright}p{2.9cm}>{\raggedright}p{6.4cm}>{\raggedright\arraybackslash}p{4.9cm}@{}}
\hline
quantity & definition in the proof & data it depends on \\
\hline
$\rho_0$ & $2\sqrt{a_1a_2}$, with $a_i=\lim_{t\uparrow s}a_i(t)$
& the two well areas at the saddle value $s$ \\
$\delta$, $\eta_0$ & separation condition; $\eta_0=\delta/4$
& $\rho_0$ and $\psi_0$ near $s$ \\
$\alpha,\beta,\gamma_\delta$ & $\inf_I\min_i\sqrt{a_i}$,
$\inf_I\sum_i\sqrt{a_i}$, $\sup_I A_0'$ & $\psi_0$ and $\delta$ \\
$C_1$ & $4\pi c_F\alpha\beta/\gamma_\delta$ & $c_F$, $\psi_0$, $\delta$ \\
$r_{\min},r_{\max}$ & extrema of $r_i(t)=\sqrt{a_i(t)/\pi}$ over
$t\in[s-\delta,s]$ & $\psi_0$, $\delta$ \\
$R_{\max}$ & $\max_{t\in[s,s+\delta]}\sqrt{|E_t|/\pi}$ & $\psi_0$, $\delta$ \\
$\rho_*$, $C_2^0$, $K$ & $\rho_*$ as in the proof,
$C_2^0=(4\sqrt2)^{2/3}$, $K=2^{1/3}r_{\max}r_{\min}^{-1/3}R_{\max}^{-2/3}$
& $r_{\min},r_{\max},R_{\max}$ \\
$C_2$ & $\max\{C_2^0,KC_2^0,\rho_0/(6\pi R_{\max}^2(2\rho_*)^{2/3})\}$,
see \eqref{e:C2} & the rows above \\
\hline
\end{tabular}
\end{center}
\end{samepage}
The three packing constants are then estimated at once,
$\kappa_j\le C_2^0\rho^{2/3}R_{\max}^{4/3}r_{\min}^{-1/3}$ for $j=1,2,4$,
where $\rho^2=\delta_E/(2C_1\eta_0)$.  The smallness hypothesis
$\varepsilon+\varepsilon_i<\tfrac\pi2r_i^2$ of Proposition~\ref{p:qpack}
follows from the first alternative in the definition of $\rho_*$ together
with $\rho<\rho_*$ of \eqref{e:csmall}, and the center separation
$|c_1-c_2|\ge\max(r_1,r_2)$ from the second alternative, through
$\kappa_4<r_{\min}$ and the contrapositive argument of Step~3.  Since
$C_2\ge KC_2^0$, the last chain of inequalities in the proof gives
$c_*=\tfrac12C_1\eta_0\bigl(\rho_0/(6\pi C_2R_{\max}^2)\bigr)^3$ as stated.

\subsection{Interior local maxima: the general lower bound}%
\label{ss:gen}

Hypothesis (H2) excludes interior local maxima of $\psi_0$, and with them
sublevel components that are not simply connected.  We now remove (H2) and
obtain a bound when all interior critical points are nondegenerate and the
boundary gradient is nonzero.

Let $c_0=\psi_0|_{\partial\Om_0}$ be the boundary value, and write
$m=\min_{\overline{\Om_0}}\psi_0$ and $M=\max_{\overline{\Om_0}}\psi_0$.  Assume
that $\nabla\psi_0\neq0$ on $\partial\Om_0$ and that every interior critical
point of $\psi_0$ is nondegenerate.  The boundary regularity gives a collar
without critical points; the interior critical points are isolated and hence,
by compactness, finite in number.  Therefore the excluded values, consisting
of $c_0$ and the interior critical values, form a finite set.  For
$t\in(m,M)\setminus\{c_0\}$ with $\nabla\psi_0\neq0$ on $\{\psi_0=t\}$, call
$t$ a \emph{regular} value.  For such a $t$ define the \emph{compact side} of
the level,
\begin{equation}\label{e:side}
\mathcal S_0(t)=
\begin{cases}
\{x\in\Om_0:\psi_0(x)<t\}, & t<c_0,\\[2pt]
\{x\in\Om_0:\psi_0(x)>t\}, & t>c_0 ,
\end{cases}
\end{equation}
which is compactly contained in $\Om_0$ (near $\partial\Om_0$ one has
$\psi_0=c_0$, which is $>t$ when $t<c_0$ and $<t$ when $t>c_0$), and whose boundary is
exactly the level set.  Since $t$ is regular and $\Om_0$ is simply connected,
$\partial\mathcal S_0(t)=\{\psi_0=t\}$ is a finite disjoint union of Jordan
curves; write
\[
\mathcal C_0(t)=\{C_1(t),\dots,C_{N(t)}(t)\}
\]
for them, and for each $C\in\mathcal C_0(t)$ let $D(C)$ be the open region it
encloses (a topological disk contained in $\Om_0$, by the Jordan curve
theorem).

\begin{lemma}[The circle data are invariant]\label{l:circinv}
Let $\eta\in\widehat{\mathcal G}_2$, $\psi=\psi_0\circ\eta^{-1}$, and let
$\mathcal S(t)$ be defined by \eqref{e:side} with $\psi$ in place of
$\psi_0$.  Then $\eta$ maps $\mathcal S_0(t)$ onto $\mathcal S(t)$ and
$\mathcal C_0(t)$ bijectively onto $\mathcal C(t)$, with
$\eta\bigl(D(C)\bigr)=D\bigl(\eta(C)\bigr)$; in particular the number
$N(t)$, the mutual inclusion relations among the enclosed regions, and the
areas $|D(C)|$ are the same for $\psi_0$ and for $\psi$.  One also has
\begin{equation}\label{e:ellcircles}
\ell(t):=\mathcal H^1(\{\psi=t\})=\sum_{C\in\mathcal C(t)}\mathcal H^1(C) .
\end{equation}
\end{lemma}

\begin{proof}
Since $\psi=\psi_0\circ\eta^{-1}$ we have $\{\psi\lessgtr t\}
=\eta(\{\psi_0\lessgtr t\})$, so $\eta(\mathcal S_0(t))=\mathcal S(t)$; and $\eta$ maps $\partial\Om_0$ onto
$\partial\Om$, on which $\psi=c_0$, so $\mathcal S(t)$ is again the compact
side.  Being a homeomorphism of $\overline\Om_0$ onto
$\overline\Om$, $\eta$ maps the boundary curves of $\mathcal S_0(t)$ onto those
of $\mathcal S(t)$, and each enclosed disk $D(C)\subseteq\Om_0$ onto a
topological disk bounded by $\eta(C)$, i.e.\ onto $D(\eta(C))$; areas are
preserved because $\det D\eta=1$ a.e.  Finally $\mathcal S(t)\Subset\Om$
gives $\partial\mathcal S(t)=\{\psi=t\}$, so
$\ell(t)=\mathcal H^1(\partial\mathcal S(t))=\sum_C\mathcal H^1(C)$.
\end{proof}

\begin{theorem}[The general lower bound]\label{t:gen}
Assume the standing assumptions of Section~\ref{ss:rearr}, that
$\nabla\psi_0\neq0$ on $\partial\Om_0$, and that every interior critical
point of $\psi_0$ is nondegenerate.  Then, for every
$\eta\in\widehat{\mathcal G}_2$,
\begin{equation}\label{e:Egen}
\int_{\Om_0}|D\eta\,\Bv_0|^2\,\dd x\ \ge\
\int_m^M\frac{4\pi\,\Lambda(t)^2}{A_0'(t)}\,\dd t,\qquad
\Lambda(t):=\sum_{C\in\mathcal C_0(t)}\sqrt{|D(C)|} ,
\end{equation}
and consequently $E(\eta)\ge\tfrac12\int_m^M\frac{4\pi\Lambda(t)^2}{A_0'(t)}
\,\dd t+2\kap\sqrt{\pi V}$.
\end{theorem}

The quantity $\Lambda(t)$ is determined by the reference data: by
Lemma~\ref{l:circinv} the circle data of $\psi$ coincide with those of
$\psi_0$.

\begin{proof}
Fix a regular value $t\neq c_0$ and let $C\in\mathcal C(t)$.  By the
isoperimetric inequality applied to the Jordan curve $C$, which encloses the
region $D(C)$ of area $|D(C)|$,
\[
\mathcal H^1(C)\ \ge\ 2\sqrt{\pi\,|D(C)|} .
\]
Summing over $C\in\mathcal C(t)$ and using \eqref{e:ellcircles},
\[
\ell(t)\ \ge\ 2\sqrt\pi\sum_{C\in\mathcal C(t)}\sqrt{|D(C)|}
\ =\ 2\sqrt\pi\,\Lambda(t) ,
\]
the last step by Lemma~\ref{l:circinv}, which identifies the circle data of
$\psi$ with those of $\psi_0$.  Substituting this into \eqref{e:csleaves} and
integrating over $t$ as in the proof of Lemma~\ref{l:ps} gives
\eqref{e:Egen}.  The perimeter term is again bounded by the isoperimetric inequality
$\Per(\eta(\Om_0))\ge2\sqrt{\pi V}$.
\end{proof}

\section{The three-dimensional problem}\label{s:3d}

Section~\ref{s:2d} solved the two-dimensional problem by reducing it to a
scalar one: on a simply connected plane domain a divergence-free field has a
stream function, area-preserving maps transport it by composition, and the
magnetic energy is the Dirichlet energy of that stream function.  In three
dimensions the magnetic energy does not control the full gradient of
admissible maps (Theorem~\ref{t:cav3d}); a scalar invariant survives only in
the axisymmetric reduction, where the poloidal flux function plays the role
of the stream function, and a conditional lower bound is proved
(Proposition~\ref{p:axisbound}); that reduction is treated in Section~\ref{ss:ax}
below.  The bound is conditional in two respects --- uniform radial confinement
$\Om\subset\{r\le r_1\}$ and compact meridian levels --- and its equality
case is defective (Remark~\ref{r:axisfail}); it is to be read as a necessary
condition rather than as a sharp variational characterization.

\subsection{Obstructions and conditional bounds}

\begin{theorem}[Energy does not control the full gradient in three dimensions]\label{t:cav3d}
Let $\Om_0\subset\R^3$ be bounded with smooth boundary and let
$\Bv_0\in C^\infty(\overline\Om_0;\R^3)$ satisfy $\dive\Bv_0=0$,
$\Bv_0\cdot\nv_0=0$ on $\partial\Om_0$, and $\Bv_0\not\equiv0$.  Then there
are $\eta_k\in\mathcal G$, $k\ge1$, with
\[
\sup_k\int_{\Om_0}|D\eta_k\,\Bv_0|^2\,\dd x<\infty,\qquad
\int_{\Om_0}|D\eta_k|^2\,\dd x\ \longrightarrow\ \infty .
\]
Consequently an energy bound alone does not give a uniform $L^2$ gradient
bound for arbitrary admissible maps in three dimensions.
\end{theorem}

The proof is a shear supported on a flow box and is deferred to
Appendix~\ref{app:3d}; it is the three-dimensional counterpart of
Theorem~\ref{t:cavgen}, with a Hamiltonian (rather than a leaf-wise)
shear on the transverse disk.  The mechanism is elementary and physical: on a
flow box of $\Bv_0$ the shear advances the field lines by an amount proportional
to the transverse coordinate, keeping the map divergence-free, tangential to the
boundary, and $L^\infty$-bounded; the field is left almost invariant while the
embedding acquires an $O(Ak)$ gradient, which is the obstruction.

\subsection{The axisymmetric reduction}\label{ss:ax}

Axial symmetry is the one three-dimensional reduction that retains a scalar
flux function, and with it part of the structure of Section~\ref{s:2d}.  Here
$\Psi$ denotes the poloidal flux
function; we reserve $\psi$ for the planar stream function of
Section~\ref{s:2d}.

\emph{The class and the energy.}  Let $(r,\varphi,z)$ be cylindrical
coordinates with frame $(e_r,e_\varphi,e_z)$, and let $\Om$ and $\Bv$ be
invariant under rotations about the $z$-axis.  An axisymmetric
divergence-free field has the Grad--Shafranov form
\begin{equation}\label{e:gs}
\Bv=\frac1r\bigl(\nabla\Psi\times e_\varphi\bigr)+B_\varphi\,e_\varphi,
\qquad \nabla\Psi=(\partial_r\Psi)\,e_r+(\partial_z\Psi)\,e_z ,
\end{equation}
with $\Psi=\Psi(r,z)$ the poloidal flux function and $B_\varphi=B_\varphi(r,z)$;
conversely every such pair defines a divergence-free axisymmetric field.  Since
$\Bv\cdot\nabla\Psi=0$, the level curves of $\Psi$ generate the flux surfaces,
and $\Bv\cdot\nv=0$ on $\partial\Om$ holds if and only if each component of
$\partial D$ is a level curve of $\Psi$, where
$D=\{(r,z):(r,\varphi,z)\in\Om\}$ is the meridian domain; we take $\Om$ to be
diffeomorphic to a ball or to a solid torus.  In both cases $D$ is simply
connected: a ball-like body is generated by a half-disk meeting the axis, and
a solid torus by an off-axis disk.  (An annular $D$ would generate a toroidal
shell with two torus boundary components; there $\Psi$ is globally defined
only if the period of the closed $1$-form $rB_z\,\dd r-rB_r\,\dd z$ around
the hole vanishes, that is, if the flux through the inner boundary component
is zero, which tangency guarantees.)  With
$|\Bv|^2=|\nabla\Psi|^2/r^2+B_\varphi^2$ and $\dd y=r\,\dd r\,\dd z\,\dd\varphi$,
\begin{equation}\label{e:gsE}
E=\pi\int_D\Bigl(\frac{|\nabla\Psi|^2}{r}+B_\varphi^2\,r\Bigr)\dd r\,\dd z
+2\pi\kap\int_{\partial D}r\,\dd s ,
\end{equation}
the last term being $\kap\Area(\partial\Om)$ for a surface of revolution, with
$s$ the arclength in the meridian plane.

\emph{The frozen-in law.}  Let $\mathcal G_{\rm ax}\subset\mathcal G$ denote the
subset of axisymmetric volume-preserving embeddings; in cylindrical
coordinates these read $(r,z)\mapsto(R,Z)$ with
$R\,\dd R\,\dd Z=r\,\dd r\,\dd z$, composed with a toroidal shift
$\varphi\mapsto\varphi+\alpha(r,z)$.  They preserve the weighted measure
$r\,\dd r\,\dd z$, and the poloidal flux is transported as a scalar: if
$\eta\in\mathcal G_{\rm ax}$ and $\Psi$ is the flux function of $\Bv_\eta$,
then $\Psi\circ\eta=\Psi_0+C$, and we normalize $C=0$.  To see this, consider
the flux $2$-form
\[
\beta=\iota_{\Bv}\bigl(r\,\dd r\wedge\dd\varphi\wedge\dd z\bigr)
=\dd\Psi\wedge\dd\varphi-B_\varphi\,\dd r\wedge\dd z ,
\]
the identity following from \eqref{e:gs}: with
$\Bv=\Bv_r\partial_r+(B_\varphi/r)\partial_\varphi+\Bv_z\partial_z$,
\[
\beta=r\Bv_r\,\dd\varphi\wedge\dd z-B_\varphi\,\dd r\wedge\dd z
+r\Bv_z\,\dd r\wedge\dd\varphi ,
\]
and $r\Bv_r=-\partial_z\Psi$, $r\Bv_z=\partial_r\Psi$.  The
$\dd r\wedge\dd z$ term carries no $\dd\varphi$, and $B_\varphi$ is not
transported as a scalar.  The frozen-in law
$\Bv_\eta=\eta_*\Bv_0$ and $\det D\eta=1$ give $\eta^*\beta=\beta_0$.  Write
an axisymmetric embedding as
$\eta(r,\varphi,z)=(R(r,z),\varphi+\alpha(r,z),Z(r,z))$.  The terms of
$\eta^*\beta$ containing $\dd\varphi$ are then
$\dd(\Psi\circ\eta)\wedge\dd\varphi$, while those of $\beta_0$ are
$\dd\Psi_0\wedge\dd\varphi$; hence
$\dd(\Psi\circ\eta-\Psi_0)=0$ on the connected meridian domain and
$\Psi\circ\eta=\Psi_0+C$. 

Since $\eta$ preserves $r\,\dd r\,\dd z$, the
distribution function
\begin{equation}
A(t)=\int_{\{\Psi<t\}}r\,\dd r\,\dd z
\end{equation}
is invariant, exactly as in Lemma~\ref{l:conj}(ii).

\begin{lemma}[Weighted isoperimetric inequality]\label{l:axisiso}
Let $\Om\subset\R^3$ be a bounded body of revolution generated by a meridian
domain $D\subset\{r\ge0\}$ with piecewise $C^1$ boundary, the segment of
$\partial D$ on the axis having zero weight in
$2\pi\int_{\partial D}r\,\dd s$.  Then
\begin{equation}\label{e:axisiso}
2\pi\int_{\partial D}r\,\dd s\ \ge\ (36\pi)^{1/3}
\Bigl(2\pi\int_D r\,\dd r\,\dd z\Bigr)^{2/3},
\end{equation}
with equality if and only if $\Om$ is a ball.
\end{lemma}

\begin{proof}
For a surface and a body of revolution the area and the volume are
$2\pi\int_{\partial D}r\,\dd s$ and $2\pi\int_Dr\,\dd r\,\dd z$, the axis
segment contributing zero to the area because $r=0$ there.  The assertion
is therefore the isoperimetric inequality in $\R^3$ applied to the body of
revolution; it holds for sets of finite perimeter, hence for piecewise $C^1$
meridian domains, and its equality case is the ball.
\end{proof}

Lemma~\ref{l:axisiso} settles the shape term: surface tension alone selects a
ball, with the same constant as in $\R^3$.  For the magnetic term, the planar
coarea argument gives the conditional lower bound below.  Its equality
conditions have a geometric defect for nondegenerate closed meridian levels,
as explained in Remark~\ref{r:axisfail}.

\begin{proposition}[Conditional lower bound]\label{p:axisbound}
Let $\eta\in\mathcal G_{\rm ax}$, $\Om=\eta(\Om_0)$ and let $\Psi$ be the
flux function of $\Bv_\eta$, and assume
\begin{enumerate}[\rm (a)]
\item $\Om\subset\{r\le r_1\}$;
\item \emph{Compact meridian levels.} $\Psi$ is constant on each component
of $\partial\Om$, and for a.e.\ $t\in(m,M)$ the set $\{\Psi<t\}$ is
compactly contained in the meridian domain with
$\partial\{\Psi<t\}=\{\Psi=t\}$ a finite union of closed curves,
\end{enumerate}
where $m=\min_{\overline\Om}\Psi$ and $M=\max_{\overline\Om}\Psi$.  Then
\begin{equation}\label{e:axisbound}
E(\eta)\ \ge\ \frac{4\pi^2}{r_1}\int_m^M\frac{A(t)}{A'(t)}\,\dd t
\ +\ \kap\,(36\pi)^{1/3}V^{2/3}.
\end{equation}
\end{proposition}

Hypothesis {\rm (b)}, together with the continuity of $\Psi$, forces $\Psi=M$ on
every component of $\partial\Om$; otherwise, for $t$ close to $M$, the closure of
$\{\Psi<t\}$ would touch $\partial\Om$ and $\{\Psi<t\}$ would not be compactly
contained in the meridian domain.  Thus \eqref{e:axisbound}
is a conditional lower bound, valid under the radial confinement {\rm (a)} and the
compact-regular-level hypothesis {\rm (b)}.

\begin{proof}
The toroidal term in \eqref{e:gsE} is nonnegative and is dropped.  For a.e.\
$t$ put $\ell(t)=\mathcal H^1(\{\Psi=t\})$.  The coarea formula with weight $r$
gives $A'(t)=\int_{\{\Psi=t\}}r|\nabla\Psi|^{-1}\,\dd\mathcal H^1$ and
\[
\int_D\frac{|\nabla\Psi|^2}{r}\,\dd r\,\dd z
=\int_m^M\Bigl(\int_{\{\Psi=t\}}\frac{|\nabla\Psi|}{r}\,\dd\mathcal H^1\Bigr)\dd t ,
\]
and the Cauchy--Schwarz inequality applied on $\{\Psi=t\}$ to
$\bigl(|\nabla\Psi|/r\bigr)^{1/2}$ and $\bigl(r/|\nabla\Psi|\bigr)^{1/2}$ gives
\[
\int_{\{\Psi=t\}}\frac{|\nabla\Psi|}{r}\,\dd\mathcal H^1
\ \ge\ \frac{\ell(t)^2}{A'(t)} .
\]
Next, the components $E_i(t)$ of $\{\Psi<t\}$ obey
$\sum_i\Per(E_i(t))\le \ell(t)$, and the planar isoperimetric inequality gives
$|E_i(t)|\le\Per(E_i(t))^2/(4\pi)$; since $r\le r_1$ on $\Om$ by (a),
\[
A(t)=\sum_i\int_{E_i(t)}r\,\dd r\,\dd z\ \le\ r_1\sum_i|E_i(t)|
\ \le\ \frac{r_1}{4\pi}\sum_i\Per(E_i(t))^2\ \le\ \frac{r_1}{4\pi}\,\ell(t)^2 ,
\]
that is, $\ell(t)^2\ge4\pi A(t)/r_1$.  Hence the magnetic part is at least the
first term of \eqref{e:axisbound}, and the surface term is estimated by
Lemma~\ref{l:axisiso}.
\end{proof}

\begin{remark}[Defect of \eqref{e:axisbound}]%
\label{r:axisfail}
(i) \emph{Dual weights.}  The energy density in \eqref{e:gsE} is
$|\nabla\Psi|^2/r$, whereas the constraint is the invariance of the
distribution with respect to $r\,\dd r\,\dd z$.  Equality in the
Cauchy--Schwarz step above requires $|\nabla\Psi|=c\,r$ on each level curve,
whereas for a nested family of circles about a magnetic axis $|\nabla\Psi|$
is constant on each circle while $r$ is not.  In addition, equality in
$A(t)\le r_1\sum_i|E_i(t)|$ would require $r=r_1$ almost everywhere on a
positive-area sublevel set, which is impossible.  Thus the bound is not
attained by a nondegenerate family of positive-area levels.  The
P\'olya--Szeg\H{o} argument of Lemma~\ref{l:ps} has no such defect, because
there the level-set measure and the energy density are the same.

(ii) \emph{Different optimal shapes.}  The equality case of \eqref{e:axisiso}
is a half-disk, that is, a ball, whereas a single-axis field has closed level
curves in the meridian plane, that is, tori.  The magnetic and surface terms
therefore have different equality structures.  The planar coincidence (that
surface tension and the magnetic energy are both minimized by the disk)
therefore has no axisymmetric counterpart.

(iii) \emph{Strength and non-emptiness.}  The factor $4\pi^2/r_1$ degenerates
as $r_1\to\infty$, so \eqref{e:axisbound} gives no uniform information without
radial confinement; by (i) and (ii) it is never sharp on a nondegenerate
family of closed levels, and at present it is a necessary condition, not a
variational characterization.  The hypotheses are nevertheless non-vacuous.
Let $\Om$ be the solid torus generated by rotating
$D=\{(r,z):(r-r_0)^2+z^2<\epsilon^2\}$ about the $z$-axis, with
$0<\epsilon<r_0$, and take $B_\varphi=0$ and
$\Psi=\frac12[(r-r_0)^2+z^2]$.  Then $m=0$, $M=\epsilon^2/2$, each level
$\{\Psi=t\}$ with $0<t<M$ is one closed curve compactly contained in $D$,
and $r_1=r_0+\epsilon$ is admissible in (a).  Writing
$s^2=(r-r_0)^2+z^2$ one has $\int_{\{\Psi<t\}}r\,\dd r\,\dd z=2\pi r_0t$,
hence $A(t)=2\pi r_0t$, $A'(t)=2\pi r_0$, and the magnetic part of
\eqref{e:axisbound} is
\[
\frac{4\pi^2}{r_0+\epsilon}\int_0^{\epsilon^2/2}t\,\dd t
=\frac{\pi^2\epsilon^4}{2(r_0+\epsilon)}>0 .
\]
The actual magnetic energy is
$\pi\int_D|\nabla\Psi|^2r^{-1}\,\dd r\,\dd z
=\pi\int_D s^2/r\,\dd r\,\dd z$; expanding $1/r$ about $r_0$ (the odd term
integrates to zero) gives $\frac{\pi^2\epsilon^4}{2r_0}(1+O((\epsilon/r_0)^2))$,
so the lower bound captures the leading order as $\epsilon/r_0\to0$ even
though equality is unattainable.  Thus \eqref{e:axisbound} is a usable
quantitative necessary condition for tightly confined toroidal
configurations, but not a sharp inequality in the sense of
Theorem~\ref{t:rearr}.
\end{remark}

\appendix

\section{The three-dimensional shear}\label{app:3d}

Throughout, $(t,u,v)$ are the flow-box coordinates of Step~1 below, with
$\Bv_0=\partial_t$ and volume element $\mathcal J(u,v)\,\dd t\,\dd u\,\dd v$;
$X_k$ is the Hamiltonian field of $h_k=\frac Ak\sin(ku)\chi$ with respect to
$\omega=\mathcal J\,\dd u\,\dd v$, $\Theta_\tau$ its flow, and $M_\tau$, $a_\tau,
b_\tau,c_\tau,d_\tau$ the entries of $DX_k$ and $D\Theta_\tau$ in the
$(u,v)$-block.

\begin{proof}[Proof of Theorem~\ref{t:cav3d}]
\emph{1. Flow-box coordinates.}  Choose $x_0$ with $\Bv_0(x_0)\neq0$ and a
small two-dimensional disk $\Sigma$ through $x_0$, transverse to
$\Bv_0(x_0)$, with coordinates $(u,v)$ and parametrization
$y(u,v)$.  Let $\Phi_t$ be the flow of $\Bv_0$.  For small $\delta>0$ the map
\[
(t,u,v)\ \longmapsto\ \Phi_t\bigl(y(u,v)\bigr),\qquad
(t,u,v)\in U:=(-\delta,\delta)\times\Sigma ,
\]
is a diffeomorphism onto a neighborhood of $x_0$, and in these coordinates
$\Bv_0=\partial_t$.  Because $\dive\Bv_0=0$ the flow preserves the Euclidean
volume, so the pull-back of $\dd x$ is $\mathcal J(u,v)\,\dd t\,\dd u\,\dd v$
with $\mathcal J>0$ independent of $t$.  Write
$\omega=\mathcal J(u,v)\,\dd u\,\dd v$ for the
corresponding area form on $\Sigma$.

\emph{2. A rapidly oscillating area-preserving map of $\Sigma$.}  Choose
rectangles $I\times K\Subset\Sigma$ and a nonnegative
$\chi\in C_c^\infty(\Sigma)$ with $\chi=1$ on $I\times K$.  Fix a small
amplitude $A>0$.  For $k\ge1$ put
\[
h_k(u,v)=\frac{A}{k}\sin(ku)\,\chi(u,v)\ \in C_c^\infty(\Sigma) ,
\]
and let $X_k$ be the Hamiltonian vector field of $h_k$ with respect to
$\omega$, i.e.\ $\iota_{X_k}\omega=\dd h_k$.  Then
$\mathcal L_{X_k}\omega=\dd\iota_{X_k}\omega+\iota_{X_k}\dd\omega
=\dd^2h_k=0$, so $\dive_\omega X_k=0$ and the flow
$\Theta_\tau$ of $X_k$ preserves $\omega$ for every $\tau$; and $\Theta_\tau$ is
defined for all $\tau$ and preserves $\Sigma$, because $X_k$ is compactly
supported in $\Sigma$.  Since $|\nabla h_k|\le C\,A$ uniformly in $k$ we have
\begin{equation}\label{e:Xkbd}
\sup_\Sigma|X_k|\ \le\ C_0A ,
\end{equation}
while the entry $\partial_uX_k^v=\frac{A}{\mathcal J}k\sin(ku)+O(A)$ on
$I\times K$ satisfies
\begin{equation}
|DX_k|\ \ge\ c_0Ak\qquad\text{on }
G_k:=\{(u,v)\in I\times K:|\sin(ku)|\ge\tfrac34\}
\end{equation}
for all sufficiently large $k$.  The measure of $G_k$ is bounded below by a
positive constant independent of $k$; the bound is asserted only on $G_k$,
not on all of $I\times K$.

\emph{3. The map.}  Let $\theta\in C_c^\infty(-\delta,\delta)$ satisfy
$\theta=1$ on $(-\tfrac\delta2,\tfrac\delta2)$.  Define on $U$
\[
\tilde\eta_k(t,u,v)=\bigl(t,\ \Theta_{\theta(t)}(u,v)\bigr) ,
\]
and $\tilde\eta_k=\mathrm{id}$ on $\overline\Om_0\setminus U$.  Because
$\theta=0$ near $t=\pm\delta$ and $\Theta_0=\mathrm{id}$, and because
$X_k$ is supported inside $\Sigma$ so that $\Theta_\tau=\mathrm{id}$ near
$\partial\Sigma$, the map $\tilde\eta_k$ agrees with the identity near
$\partial U$; it is therefore a smooth diffeomorphism of
$\overline\Om_0$ onto itself.

\emph{4. It is volume-preserving.}  In the coordinates $(t,u,v)$ the
Jacobian matrix of $\tilde\eta_k$ is block triangular, with $(t,t)$ entry
$1$ and $(u,v)$-block $D_{(u,v)}\Theta_{\theta(t)}$.  That block does not
have coordinate determinant $1$: since $\Theta_\tau$ preserves
$\omega=\mathcal J\,\dd u\wedge\dd v$,
\[
\mathcal J(\Theta_\tau)\,\det D\Theta_\tau=\mathcal J,\qquad
\det D_{(u,v)}\Theta_\tau=\frac{\mathcal J(u,v)}
{\mathcal J(\Theta_\tau(u,v))}.
\]
What is preserved is the $3$-form $\dd t\wedge\omega$: $\tilde\eta_k$ fixes
$t$, so
\[
\tilde\eta_k^*(\dd t\wedge\omega)=\dd t\wedge\Theta_{\theta(t)}^*\omega
=\dd t\wedge\omega ,
\]
the $\theta'(t)\,\dd t$ terms vanishing when wedged with $\dd t$.  Because
$\dd t\wedge\omega=\mathcal J\,\dd t\,\dd u\,\dd v$ is the pull-back of the
Euclidean volume under the flow-box chart, $\det D\eta_k=1$ in Cartesian
coordinates; thus
$\tilde\eta_k\in\mathcal G$; write $\eta_k:=\tilde\eta_k$ for this map.

\emph{5. The magnetic energy stays bounded.}  From the definition,
\[
D\tilde\eta_k(\partial_t)
=\partial_t+\theta'(t)\,X_k\bigl(\Theta_{\theta(t)}(u,v)\bigr) .
\]
By \eqref{e:Xkbd} the second term is bounded in modulus by
$C_0A\,\|\theta'\|_\infty$, uniformly in $k$, and vanishes off $U$; since
$\Bv_0=\partial_t$ on $U$, it follows that $|D\tilde\eta_k\,\Bv_0|\le C$
there, while $|D\tilde\eta_k\,\Bv_0|=|\Bv_0|$ off $U$.  Hence
$\sup_k\int_{\Om_0}|D\eta_k\Bv_0|^2\,\dd x<\infty$.  Since $\eta_k$ is a
diffeomorphism of $\overline\Om_0$ onto itself, $\Om_{\eta_k}=\Om_0$ and
$\Area(\partial\Om_{\eta_k})$ is independent of $k$, so
$\sup_k E(\eta_k)<\infty$ as well: an energy bound alone does not control
$\int_{\Om_0}|D\eta_k|^2\,\dd x$.

\emph{6. The gradient blows up.}  On
$(-\tfrac\delta2,\tfrac\delta2)\times\Sigma$ we have $\theta=1$, so
$D\tilde\eta_k$ contains the block $D\Theta_1$.  In the coordinates
$(u,v)$ the components of $X_k$ are
\[
X_k^u=\tfrac1{\mathcal J}\partial_vh_k=O(A/k),\qquad
X_k^v=-\tfrac1{\mathcal J}\partial_uh_k
=-\tfrac{A}{\mathcal J}\cos(ku)\,\chi+O(A/k) ,
\]
so the flow displaces points by $O(A)$ in $v$ and by $O(A/k)$ in $u$.
Choose $I'\Subset I$ and $K'\Subset K$, then choose $A$ small enough that
trajectories starting in $I'\times K'$ stay in $I\times K$ and their phase
change $k|\Delta u|$ is small.  Start on
$G'_k=\{(u,v)\in I'\times K':\sin(ku)\ge7/8\}$; this set has measure
bounded below independently of large $k$.  Along each such trajectory,
$\partial_uX_k^v=\tfrac{A}{\mathcal J}k\sin(ku)+O(A)$ is positive and
has magnitude at least $cAk$, since $\chi=1$ on $I\times K$.  Restricting
to the positive branch of $\sin(ku)$ is essential: on the negative branch
the same entry is $\le-cAk$.

The derivative of the flow obeys
$\frac{\dd}{\dd\tau}D\Theta_\tau=DX_k(\Theta_\tau)\,D\Theta_\tau$ with
$D\Theta_0=I$.  Write $M_\tau:=DX_k(\Theta_\tau)$ and
$D\Theta_\tau=\begin{pmatrix}a_\tau&b_\tau\\c_\tau&d_\tau\end{pmatrix}$ in the
$(u,v)$ block, so that $a_0=d_0=1$ and $b_0=c_0=0$.  Differentiating the
formula $X^u_k=\frac{1}{\mathcal J}\partial_vh_k$,
$X^v_k=-\frac1{\mathcal J}\partial_uh_k$ with
$h_k=\frac Ak\sin(ku)\chi$ gives the entry sizes
\[
M^v_u=\partial_uX^v_k=\tfrac{A}{\mathcal J}k\sin(ku)+O(A),\qquad
M^u_u=O(A),\quad M^v_v=O(A),\quad M^u_v=\partial_vX^u_k=O(A/k).
\]
Thus the only large entry is $M^v_u$, and it sits in a strictly lower
triangular position, so the growth it drives is only linear in $Ak$ rather
than exponential.  Quantitatively, the columns obey
$\dot a=M^u_ua+M^u_vc$, $\dot c=M^v_ua+M^v_vc$ and the identical equations for
$b,d$; hence, writing $Y:=|a|+|b|$ and $Z:=|c|+|d|$,
\[
|\dot Y|\le C\Bigl(AY+\tfrac AkZ\Bigr),\qquad
|\dot Z|\le C\Bigl(Ak\,Y+AZ\Bigr),
\]
and Gronwall's inequality applied to $W:=Y+Z/(Ak)$ yields $W(\tau)\le W(0)e^{C'\tau}$, so
$Y=O(1)$ and $Z=O(Ak)$ uniformly for $\tau\in[0,1]$ and all $k$.  In
particular $|c_\tau|=O(Ak)$ and $|a_\tau|=O(1)$, so the equation for $a$
gives $\dot a=M^u_ua+M^u_vc=O(A)\cdot O(1)+O(A/k)\cdot O(Ak)=O(A)$ with
$a_0=1$, i.e.\ $(D\Theta_\tau)^u_u=a_\tau=1+O(A)$; also
$\int_\tau^1M^v_v\,\dd s=O(A)$, so $e^{\int_\tau^1M^v_v\,\dd s}=1+O(A)$.

The $(v,u)$ entry of the
variational equation is
\[
\frac{\dd}{\dd\tau}(D\Theta_\tau)^v_u
=M^v_u\,(D\Theta_\tau)^u_u+M^v_v\,(D\Theta_\tau)^v_u ,
\]
whose solution is
\[
(D\Theta_1)^v_u
=\int_0^1 e^{\int_\tau^1 M^v_v(\Theta_s)\,\dd s}\,
M^v_u(\Theta_\tau)\,(D\Theta_\tau)^u_u\,\dd\tau .
\]
Along a trajectory in $G'_k$ the exponent satisfies
$e^{\int_\tau^1 M^v_v\,\dd s}=1+O(A)$, the factor $(D\Theta_\tau)^u_u=1+O(A)$,
and $M^v_u(\Theta_\tau)$ is positive with $\int_0^1M^v_u\ge cAk$ (the
phase change being small).  Hence $(D\Theta_1)^v_u\ge c'Ak$ for a constant
$c'>0$ independent of $k$.  Since the coordinate frame is comparable to an
orthonormal one on the compact set $U$, this gives $|D\Theta_1|\ge c''Ak$, and
hence $|D\eta_k|\ge c''Ak$ on the subset
$(-\tfrac\delta2,\tfrac\delta2)\times G'_k$ of $\Om_0$, whose measure is
bounded below independently of $k$.  Therefore
\[
\int_{\Om_0}|D\eta_k|^2\,\dd x\ \ge\ c'\,A^2k^2\ \to\ \infty ,
\]
which completes the proof.
\end{proof}

\section{The elliptic example}\label{app:ellipse}

We verify Theorem~\ref{t:rearr} on the datum of Example~\ref{ex:ellipse} in
full.  Throughout $a>b>0$,
\[
\Om_0=\{x_1^2/a^2+x_2^2/b^2<1\},\qquad
\psi_0=x_1^2/a^2+x_2^2/b^2-1,\qquad \Bv_0=\nabla^\perp\psi_0 .
\]
Then $\psi_0=0$ on $\partial\Om_0$, the only critical point is the
nondegenerate minimum at the origin with value $-1$, and no other critical
level exists, so {\rm (H1)} holds with $m=-1$ and $M=0$.

\emph{(1) Area function.}  For $t\in[-1,0]$ the set $\{\psi_0<t\}$ is the
ellipse with semiaxes $a\sqrt{1+t}$ and $b\sqrt{1+t}$, of area
$\pi ab(1+t)$.  Hence
\[
A_0(t)=\pi ab(1+t),\qquad A_0'(t)=\pi ab,\qquad V=\pi ab .
\]

\emph{(2) The rearranged bound.}
\[
M_*=\tfrac12\int_{-1}^{0}\frac{4\pi A_0(t)}{A_0'(t)}\,\dd t
=\tfrac12\int_{-1}^{0}4\pi(1+t)\,\dd t=2\pi\cdot\tfrac12=\pi .
\]

\emph{(3) Magnetic energy of the identity.}  Here
$\Bv_0=(-2x_2/b^2,\,2x_1/a^2)$, so
$|\Bv_0|^2=4x_2^2/b^4+4x_1^2/a^4$.  Using
$\int_{\Om_0}x_1^2\,\dd x=\pi ab\cdot a^2/4$ and
$\int_{\Om_0}x_2^2\,\dd x=\pi ab\cdot b^2/4$,
\[
\int_{\Om_0}|\Bv_0|^2\,\dd x
=\frac{4}{b^4}\cdot\frac{\pi ab^3}{4}+\frac{4}{a^4}\cdot\frac{\pi a^3b}{4}
=\pi ab\Bigl(\frac1{a^2}+\frac1{b^2}\Bigr),
\]
that is,
\[
\tfrac12\int_{\Om_0}|\Bv_0|^2\,\dd x=\frac{\pi(a^2+b^2)}{2ab}
=\pi+\frac{\pi(a-b)^2}{2ab}\ \ge\ \pi ,
\]
with equality if and only if $a=b$.

\emph{(4) The rearrangement and its image.}  By \eqref{e:psistar},
$\varrho(t)=\sqrt{A_0(t)/\pi}=\sqrt{ab(1+t)}$, hence
$\psi_*(r)=r^2/(ab)-1$ and $\Bv_*=\psi_*'(r)\,e_\theta
=\frac{2r}{ab}e_\theta$ on the disk of radius $\sqrt{ab}$; moreover
$\Delta\psi_*=4/(ab)$ is constant, in agreement with
Proposition~\ref{p:scalar}: here all local primitives coincide with the
global constant $f\equiv4/(ab)$.  The image is a rigid
rotation, although $\Bv_0$ is not.

\emph{(5) The strict inequality.}  Combining (2)--(4),
\[
E(\mathrm{id})-E(\eta_*)=\frac{\pi(a-b)^2}{2ab}
+\kap\bigl(\Per(\Om_0)-2\pi\sqrt{ab}\bigr)>0
\qquad(a\ne b),
\]
by the isoperimetric inequality $\Per(\Om_0)>2\pi\sqrt{ab}=2\sqrt{\pi V}$
for $a\ne b$.  The first term is of order $(a-b)^2$ for $a$ close to $b$: the
strict inequality of Example~\ref{ex:ellipse} is sharp to second order, since
at $a=b$ the datum is circular and the identity becomes a minimizer.

\section*{Acknowledgment}

Hao's research was partially supported by the National Natural Science Foundation of China (NSFC) under Grant No. 12171460, the CAS Project for Young Scientists in Basic Research under Grant No. YSBR-031 and the National Key R\&D Program of China under Grant No. 2021YFA1000800.

\medskip
\noindent\textbf{AI declaration.}
The author employed generative AI tools to assist with language polishing and manuscript organization. These tools were also utilized to generate and refine the TikZ code for schematic figures, audit mathematical arguments, and aid in revising portions of the proofs. The author retain full responsibility for the accuracy and integrity of all mathematical content and figures.

\end{document}